%% file: SISurfaceTension.tex
\RequirePackage{fix-cm}
\documentclass[smallextended]{svjour3}       
\smartqed  
\usepackage{url}
\usepackage{xcolor}
\usepackage{rotating}
\definecolor{newcolor}{rgb}{.8,.349,.1}

\usepackage{amsmath}	
\usepackage{subfigure}

\usepackage{latexsym}
\usepackage{graphicx}  

\usepackage{url}
\usepackage{multirow}

\usepackage{listings}
\usepackage{verbatim}
\usepackage{underscore}

\usepackage{amssymb} 
\usepackage{amsfonts} 
\usepackage{mathtools}

\usepackage{tabularx}
\usepackage{booktabs}

   \newcolumntype{C}{>{\centering\let\newline\\\arraybackslash}X}
   \newcolumntype{R}{>{\raggedleft\let\newline\\\arraybackslash}X}
   \newcolumntype{L}{>{\raggedright\let\newline\\\arraybackslash}X}
   \newcolumntype{Z}{>{\raggedright\let\newline\\\arraybackslash}X}

\usepackage{cite}

\newcommand{\aposteriori}{\textit{a posteriori}}

\newcommand{\fls}{\vphantom{l^s}}

\newcommand{\cshear}{c_\up{s}}
\newcommand{\reynolds}{\mathbb{R}\mathrm{e}}

\newcommand{\weber}{\mathbb{W}\mathrm{e}}

\newcommand{\pnpm}[2]{\ensuremath{\mathbb{P}_{#1}\mathbb{P}_{#2}}}

\newcommand{\vrec}{\vec{V}^\up{r}}
\newcommand{\qrec}{\vec{Q}^\up{r}}

\usepackage{bm}
   \renewcommand{\vec}[1]{\bm{\mathrm{#1}}}
   \newcommand{\uvec}[1]{\bm{\mathrm{\hat{#1}}}}
 
\DeclarePairedDelimiter{\norm}{\lVert}{\rVert}

\newcommand{\up}[1]{\ensuremath{\mathrm{#1}}}

\newcommand{\od}[2]{\dfrac{\de{#1}}{\vphantom{l^l}\de{#2}}}

\DeclareMathOperator{\de}{d\!}
\DeclareMathOperator{\tr}{tr}
\DeclareMathOperator{\dev}{dev}
\DeclareMathOperator{\abs}{abs}
\DeclareMathOperator{\sign}{sign}

\DeclareMathAlphabet{\mathsfbi}{OT1}{\sfdefault}{bx}{sl}
\DeclareMathVersion{sfletters}
\SetSymbolFont{letters}{sfletters}{OML}{ntxsfmi}{b}{it}

\usepackage{scalefnt}
\usepackage{relsize}
\newcommand{\transpose}{{\mathsf{\mathsmaller{T}}}}
\newcommand{\Ge}{\mathbf{G}}

\newcommand{\Lstar}{\mathbf{L}_\ast}
\newcommand{\Pstar}{\vec{P}_\ast}
\newcommand{\sour}{\vec{Z}}
\newcommand{\gradflux}{\vec{K}_v}

\newcommand{\fll}{\vphantom{l^l}} 
\newcommand{\flll}{\vphantom{{l^l}^l}}

\newcommand{\devi}{\widetilde{\Delta\vec{V}}_i}
\newcommand{\devic}{\widetilde{\Delta\vec{V}}_i^3}
\newcommand{\dgrad}{\nabla_h}

\newcommand{\taut}{\tau_\textsc{h}}

\newcommand{\Gh}{\hat{\vec{G}}}
\newcommand{\ege}{\vec{E}\,\Gh\,\vec{E}^{-1}}
\newcommand{\erge}{\vec{E}\,\rGh\,\vec{E}^{-1}}
\newcommand{\egdge}{\vec{E}\,\Gh\,(\dev\Gh)\,\vec{E}^{-1}}
\newcommand{\ergdge}{\vec{E}\,\rGh\,(\dev\Gh)\,\vec{E}^{-1}}
\newcommand{\rGh}{\Gh^{1/2}}
\newcommand{\ddt}{\frac{\de}{\de t}}

\newcommand{\pd}[3][0pt]{\frac{\raisebox{#1}{$\partial#2$}}{\vphantom{l^l}\partial#3}}

\begin{document}

\title{An exactly curl-free staggered semi-implicit finite volume scheme for a first order hyperbolic model of viscous flow with surface tension}

\titlerunning{Curl-free semi-implicit scheme for a hyperbolic model of surface tension}        

\author{Simone Chiocchetti \and Michael Dumbser} 


\institute{
S. Chiocchetti \at
Department of Civil, Environmental and Mechanical Engineering, University of Trento, Via Mesiano 77, 38123 Trento, Italy \\
\email{simone.chiocchetti@unitn.it}         
\and
M. Dumbser \at
Department of Civil, Environmental and Mechanical Engineering, University of Trento, Via Mesiano 77, 38123 Trento, Italy 
\email{michael.dumbser@unitn.it}  
}

\date{Received: date / Accepted: date}

\maketitle

\begin{abstract} 
In this paper, we present a semi-implicit numerical solver for a \textit{first order 
hyperbolic} formulation of \textit{two-phase flow with surface tension and viscosity}. 
The numerical method addresses several complexities presented by the PDE system
in consideration: (i) The presence of involution constraints of curl type in the 
governing equations requires explicit enforcement of the zero-curl property of
certain vector fields (an interface field and a distortion field); the problem is 
eliminated entirely  by adopting a set of \textit{compatible curl and gradient discrete differential operators on a staggered grid}, 
allowing to \textit{preserve the Schwartz identity of cross-derivatives exactly}. (ii) The 
evolution equations feature \textit{highly nonlinear stiff algebraic source terms} which are used for 
the description of viscous interactions as emergent behaviour of an elasto-plastic solid in the stiff strain relaxation limit; 
such source terms are reliably integrated
with an efficient semi-analytical technique. (iii) In the low-Mach number regime, 
standard explicit Godunov-type schemes lose efficiency and accuracy; the issue is addressed
by means of a simple semi-implicit, pressure-based, split treatment of acoustic and non-acoustic waves, 
again using staggered grids that recover the implicit solution for a single scalar field (the pressure) through a sequence of 
symmetric-positive definite linear systems that can be efficiently solved via the conjugate gradient method.

\end{abstract}

\section{Introduction}
\label{sec.intro}

This work is part of a sequence of papers \cite{GPRmodel, GPRmodelMHD, hstglm, sigpr, ilyanonnewtonian, tavellicrack, ptrsa, godbook}
dedicated to the extension and numerical solution of a first order hyperbolic Unified Model of Continuum Mechanics (UMCM), 
sometimes called in the literature HPR (Hyperbolic--Peshkov--Romenski) or GPR (from Godunov--Peshkov--Romenski).
Besides its first order hyperbolic character, the model in consideration has the very notable
feature of describing any continuum, ranging from elastic solids to viscous and inviscid fluids, 
in a unified formulation (hence the designation of Unified Model of Continuum Mechanics), with
only the material parameters mediating the distinction between solids and fluids within a single 
set of partial differential equations.

The model draws its origin in the work of Godunov and Romenski \cite{GodunovRomenski72,
GodRom1998,
GodRom2003, God1978, Godunov1996,Rom2001,SHTC-GENERIC-CMAT} on a formulation of Hyperelasticity in Eulerian coordinates, 
rather than in the Lagrangian frame more commonly adopted in the field. 
In \cite{PeshRom2014}, Peshkov and Romenski advanced the key insight  
that the Godunov--Romenski model may be applied not only to elasto-plastic solids, 
but to fluid flows as well, up to then unexplored in this framework, with the very notable exception of the 
early paper of Besseling \cite{Besseling1968}, in which an equivalent formalism was introduced. 
The same modelling approach is employed in 
\cite{FavrGavr2012,
Ndanou2014,
favrie2009solid,
Iollo2017, Barton2019,
hankplasticity,
Gavrilyuk2008,
FavrieGavrilyukSaurel,
NdanouFavrieGavrilyuk} for the study of elasto-plastic solids, and 
in \cite{Jackson2019, Jackson2019a} with explicit reference to the fluid applications
proposed by Peshkov and Romenski in \cite{PeshRom2014}.

In this paper, 
the Eulerian equations of hyperelasticity of Godunov and Romenski are 
employed for the purpose of introducing viscous-type forces
in the hyperbolic two-phase flow model with surface tension of Gavrilyuk and collaborators
\cite{Schmidmayer2017,Berry2008a}. Furthermore, in \cite{hstglm} the authors have shown that
the weakly hyperbolic model of surface tension given in \cite{Schmidmayer2017,Berry2008a}
cannot be solved by means of general purpose Godunov-type schemes without explicitly accounting 
for the involutive constraints present in the model, and, in the process, restoring hyperbolicity
through Godunov--Powell \cite{God1972MHD, Powell1997, powell1999} nonconservative symmetrising terms 
or by means of GLM \cite{MunzCleaning, Dedneretal} curl-cleaning.
Here the problem is tackled by adopting special \textit{compatible} discrete 
gradient and curl operators on a staggered grid like in \cite{boscherisigpr}, so that curl involutions are preserved
exactly up to machine precision by construction.
In this way (i.e. whenever curl-constraints are satisfied \textit{exactly}), 
the curl-cleaning and Godunov--Powell formulations of the model given in \cite{hstglm}
simply coincide exactly with the weakly hyperbolic one and viceversa, since only \textit{discrete} zeros
are introduced in either formulation.

With respect to the one proposed in \cite{boscherisigpr}, besides the application to a new PDE system, 
the scheme presents several novel aspects aimed at improving its robustness and accuracy. These concern 
in particular 
the implementation of momentum and energy updates due to the implicit pressure solver
included in the method, that interpolate the pressure field at the edge locations of the main
Finite Volume grid before computing pressure fluxes, as opposed to a two-way interpolation procedure from and to the main or dual grids.
Nonetheless, these modifications are introduced while at the same time preserving the 
staggered
semi-implicit split discretisation which enables 
the efficient application of the method to low Mach number flows.

Additionally, the scheme employs a 
semi-analytical time integration method for stiff nonlinear source governing strain relaxation, 
which represents a rather challenging task, especially in the context of two-phase flows.
This time integration approach extends the ideas introduced in \cite{chiocchettimueller, tavellicrack, boscheriulaggpr}
to the full equations of the Unified Model of Continuum Mechanics in the fluid regime, thanks
to a convenient polar decomposition of the stretch and rotation components of the distortion (or cobasis) field.

Moreover, the scheme is constructed in such a way that it can discretely preserve 
multiphase flows having uniform pressure and velocity, regardless of the distribution of density or volume fraction.

The paper is structured as follows: 
in Section~\ref{sec.model}, we briefly introduce the first order hyperbolic mathematical model 
for viscous two phase-flow with surface tension
considered in this work.
In Section~\ref{sec.method}, we detail the numerical method developed for the solution
of such a model, with particular attention given to the
discretisation of involution constraints and the split treatment of convective and acoustic phenomena.
In Section~\ref{sec:semianalytical}, we present a semi-analytical technique developed for the integration
of the relaxation sources at the basis of the viscous-like behaviour present in the Unified Model of Continuum Mechanics, 
and at the same time, we provide some details concerning the mathematical structure of the governing equations.
In Section~\ref{sec.results}, we show an extensive selection of numerical experiments, including convergence
results for the numerical method,
an experimental verification of the uniform flow compatibility of the scheme, formally proven in Section~\ref{sec.method}, 
an experimental verification of the curl involution compatibility, and large scale simulations of 
colliding droplets and gravity-driven \textit{multiphase} Rayleigh--Taylor instabilities.

In Section~\ref{sec.conclusion}, we summarise the conclusions that can be drawn from the results included in 
the paper, and we point towards some possible future research directions.

\section{First order hyperbolic model for two-phase viscous flow with surface tension} 
\label{sec.model} 

The basis of the governing 
equations studied in this paper formally resembles the Kapila system \cite{kapila2001} (sometimes called
five-equation Baer--Nunziato\cite{BaerNunziato1986} model) which describes two-phase flows
under pressure and velocity equilibrium hypothesis. In addition to the Kapila (sub)-system, 
a vector-valued equation for an interface field $\vec{b}$ is employed in order to track interfaces
and provide a hyperbolic description of surface tension \cite{Schmidmayer2017,hstglm,Berry2008a}. 
Furthermore, another matrix-valued field
$\vec{A}$, called \textit{distortion} field or \textit{cobasis}, 
is used to model the deformation of the fluid in consideration, 
which is actually described as a visco-elastoplastic solid with stiff strain relaxation, that
is, a stiff relaxation source is included in the governing equations such that the strain 
encoded in the components of the cobasis matrix $\vec{A}$.

The first order hyperbolic system for single velocity, single-pressure, compressible two-phase flow with
surface tension and hyperbolic viscosity then reads  
\begin{subequations} 
	\label{eq:gavrilyuk}
	\begin{align}
	&\partial_t \left(\alpha_1\,\rho_1\right) 
	+ \nabla \cdot \left(\alpha_1\,\rho_1\, \vec{u} \right) = 0, \label{gavr.mass1}\\[1mm] %
	&\partial_t\left(\alpha_2\,\rho_2\right) 
	+ \nabla \cdot \left(\alpha_2\,\rho_2\, \vec{u} \right) = 0, \label{gavr.mass2}\\[1mm] %
	&\partial_t\left(\rho\,\vec{u} \right) 
	+ \nabla \cdot \left(\rho\, \vec{u}\otimes\vec{u} + p\,\vec{I} - \vec{\Sigma}_\up{t} - \vec{\Sigma}_\up{s} \right) = 
		\rho\,\vec{g}, \label{gavr.mom}\\[1mm] 
	&\partial_t\left(\rho\,E\right) 
	+ \nabla \cdot \left[\left(\rho\,E + p\right) \vec{u} - \left(\vec{\Sigma}_\up{t} + \vec{\Sigma}_\up{s}\right)\cdot\vec{u} \right] = \rho\,\vec{g}\cdot\vec{u}, \label{gavr.energy}\\[1mm] 
	&\partial_t\left(\alpha_1\right) + \vec{u} \cdot \nabla \alpha_1 - K\, \nabla\cdot\vec{u} = 0,\label{gavr.alpha}\\[1mm] 
	&\partial_t\left(\vec{b}\right) + \nabla \left(\vec{b}\cdot\vec{u} \right) + \left( \nabla\vec{b} - \nabla\vec{b}^\transpose \right)\cdot\vec{u}  = \vec{0},\label{gavr.interface}\\[1mm] 
   &\partial_t\left(\vec{A}\right) + \nabla \left(\vec{A}\cdot\vec{u} \right) + \left( \nabla\vec{A} - \nabla\vec{A}^\transpose \right)\cdot\vec{u}  = \sour,\label{gavr.A}  
	\end{align}
\end{subequations}
where $\rho_1$ and $\rho_2$ are the mass density of phase number $1$ and $2$ respectively, $\alpha_1$ and $\alpha_2$ 
their volume fraction; the velocity vector is denoted by $\vec{u}$; 
$\rho E$ is the total energy density, 
$p$ is the mixture pressure and the mixture density is given by $\rho = \alpha_1\,\rho_1 + \alpha_2\,\rho_2$.

The first two equations (\ref{gavr.mass1}, \ref{gavr.mass2}) state mass conservation for each phase, 
Equation \eqref{gavr.mom} is the momentum balance for the mixture, with the source term accounting
for constant gravitational acceleration $\vec{g}$. The energy balance law is \eqref{gavr.energy} and 
Equation \eqref{gavr.alpha} governs the evolution of the volume fraction of the first phase, which due to 
the saturation constraint $\alpha_1 + \alpha_2 = 1$ is sufficient to characterise the dynamics of both volume fractions.

The vector $\vec{b}$ in \eqref{gavr.interface} is an \textit{interface field} which 
represents the gradient of a scalar colour function $c$ and its components are evolved \emph{independently}
as three separate state variables, rather than computed as the discrete gradient of $c$.
Due to the field $\vec{b}$ representing the gradient of a scalar, 
it must be curl-free, i.e. it is required to satisfy 
\begin{equation}
   \nabla \vec{b} - \nabla \vec{b}^\transpose = \vec{0}.
   \label{eqn.curlfree} 
\end{equation}
It can easily be checked that if \eqref{eqn.curlfree} holds at the initial time, then due 
to \eqref{gavr.interface} it remains curl-free for all times.

The vectors $\vec{a}_{1}$, $\vec{a}_{2}$, and $\vec{a}_{3}$ are the rows of a three-by-three 
nonsymmetric matrix $\vec{A}$, here called \emph{distortion}, 
are the rows of a
and $\vec{Z}$ is a three-by-three \emph{strain relaxation} source term 
\begin{equation}
   \sour = - \frac{3}{\tau}\,{\left(\det{\vec{A}}\right)}^{5/3}\,\vec{A}\,\dev{\left(\vec{A}^\transpose\,\vec{A}\right)}.
\end{equation}
The homogeneous part of \eqref{gavr.A} shares the exact same 
structure with Equation \eqref{gavr.interface} and thus each one of the rows or $\vec{A}$ 
can be discretised in the same manner.
The presence of the source term, however, couples the three vector equations for the distortion matrix, 
which means that, unlike the spatial discretisation, the complete time 
integration of these equations cannot be carried out independently of each other.

The parameter $\tau$ is called \emph{strain relaxation time} for the mixture and is computed from a 
logarithmic interpolation $\tau = \tau_1^{\alpha_1}\,\tau_2^{1 - \alpha_1}$, where the relaxation
times of each phase are set as
$\tau_1 = 6\,\nu_1/c_\up{s}^2$ and $\tau_2 = 6\,\nu_2/c_\up{s}^2$ to fit the kinematic viscosities $\nu_1$ and $\nu_2$ of
the two fluids.
For simplicity, in this work we take $c_\up{s}$, a parameter representing the propagation speed of 
small-amplitude shear waves, to be common for both phases.

For each phase we assume 
that the equation of state has the following form 
\begin{equation}
p_1 = \left(\gamma_1 - 1\right)\,\rho_1\,e_1 - \gamma_1\,\Pi_1, \qquad p_2 = 
\left(\gamma_2 - 1\right)\,\rho_2\,e_2 - \gamma_2\,\Pi_2,
\end{equation}
with $\gamma_1$, $\gamma_2$, $\Pi_1$, $\Pi_2$ given parameters of the equation of state
and $\rho\,e_1$ and $\rho\,e_2$ the internal energy densities of the two phases. Due to the
pressure equilibrium assumption $p_1 = p_2 = p$, the mixture equation of state then reads 
\begin{equation}
p = \frac{
	\rho\,e\,\left(\gamma_1 - 1\right)\left(\gamma_2 - 1\right) - \alpha_1\,\gamma_1\,\Pi_1\,
	\left(\gamma_2 - 1\right) - \alpha_2\,\gamma_2\,\Pi_2\,\left(\gamma_1 - 1\right)
}{
	\alpha_1\,\left(\gamma_2 - 1\right) + \alpha_2\,\left(\gamma_1 - 1\right)
},
\end{equation}
where $\rho\,e$ is the internal energy density of the mixture. 
Furthermore, for this choice of closure relation, we have
\begin{equation}
        K = \frac{\alpha_1\,\alpha_2\,\left(\rho_2\,a_2^2 - \rho_1\,a_1^2\right)}{\alpha_1\,
	\rho_2\,a_2^2 + \alpha_2\,\rho_1\,a_1^2},
 \label{eqn.K} 
\end{equation}
with 
\begin{equation*}
a_1 = \sqrt{\frac{\gamma_1\,\left(p + \Pi_1\right)}{\rho_1}}, \qquad \textnormal{and} \qquad 
a_2 = \sqrt{\frac{\gamma_2\,\left(p + \Pi_2\right)}{\rho_2}}.
\end{equation*}
According to \cite{Schmidmayer2017} the capillary stress tensor is given by  
\begin{equation}
  \vec{\Sigma}_\up{t} = -\sigma\,\norm{\vec{b}}\left( \dfrac{\vec{b} \otimes \vec{b}}{\norm{\vec{b}}^2} - \vec{I} \right), 
 \label{eqn.Omega} 
\end{equation}
where $\sigma$ is a constant that characterizes surface tension. 
Furthermore, the total energy density reads in terms of the other state variables 
\begin{equation}
\begin{aligned}
   \rho\,E &= \rho\,e + \rho\,e_\up{s} + \rho\,e_\up{t} + \rho\,e_\up{k} = \rho\,e + \rho\,e_2 + \rho\,e_\up{k} = \\
   &= 
   \rho\,e + \frac{1}{4}\,\rho\,c_\up{s}^2\,\tr{\left(\dev\vec{G}\,\dev{\vec{G}}\right)} + 
   \sigma\,\norm{\vec{b}} + \frac{1}{2}\,\rho\,\norm{\vec{u}}^2,
\end{aligned}
\end{equation}
with $\rho\,e_\up{s} = \rho\,c_\up{s}^2\,\tr{\left(\dev\vec{G}\,\dev{\vec{G}}\right)}/4$ the energy 
associated with elastic/shear stress, $\rho\,e_\up{t} = \sigma\,\norm{\vec{b}}$ the surface energy 
density, and $\rho\,e_\up{k} = \rho\,\norm{\vec{u}}^2/2$ the kinetic energy density.
For this choice of elastic energy potential, the 
elastic stress tensor reads
\begin{equation}
   \vec{\Sigma}_\up{s} = -\rho\,c_\up{s}^2\,\vec{G}\,\dev{\vec{G}},
\end{equation}
where $\vec{G} = \vec{A}^\transpose\,\vec{A}$ is the \emph{metric tensor} for the fluid mixture,
which in the SHTC formalism is guaranteed to be symmetric and positive definite
by construction at the discrete level, 
since it is computed from the distortion matrix $\vec{A}$ rather 
than evolved directly.

\section{Numerical method} 
\label{sec.method}
\input{method}
\input{methodsource}



\section{Numerical results} 
\label{sec.results}

\input{results}
\section{Conclusion}
\label{sec.conclusion}
In this paper, we have presented a novel semi-implicit scheme for the solution of two-phase flow
with surface tension and viscous effects. Both surface tension and viscosity are modelled within the 
framework of Symmetric Hyperbolic Thermodynamically Compatible (SHTC) systems and thus the governing
equations are of first order in space and time and have hyperbolic character. 

The scheme makes use of a multiply-staggered Cartesian grid, featuring
traditional Finite Volume cell average 
quantities which can be interpreted (up to second order errors) as collocated cell center values, as well as edge-based
quantities (momentum components) that are used for the discretisation of pressure forces (the acoustic subsystem), 
making the scheme \textit{well suited for simulating low Mach number flows}, since such a staggered configuration
allows the computation of the pressure field, along with its associated forces and work terms, 
as the solution of a sequence of symmetric positive definite linear 
systems \cite{MunzPark, CasulliCompressible, DumbserCasulli2016}.

Additionally, since the governing equations must obey a set of curl-type involutions 
(notably the interface field $\vec{b}$ is
to remain curl-free if initially curl-free), 
several corner quantities are employed in order to achieve a curl-free discretisation 
of the distortion field $\vec{A}$ (also called cobasis in solid mechanics), like in \cite{boscherisigpr}, 
and in this work the same discretisation is also applied to an interface field $\vec{b}$ that
tracks interfaces and is used for the computation of surface tension forces.

A major result of this work is that curlfree discretisation can be used to solve the surface tension 
model of S. Gavrilyuk \cite{Berry2008a, Schmidmayer2017} in its original 
weakly hyperbolic formulation, without GLM curl-cleaning or Godunov--Powell symmetrising terms \cite{hstglm}.
With the curlfree approach adopted in this paper, the interface field $\vec{b}$ can be evolved \textit{directly}
by means of the original weakly hyperbolic formulation of S. Gavrilyuk and collaborators \cite{Berry2008a, Schmidmayer2017}
without the instabilities associated with weak hyperbolicity, highlighted and solved in \cite{hstglm} by
providing an alternative (hyperbolic) formulation of the governing equations. It is worth noting that, 
if the constraint $\nabla\times\vec{b} = \vec{0}$ is satisfied \textit{exactly}, all formulations of 
the surface tension model of Gavrilyuk collapse onto the same set of partial differential equations, since
the non-conservative Godunov--Powell terms vanish, and so does the auxiliary GLM curl cleaning field 
introduced in \cite{hstglm}.

Several computational experiments, including very high resolution simulations of Rayleigh--Taylor 
instabilities at different Weber and Reynolds number, have been carried out and employed to highlight 
the efficacy of the numerical scheme.

Future applications will include the extension of the method to uniform high order in space and time, 
by means of IMEX timestepping \cite{imex1, imex2, imex3, imex4, imex5, imex6, imex7}, as well as further implicit discretization of 
shear waves, aimed at eliminating the timestep restriction (roughly proportional to the 
shear wavespeed $\cshear$) which can be very high in certain fluids, for example when the model is
applied to lava flows.


\section*{Acknowledgments}

The research presented in this paper has been funded by the Italian Ministry of Education, University and Research (MIUR) 
in the frame of the Departments of Excellence Initiative 2018--2022 attributed to DICAM of the University of Trento (grant L. 232/2016) and in the frame of the 
PRIN 2017 project \textit{Innovative numerical methods for evolutionary partial differential equations and  applications}. 
The research was also funded by the 
Deutsche Forschungsgemeinschaft (DFG) via the project DROPIT, grant no. GRK 2160/2. 

Furthermore, M.~D. has also received funding from the University of Trento via the Strategic 
Initiative \textit{Modeling and Simulation}.

%
%

\bibliographystyle{spmpsci}      
\bibliography{biblio}  

\end{document}

%% file: method.tex
The scheme proposed in this work is based on a \textit{multiply staggered} 
Cartesian discretisation introduced in \cite{boscherisigpr}
that employs special gradient and curl operators that can be used to evolve sensitive involution-constrained
quantities such as the interface field that is used to track material interfaces and compute surface tension forces.
This allows to solve the weakly hyperbolic surface tension model given in \cite{Schmidmayer2017}
without curl cleaning procedures, due to the fact that, if curl involutions are enforced 
exactly, the weakly hyperbolic model and its augmented hyperbolic variants collapse onto the same stable system.

The staggered grid also allows to discretely recover a numerical scheme for 
the incompressible Navier--Stokes equations in the low Mach limit \cite{KlaMaj,KlaMaj82,Klein2001,Munz2003,MunzDumbserRoller}, 
and at the same time 
eliminate the accuracy and efficiency issues of Godunov-type methods: the timestep restriction
due to acoustic waves is lifted at the rather limited cost of solving a sequence of 
sparse symmetric-positive-definite systems for the scalar pressure at each timestep.
Together with a suitable low-dissipative discretisation of convective fluxes, this enables
the second-order spatially accurate method to achieve results that are comparable with 
schemes that formally feature a much higher order of accuracy.

Finally, the method employs a semi-analytical integration technique for the 
strain relaxation sources of the unified model of continuum mechanics, introduced in Section \ref{sec:semianalytical}, 
in order to robustly compute viscous forces in flows with high spatial and parametric variability.  

\subsection{Flux-splitting approach}
The semi-implicit curl-preserving numerical method presented in this work
is based on a three way split of the governing equations, such that convection, 
capillarity, strain, and acoustic effects are treated separately, each with an ad-hoc discretisation.
In particular: (\textit{i}) a MUSCL--Hancock scheme with primitive variable reconstruction and positivity preserving limiting
is adopted for the solution of the convective part of the system;
(\textit{ii}) an implicit staggered conservative scheme is used to compute the unknown pressure field 
as the solution of a simple discrete wave equation leading to a well-behaved symmetric positive definite system
of linear equations. This lifts timestep restrictions due to acoustic waves and preserves the accuracy
of the method in the low-Mach regime; (\textit{iii}) the evolution of geometric involution constrained fields
associated with material distortion $\vec{A}$ and material interfaces $\vec{b}$ is carried out with
ad-hoc discrete differential operators that can preserve the curl involutions of the governing equations 
exactly up to machine precision.

In sequence, at each timestep, first, the convective update of the conserved variable is computed via a path-conservative MUSCL--Hancock 
method, and at the same time the quantities $\vec{A}$ and $\vec{b}$, that are endowed with curl constraints, are evolved in time
with a simple two-stage Runge--Kutta scheme, which adopts the semi-analytical solver introduced in Section \ref{sec:semianalytical}
for strain relaxation.
Following this, corner fluxes due to viscous forces and capillarity can be computed.
Then, a discrete wave equation, derived from a staggered discretisation of the momentum-energy system,
can be solved for the unknown scalar pressure and finally, since a robust predictor for all state
variables has been obtained, momentum and energy interface fluxes can be computed and used to update
the conserved variables to the next time level in a conservative fashion.

The system \eqref{eq:gavrilyuk} can be written with matrix-vector notation as    
\begin{equation}
\label{eqn.pde} 
\partial_t \vec{Q} + \nabla \cdot \vec{F}(\vec{Q}) + \vec{B}(\vec{Q}) \cdot \, \vec{Q} = \vec{S}(\vec{Q}), 
\end{equation}
with the state vector
\begin{equation}
    \vec{Q} = \left(\alpha_1\,\rho_1,\ \alpha_2\,\rho_2,\ \rho\,\vec{u},\ \rho\,E,\ \alpha_1,\ \vec{b},\ \vec{A}\right)^\transpose,
\end{equation}
a flux tensor $\vec{F}(\vec{Q})$, a non-conservative product $\vec{B}(\vec{Q}) \, \nabla \vec{Q}$, 
and an algebraic source term $\vec{S}(\vec{Q})$. 
As proposed in \cite{DumbserCasulli2016,SIMHD,boscherisigpr} the flux tensor is \emph{split} 
into a pressure part, and a convective part, to be treated partially 
by means of a path-conservative 
MUSCL--Hancock scheme, partially with a special \textit{compatible} and 
structure-preserving discretization using a vertex-based staggered grid. 
Hence, 
eqn.~\eqref{eqn.pde} is rewritten as  
\begin{equation}
\begin{aligned}
\label{eqn.pde.split} 
\partial_t \vec{Q}+ \nabla \cdot &\left[\fll\mathbf{F}_c(\vec{Q}) + \mathbf{F}_p(\vec{Q}) + 
\mathbf{F}_v(\vec{Q}) \right] + \nabla \left[\fll\gradflux(\vec{Q})\right]\, + \\
+&\left[ \fll\vec{B}_c(\vec{Q})  + \vec{B}_v(\vec{Q}) \right] \cdot \nabla \vec{Q}= \vec{S}_p(\vec{Q}) + \vec{S}_v(\vec{Q}), 
\end{aligned}
\end{equation}
where $\vec{F}_c(\vec{Q}) $ and $\vec{B}_c(\vec{Q})$ contain the convective terms that will be  discretized explicitly; 
$\vec{F}_p(\vec{Q})$ are the pressure fluxes that will be discretized implicitly using an 
edge-based staggered grid. The resulting splitting into pressure and convective fluxes is 
identical to the flux-vector splitting scheme of Toro and V\'azquez-Cend\'on recently 
forwarded in \cite{ToroVazquez}. The terms $\mathbf{F}_v(\vec{Q}) $, $\nabla \gradflux(\vec{Q})$ and $\vec{B}_v(\vec{Q}) 
 \nabla \vec{Q}$ are discretized in a structure-preserving manner using an explicit scheme on a vertex-based staggered mesh. 
The split fluxes read 
\begin{equation} 
\label{eqn.split.def} 
\mathbf{F}_c = \left( \begin{array}{c} \alpha_1\,\rho_1\,\vec{u} \\
 \alpha_2\,\rho_2\,\vec{u}  \\
 \rho\,\vec{u}\otimes\vec{u} \\
 \left(\rho\,E - \rho\,e\right)\,\vec{u} \\
 0 \\
 \vec{0} \\
 \vec{0}  \end{array} \right),  \  
\mathbf{F}_p = \left( \begin{array}{c} 0 \\
 0 \\
 p\,\vec{I} \\
 (\rho\,e + p)\,\vec{u} \\
 0 \\
 \vec{0} \\
 \vec{0} \end{array} \right), \  
\mathbf{F}_v = \left( \begin{array}{c} 0 \\
 0 \\
 -\vec{\Sigma}_\up{t} - \vec{\Sigma}_\up{s}  \\
 \left(-\vec{\Sigma}_\up{t} - \vec{\Sigma}_\up{s}\right)\,\vec{u} \\
 0 \\
 \vec{0} \\
 \vec{0}  \end{array} \right),
\end{equation} 
The convective part of the non-conservative product is given by 

\begin{equation}
    \vec{B}_c(\vec{Q})\,\nabla\vec{Q} = \left(0,\ 0,\ \vec{0},\ 0,\ \vec{u}\cdot\nabla\alpha_1 -K\,\nabla\cdot\vec{u},\ \vec{0},\ \vec{0}\right)^\transpose,
\end{equation}
The evolution of the curl-free vector field $\mathbf{b}$ and the distortion matrix $\vec{A}$ is governed by the terms 
$\nabla \gradflux(\vec{Q})$, $\vec{B}_v(\vec{Q}) \, \nabla \vec{Q}$, and $\vec{S}_v(\vec{Q})$, with 
\begin{equation} 
\gradflux(\vec{Q}) = \left( \begin{array}{c} 0 \\ 0 \\ \vec{0} \\ 0 \\ 0 \\ 
    \vec{b}\cdot\vec{u}  \\ \vec{A}\,\vec{u}\end{array} \right),  
\quad 
\vec{B}_v(\vec{Q}) \, \nabla \vec{Q} = \left( \begin{array}{c} 0 \\ 0 \\ \vec{0} \\ 
0 \\ 0 \\ 
\left(\nabla\vec{b} - \nabla\vec{b}^\transpose\right)\,\vec{u}  \\ 
\left(\nabla\vec{A} - \nabla\vec{A}^\transpose\right)\,\vec{u}  \\ 
\end{array} \right), 
\end{equation} 
and 
\begin{equation}
    \vec{S}_v(\vec{Q}) = \left(0,\ 0,\ \vec{0},\ 0,\ 0 ,\ \vec{0},\ \vec{Z}\right)^\transpose,
\end{equation}
which accounts for the strain relaxation effects.
For clarity, we specify that the tensor 
notation $\left(\nabla\vec{A} - \nabla\vec{A}^\transpose\right)\,\vec{u}$ yields a three by three matrix whose
rows $\vec{a}_k$ are the obtained by computing $\left(\nabla\vec{a}_k - \nabla\vec{a}_k^\transpose\right)\,\vec{u}$
for each row $\vec{a}_k$ from $\vec{a}_1$ to $\vec{a}_3$.
The remaining source terms, corresponding to gravity forces, are included in the pressure sub-system as
\begin{equation}
    \vec{S}_p(\vec{Q}) = \left(0,\ 0,\ \vec{g},\ \vec{g}\cdot\vec{u},\ 0 ,\ \vec{0},\ \vec{0}\right)^\transpose.
\end{equation}
As already mentioned before, the subsystem  
\begin{equation}
\label{eqn.pde.ex} 
\partial_t \vec{Q}+ \nabla \cdot \left[\fll \mathbf{F}_c(\vec{Q}) + \mathbf{F}_v(\vec{Q}) \right] + 
\nabla \left[\fll\gradflux(\vec{Q})\right] + \left[\fll \vec{B}_c(\vec{Q}) + \vec{B}_v(\vec{Q}) \right] \, \nabla \vec{Q}= \vec{S}_v(\vec{Q}), 
\end{equation}
will be discretized explicitly. The discretization method 
presented in the next section is a combination of a classical second order 
MUSCL--Hancock type \cite{vanleer1974, vanleer1979, torobook} TVD finite volume scheme for the convective fluxes $\vec{F}_c$ and the 
nonconservative term $\vec{B}_c\,\nabla\vec{Q}$, a 
curl-free discretization for the terms $\gradflux$, and $\vec{B}_v \, \nabla \vec{Q}$ using compatible gradient 
and curl operators, a semi-analytical integration technique for
the relaxation source $\vec{S}_v$, as well as a vertex-based discretization of the terms $\vec{F}_v$. 
The pressure subsystem 
\begin{equation}
\label{eqn.pde.im} 
\partial_t \vec{Q} + \nabla \cdot \vec{F}_p(\vec{Q}) = \vec{S}_p(\vec{Q})
\end{equation}
is formally identical to the  Toro--V\'azquez pressure system \cite{ToroVazquez}, 
with the simple addition of gravity sources $\vec{S}_p$ and is discretized implicitly.

\subsection{Eigenvalue estimates}
Due to the large size of the hyperbolic PDE system coupling convective, acoustic, thermal, shear, and capillarity effects, 
it is at the present time \textit{impossible} to explicitly compute the eigenstructure of the full system, 
that is, the eigenvalues and eigenvectors of the matrix $\vec{C}_1 = \vec{C}\cdot\uvec{n}$, called in
an abuse of terminology, Jacobian of the system in the direction specified by the unit vector $\uvec{n}$.
Formally, 
$\vec{C}$ is defined as 
\begin{equation}
    \vec{C} = \left(\pd{\vec{Q}}{\vec{V}}\right)^{-1}\,\pd{\vec{F}}{\vec{V}} + \vec{B}\,\pd{\vec{Q}}{\vec{V}}, 
\end{equation}
the quasi-linear form of the governing equations
\begin{equation}
    \partial_t\vec{V} + \vec{C}\,\nabla\vec{V} = \vec{S},
\end{equation}
recalling that
the general first order balance law reads
\begin{equation}
    \partial_t \vec{Q} + \nabla\cdot\vec{F} + \vec{B}\,\nabla\vec{Q} = \vec{S}.
\end{equation}
As a matter of fact, even the eigenvalues of $\vec{C}$ can be obtained in closed form only with several simplificative assumptions, 
such as setting the distortion matrix $\vec{A}$ identity or the thermal impulse vector $\vec{J}$ to null.
Moreover, even the numerical computation of the eigenvalues of the system can become prohibitively expensive
as the Jacobian matrix of the fully coupled system for two-phase flow with shear, surface tension, and heat conduction
is large enough that application of standard numerical eigenvalue methods whenever a Riemann solver
of Rusanov of HLL type is
evaluated is rather wasteful from a computational standpoint. 

Since the Rusanov flux only requires an estimate of the maximum absolute eigenvalue $\lambda_\up{max}$ 
of the system Jacobian, 
one might be tempted to employ a power iteration method (Von Mises iteration or Rayleigh quotient iteration) 
for the computation of the spectral radius of the matrix.
However, even the \textit{evaluation}
of the Jacobian matrix for a given state vector is arithmetically very intensive and thus it is preferable
to avoid the procedure entirely, favouring instead simpler estimates that can be computed directly from the
state vector. 

A practical and effective choice, which we found to be 
rather safe for the estimation of the spectral radius of the Jacobian matrix of the full system,
while only leading to occasional mild overestimates, is setting
\begin{align}
   &\lambda_\up{max} = \max\left(|\vec{u}\cdot\uvec{n} + \lambda|,\ |\vec{u}\cdot\uvec{n} - \lambda|\right)
\shortintertext{with}
   &\lambda = \sqrt{\lambda_\up{ps}^2 + \lambda_\up{t}^2},
\end{align}
where $\lambda_\up{ps}$ accounts for mixed pressure/shear waves, and $\lambda_\up{t}$ is an estimate of the contribution due to capillarity waves only.
In principle, the eigenvalues of the full model including two-phase flow, shear, surface tension and heat conduction, 
appear to couple all effects in mixed thermal-shear-pressure-capillarity waves, and the same is true for 
the surface tension sub-system yielding acoustic-capillary waves, the shear sub-system yielding pressure-shear waves, 
and the heat conduction sub-system yielding thermo-acoustic waves, thus it is impossible to rigorously
assign a wavespeed to only one of the effects.
Nonetheless, by means of numerical experimentations, we found that surprisingly robust and accurate estimates of the 
maximum absolute eigenvalue of the system can be achieved with appropriate choices of the 
two estimates $\lambda_\up{ps}$,
and $\lambda_\up{t}$, which are given in the following paragraphs.

\subsubsection{Wavespeed estimate for capillarity waves}
We begin discussing the eigenvalue estimates for the capillary system, 
Since the full eigenstructure of the capillarity sub-system
\begin{equation}
\left\{
   \begin{aligned}
              &\partial_t\left(\alpha_1\,\rho_1\right) 
            + \nabla\cdot\left(\alpha_1\,\rho_1\,\vec{u}\right) = 0, \\[1mm] %
        &\partial_t\left(\alpha_2\,\rho_2\right) 
            + \nabla\cdot\left(\alpha_2\,\rho_2\,\vec{u}\right) = 0, \\[1mm] %
        %
        %
        &\partial_t\left(\rho\,\vec{u}\right) 
            + \nabla\cdot\left(\rho\,\vec{u}\otimes\vec{u} + p\,\vec{I} + \vec{\Omega}\right) = 
            \vec{0}, \\[1mm] 
        &\partial_t\left(\rho\,E\right) 
            + \nabla\cdot\left[\left(\rho\,E + p\right)\,\vec{u} + \vec{\Omega}\,\vec{u}\right] = 
            0, \\[1mm] 
        &\partial_t\left(\alpha_1\right) + \vec{u} \cdot \nabla \alpha_1  - K\,\nabla\cdot\vec{u} = 
        0, \\[1mm] 
        &\partial_t\left(\vec{b}\right) + \left(\nabla\vec{b}\right)\,\vec{u} + 
        {(\nabla\vec{u})}^\transpose\,\vec{b} = \vec{0}, \qquad \nabla\times\vec{b} = \vec{0},
   \end{aligned}
   \right.
\end{equation}
can be computed analytically as shown in \cite{hstglm}.
A simple estimate for the speed associated with 
capillarity waves $\lambda_\up{t}$, can be computed analytically by recognising
different velocities (i.e. the adiabatic sound speed) in the expression of the eigenvalues for the system of two-phase flow with surface tension,
as already formally detailed in \cite{hstglm}.
For clarity, the capillarity tensor in the form
\begin{equation} 
  \vec{\Omega} = \sigma\,\left({
      {\norm{\vec{b}}^{-1}}\,\vec{b} \otimes \vec{b}} - \norm{\vec{b}}\,\vec{I}\right), 
\end{equation}
and $K$ is computed as $K = \left(\rho_2\,a_2^2 - \rho_1\,a_1^2\right)\,\alpha_1\,\alpha_2/(\alpha_1\,\rho_2\,a_2^2 + \alpha_2\,\rho_1\,a_1^2)$.
%
We generically denote with $a$ the mixture speed of sound (adiabatic pressure waves), 
which, for Kapila-type models such as those considered in 
this section, is the Wood \cite{wood1930} speed of sound which
reads

\begin{equation}
    a = \sqrt{\frac{\rho_1\,a_1^2\,\rho_2\,a_2^2}{\rho\,\left(\alpha_1\,\rho_2\,
        a_2^2 + \alpha_2\,\rho_1\,a_1^2\right)}},
\end{equation}
with the frozen soundspeeds $a_1$ and $a_2$ computed as
\begin{equation}
        a_1 = \sqrt{\rho_1^{-1}\,{\gamma_1\,\left(p + \Pi_1\right)}\fll}, \qquad
          a_2 = \sqrt{\rho_2^{-1}\,{\gamma_2\,\left(p + \Pi_2\right)}\fll},
\end{equation}
due to the stiffened gas equation of state being adopted for both phases, 

The simplest estimate for the speed associated with 
capillarity waves $\lambda_\up{t}$, is then given by
\begin{equation}
    \lambda_\up{t} = a_\sigma = \sqrt{{\fll{\rho}^{-1}\,\sigma}\,\norm{\vec{b}}\,\left[1 - {\norm{\vec{b}}^{-2}}\,{\left(\vec{b}\cdot\uvec{n}\right)^2}\right]}.
\end{equation}
We refer to \cite{hstglm} for a detailed derivation of the eigenstructure of the capillarity sub-system.

\subsubsection{Wavespeed of large amplitude pressure-shear waves}
The sub-system for two-phase flow with viscosity, which of course can also model hyperelastic solids by simply settings
the relaxation time $\tau\to\infty$, reads 
\begin{equation}
\label{eq.eigspsystem}
\left\{
   \begin{aligned}
              &\partial_t\left(\alpha_1\,\rho_1\right) 
            + \nabla\cdot\left(\alpha_1\,\rho_1\,\vec{u}\right) = 0, \\[1.5mm] %
        &\partial_t\left(\alpha_2\,\rho_2\right) 
            + \nabla\cdot\left(\alpha_2\,\rho_2\,\vec{u}\right) = 0, \\[1mm] %
        %
        %
        &\partial_t\left(\rho\,\vec{u}\right) 
            + \nabla\cdot\left(\rho\,\vec{u}\otimes\vec{u} + p\,\vec{I} + \rho\,c_\up{s}^2\,\vec{G}\,\dev\vec{G}\right) = 
            \vec{0}, \\[0.5mm] 
        &\partial_t\left(\rho\,E\right) 
            + \nabla\cdot\left[\left(\rho\,E + p\right)\,\vec{u} + \rho\,c_\up{s}^2\,\vec{G}\,\dev\vec{G}\,\vec{u}\right] = 
            0, \\[0.5mm] 
        &\partial_t\left(\alpha_1\right) + \vec{u} \cdot \nabla \alpha_1  - K\,\nabla\cdot\vec{u} = 
        0, \\[0.5mm] 
        &\partial_t\left(\vec{A}\right) + \left(\nabla\vec{A}\right)\,\vec{u} + 
        \vec{A} \,{(\nabla\vec{u})}= -3\,\tau^{-1}\,\left(\det\vec{A}\right)^{5/3}\,\vec{A}\,\dev\vec{G}.
   \end{aligned}
\right.
\end{equation}
By algebraic manipulation of the generalised Jacobian matrix appearing in the quasilinear form of \eqref{eq.eigspsystem}, 
one can find that the eigenvalues do not depend directly on the nine components of $\vec{A}$, but 
only on the metric tensor $\vec{G} = \vec{A}^\transpose\,\vec{A}$.
In three space dimensions, i.e. for a general 
metric tensor of the form
\begin{equation}
   \vec{G} = \begin{pmatrix}
   G_{11} & G_{12} & G_{13} \\
   G_{12} & G_{22} & G_{23} \\
   G_{13} & G_{23} & G_{33} \\
   \end{pmatrix},
\end{equation}
one has that the three eigenvalues $\lambda_1$, $\lambda_2$, and $\lambda_3$ associated with mixed pressure/shear waves 
are the square roots of the eigenvalues of the symmetric matrix
\begin{equation}
\label{eqn.matrixm}
   \vec{M} = \begin{pmatrix}
   m_{11} & m_{12} & m_{13} \\
   m_{12} & m_{22} & m_{23} \\
   m_{13} & m_{23} & m_{33} \\
   \end{pmatrix},
\end{equation}
with 
\begin{equation}
\begin{aligned}
    &m_{11} = c_\up{s}^2\,\left[10\,G_{11}^2 + 9\,\left(G_{12}^2 + G_{13}^2\right) - 3\,G_{11}\,\left(G_{22} + G_{33}\right)\right]/3 + a^2,\\
    &m_{22} = c_\up{s}^2\,\left\{4\,G_{12}^2 + 3\,G_{23}^2 + G_{22}\,\left[2\,\left(G_{11} + G_{22}\right) - G_{33}\right]\right\}/3,\\
    &m_{33} = c_\up{s}^2\,\left\{4\,G_{13}^2 + 3\,G_{23}^2 + G_{33}\,\left[2\,\left(G_{11} + G_{33}\right) - G_{22}\right]\right\}/3,\\
    &m_{12} = 2\,c_\up{s}^2\,\left[3\,G_{13}\,G_{23} + G_{12}\,\left(4\,G_{11} + 2\,G_{22} - G_{33}\right)\right]/3,\\[1.0mm]
    &m_{13} = 2\,c_\up{s}^2\,\left[3\,G_{12}\,G_{23} + G_{13}\,\left(4\,G_{11} + 2\,G_{33} - G_{22}\right)\right]/3,\\[1.0mm]
    &m_{23} = 2\,c_\up{s}^2\,\left[2\,G_{12}\,G_{13} + G_{23}\,\left(G_{11} + G_{22} + G_{33}\right)\right]/3,
\end{aligned}
\end{equation}
Note that the above expression for the components of $\vec{M}$ is associated with the Jacobian matrix of the system
projected along the $x$-axis direction, i.e. for $\uvec{n} = {(1,\ 0,\ 0)}^\transpose$. However, due to 
the rotational invariance of the governing equations, it is always possible to define a local reference frame in which the 
new $x$ direction is that along which the directional eigenvalue estimate is sought, that would be along any of 
the outward normal vectors of
the space-time faces when computing approximate Riemann fluxes.
The characteristic polynomial associated with $\vec{M}$ then reads
\begin{equation}
\label{eq.charpolps}
   (\lambda^2)^3 - (\lambda^2)^2\,\tr{\vec{M}} - \lambda^2\,\left[\tr\left(\vec{M}\,\vec{M}\right) - 
   \left(\tr\vec{M}\right)^2\right]/2 - \det{\vec{M}} = 0
\end{equation}
and can be solved analytically for $\lambda^2$ by means of the Del Ferro--Tartaglia--Cardano procedure.
However, from a computational standpoint, it is much more efficient and accurate to apply 
the Jacobi method to the symmetric matrix $\vec{M}$ directly.
In any case we formally set the wavespeed estimate due to mixed pressure/shear waves
to be $\lambda_\up{ps} = \max\left(\lambda_1,\ \lambda_2,\ \lambda_3\right)$, 
with $\lambda_1^2,\ \lambda_2^2,\ \lambda_3^2$ obtained by solving \eqref{eq.charpolps}.

In two space dimensions, the $G_{13}$ and $G_{23}$ components of the metric tensor $\vec{G}$ vanish, 
thus the eigenvalues of the system 
can be found as the square roots of the eigenvalues
of a simplified matrix
\begin{equation}
   \vec{M} = \begin{pmatrix}
   m_{11} & m_{12} & 0 \\
   m_{12} & m_{22} & 0 \\
   0 & 0 & m_{33} \\
   \end{pmatrix},
\end{equation}
for which closed form expression are easy to compute and write.
In particular the components of \eqref{eqn.matrixm} simplify to 
\begin{equation}
\label{eqn.matrixm2d}
   \begin{aligned}
    & m_{11} = c_\up{s}^2\,\left\{4\,G_{11}^2 + 9\,\left[G_{12}^2 + G_{11}\,\left(G_{11} - \tr\vec{G}/3\right)\right]\right\}/3 + a^2,\\
    & m_{22} = c_\up{s}^2\,\left[G_{22}\,\left(2\,\tr\vec{G} - 3\,G_{33}\right) + 4\,G_{12}\,G_{12}\right]/3,\\[1mm]
    & m_{33} = c_\up{s}^2\,G_{33}\,\left(2\,\tr\vec{G} - 3\,G_{22}\right)/3,\\[1mm]
    & m_{12} = 2\,c_\up{s}^2\,G_{12}\,\left(4\,G_{11} + 2\,G_{22} - G_{33}\right)/3,
   \end{aligned}
\end{equation}
and the eigenvalues of the sub-model for two-phase flow with shear, 
associated with pressure/shear waves are
\begin{equation}
   \lambda_1 = \sqrt{\fls m_4 + m_5},\quad
   \lambda_2 = \sqrt{\fls m_4 - m_5},\quad
   \lambda_3 = \sqrt{\fls m_{33}}, 
\end{equation}
   with
\begin{equation}
   m_4 = \left(m_{11} + m_{22}\right)/2,\quad m_5 = \sqrt{m_4^2 + m_{12}^2 - m_{11}\,m_{22}}.
\end{equation}
For small deformations, i.e. when $\vec{G}\to\vec{I}$,
it is easy to verify that the components of \eqref{eqn.matrixm2d} further simplify
\begin{equation}
    m_{11} = a^2 + 4\,c_\up{s}^2/3,\quad
    m_{22} = c_\up{s}^2,\quad
    m_{33} = c_\up{s}^2,\quad
    m_{12} = 0,
\end{equation}
and thus the linearised estimates for the eigenvalues are recovered
\begin{equation}
\begin{aligned}
   \lambda_1 = \sqrt{a^2 + 4\,c_\up{s}^2/3},\quad
   \lambda_2 = c_\up{s},\quad
   \lambda_3 = c_\up{s}, 
\end{aligned}
\end{equation}
as given and employed in \cite{GPRmodel}.
For simplicity, when adopting semi-implicit schemes, such as the one presented in this paper, 
the eigenvalues estimates used in the Rusanov dissipation and to define the timestep size, 
are taken to be the same used for explicit methods, but setting the adiabatic sound speed to $a = 0$, 
reflecting the fact that the implicit solution of the pressure subsystem eliminates the timestep restrictions
due to acoustic waves.

\subsection{Explicit discretisation of the convective subsystem}

\subsubsection{Data reconstruction and slope limiting}

In order to achieve second order spatial accuracy for the convective fluxes, a data reconstruction
yielding a piecewise first-degree polynomial representation of the state variables, denoted
$\qrec_{ij}(x,\ y)$, must be carried out.
For our PDE system, it is convenient to compute such a reconstruction in the primitive variable 
space, specified by choosing as a primitive state vector
\begin{equation}
\label{eqn.primvars}
   \vec{V} = \left(\rho_1,\ \rho_2,\ \vec{u},\ p,\ \alpha_1,\ \vec{b},\ \vec{A}\right)^\transpose,
\end{equation}
which is related to the conserved state $\vec{Q}$ by
\begin{equation}
   \vec{V}(x,\ y) = \mathcal{P}\left[\vec{Q}(x,\ y)\right], \text{ and }
   \vec{Q}(x,\ y) = \mathcal{C}\left[\vec{V}(x,\ y)\right].
\end{equation}
Then we denote the primitive-variable reconstruction polynomial as 
$\vrec(x,\ y) = \mathcal{P}\left[\qrec(x,\ y)\right]$, and, complementarily, 
the conserved variable reconstruction polynomial is $\qrec(x,\ y) = \mathcal{C}\left[\vrec(x,\ y)\right]$.
For the sake of clarity, we remark that the primitive-to-conserved operator $\mathcal{C}$ and
conserved-to-primitive operator $\mathcal{P}$, due to their nonlinear nature, 
are to be read as \emph{pointwise} operations, and the evaluation points will be explicitly stated in the
following whenever conversion of state variables is necessary.

For each cell of index $i$, the left and right differences are computed in 
the primitive variable space as $\Delta\vec{V}_L = \vec{V}_{i} - \vec{V}_{i-1}$ 
and $\Delta\vec{V}_R = \vec{V}_{i+1} - \vec{V}_{i}$ respectively. They are then combined in a nonlinear
fashion to ensure non-oscillatory property of the resulting scheme.
In particular, we employ a simple slope limiter that can be computed as 
\begin{equation}
\label{eqn.limiter}
\begin{aligned}
   \widetilde{\Delta\vec{V}}_i = 
   &\dfrac{\Delta\vec{V}_R\,\max\left[0,\ \min\left(\beta\,\Delta\vec{V}_R^2,\ \Delta\vec{V}_R\,\Delta\vec{V}_L\right)\right]}{2\,\Delta\vec{V}_R^2 + \epsilon^2} + \\
   &\dfrac{\Delta\vec{V}_L\,\max\left[0,\ \min\left(\beta\,\Delta\vec{V}_L^2,\ \Delta\vec{V}_L\,\Delta\vec{V}_R\right)\right]}{2\,\Delta\vec{V}_L^2 + \epsilon^2},
\end{aligned}
\end{equation}
where $\epsilon = 10^{-14}$ is a small constant that avoids division by zero and all operations are to be taken componentwise.


 
The slope limiter \eqref{eqn.limiter} yields the minmod slope for $\beta = 1$ and reduces to the MUSCL--Barth--Jespersen limiter
for $\beta = 3$. In all our numerical tests we set $\beta = 2$.

The preliminary (undivided) slope $\widetilde{\Delta\vec{V}}_i$ is then corrected to enforce the
positivity of the reconstructed values of density and pressure, 
as well as the unit-sum constraints on 
the volume fractions $\alpha_1$ and $\alpha_2$. This is achieved by rescaling the slope with
\begin{equation}
   \Delta\vec{V}_i = \devi\,\min\left(1,\ \phi_i^+,\ \phi_i^-\right),
\end{equation}
having set 
\begin{equation}
\begin{aligned}
    &\phi_i^+ = \dfrac{\left[
      \left(|\devi| + \devi\right)\,\left(\vec{V}_\up{max} - \vec{V}_i\right) + 
      \left(|\devi| - \devi\right)\,\left(\vec{V}_\up{min} - \vec{V}_i\right)\right]\,\devi}{2\,|\devic| + \epsilon^3},\\
   &\phi_i^- = \dfrac{\left[
      \left(|\devi| - \devi\right)\,\left(\vec{V}_i - \vec{V}_\up{max}\right) + 
      \left(|\devi| + \devi\right)\,\left(\vec{V}_i - \vec{V}_\up{min}\right)\right]\,\devi}{2\,|\devic| + \epsilon^3},
\end{aligned}
\end{equation}
where, with reference to the primitive variable state vector \eqref{eqn.primvars},
we have set the lower and upper bounds for each variable as 
\begin{equation}
\begin{aligned}
    &\vec{V}_\up{min} = \left(0,\ 0,\ -\vec{h},\ 0,\ 0,\ -\vec{h},\ -\vec{H}\right)^\transpose,\\[1mm]
    &\vec{V}_\up{max} = \left(H,\ H,\ \vec{h},\ H,\ 1,\ \vec{h},\ \vec{H}\right)^\transpose.
\end{aligned}
\end{equation}
The values of $H$, $\vec{h}$, and $\vec{H}$ are set to a 
large arbitrary scalar, vector or matrix (like $H = 10^{40}$) 
to represent the absence of
an upper or lower bound for the corresponding variable.
The same sequence of operations is carried out in the $y$-direction to compute $\Delta\vec{V}_j$.
Then the primitive reconstruction polynomial can be evaluated at any point in space as 
\begin{equation}
\vrec_{ij}(x,\ y) = \vec{V}_{ij} + 
    (x - x_{ij})\,\dfrac{\Delta\vec{V}_i}{\Delta x} + 
    (y - y_{ij})\,\dfrac{\Delta\vec{V}_j}{\Delta y}.
\end{equation}

\subsubsection{Computation of convective fluxes}
The convective terms are explicitly integrated by means of a path-conservative MUSCL--Hancock scheme.
The fully discrete one-step update formula reads
\begin{equation}
   \begin{aligned}
   \vec{Q}_{ij}^{n+1} = \vec{Q}_{ij}^n &- \dfrac{\Delta t}{\Delta x}\,
       \left(\vec{F}^{\up{\textsc{rs}}}_{i+1/2,\,j} - \vec{F}^{\up{\textsc{rs}}}_{i-1/2,\,j} + 
       \vec{D}_{i-1/2,\,j}^{+} + \vec{D}_{i+1/2,\,j}^-\right) +  \\
       &-\dfrac{\Delta t}{\Delta y}\,\left(\vec{F}^{\up{\textsc{rs}}}_{i,\,j+1/2} - \vec{F}^{\up{\textsc{rs}}}_{i,\,j-1/2} + 
       \vec{D}_{i,\,j-1/2}^{+} + \vec{D}_{i,\,j+1/2}^- \right) + \\
       &-\dfrac{\Delta t}{\Delta x}\,\vec{B}_{1}^\up{p}\left[\vec{v}_{ij}(t^{n+1/2},\ x_i,\ y_j)\right]\,\Delta\vec{V}_i +\\
       &-\dfrac{\Delta t}{\Delta y}\,\vec{B}_{2}^\up{p}\left[\vec{v}_{ij}(t^{n+1/2},\ x_i,\ y_j)\right]\,\Delta\vec{V}_j + \\
       &+\Delta t \,\vec{S}\left[\vec{v}_{ij}\left(t^{n+1/2},\ x_i,\ y_j\right)\right], 
   \end{aligned}
\end{equation}
which is then applied to the convective subsystem and used to formally define a convective state
$\vec{Q}_{ij}^\ast = \vec{Q}_{ij}^{n+1}$, which in particular is 
\begin{equation}
\label{eqn.Qstar}
         \vec{Q}_{ij}^\ast = \left((\alpha_1\,\rho_1)_{ij}^\ast,\ (\alpha_2\,\rho_2)_{ij}^\ast,\ 
         (\rho\,\vec{u})_{ij}^\ast,\ (\rho\,E)_{ij}^\ast,\ (\alpha_1)_{ij}^\ast,\ (\vec{b})_{ij}^\ast,\ (\vec{A})_{ij}^\ast\right)^\transpose.
 \end{equation} 
For the computation of the conservative numerical fluxes we employ the simple Rusanov Riemann solver 
\begin{equation}
\begin{aligned}
   &\vec{F}^{\up{\textsc{rs}}}_{i+1/2,\,j}(\vec{v}_L,\ \vec{v}_R) = \dfrac{1}{2}\,\left[\fll\vec{F}_1(\vec{v}_L) + \vec{F}_1(\vec{v}_R)\right]
    - \dfrac{1}{2}\,s_1^\up{max}\,\left[\fll\mathcal{C}(\vec{v}_R) - \mathcal{C}(\vec{v}_L)\right],\\
   &\vec{F}^{\up{\textsc{rs}}}_{i,\,j+1/2}(\vec{v}_L,\ \vec{v}_R) = \dfrac{1}{2}\,\left[\fll\vec{F}_2(\vec{v}_L) + \vec{F}_2(\vec{v}_R)\right]
    - \dfrac{1}{2}\,s_2^\up{max}\,\left[\fll\mathcal{C}(\vec{v}_R) - \mathcal{C}(\vec{v}_L)\right],
\end{aligned}
\end{equation}
where the signal speed estimates $s_1^\up{max}$ and $s_2^\up{max}$, i.e. 
the maximum absolute eigenvalues of the convective subsystem, in the first or second space direction, 
associated with a pair of states $\vec{v}_L$ and $\vec{v}_R$, are
given by
\begin{equation}
\begin{aligned}
    &s_1^\up{max} = s_1^\up{max}\left(\fls\vec{v}_L,\ \vec{v}_R\right) = 
        \max\left[\fll
        \max\left(\fls|\vec{\lambda}_1\left(\vec{v}_L\right)|\right),\ 
        \max\left(\fls|\vec{\lambda}_1\left(\vec{v}_R\right)|\right)\right],\\
    &s_2^\up{max} = s_2^\up{max}\left(\fls\vec{v}_L,\ \vec{v}_R\right) = 
        \max\left[\fll
        \max\left(\fls|\vec{\lambda}_2\left(\vec{v}_L\right)|\right),\ 
        \max\left(\fls|\vec{\lambda}_2\left(\vec{v}_R\right)|\right)\right].
\end{aligned}
\end{equation}
An important consideration is that, in order to achieve a compatible discretisation of 
density, momentum, and kinetic energy in uniform flows, the jump of conserved variables 
$\mathcal{C}(\vec{v}_R) - \mathcal{C}(\vec{v}_L)$ must \emph{exclude} the difference of internal energies.
Thus, instead of ${\Delta E = \rho\,E(\vec{v}_R) - \rho\,E(\vec{v}_L)}$, the jump in the 
energy conservation equation will be 
\begin{equation}
\begin{aligned}
   \Delta E^\prime(\vec{v}_L,\ \vec{v}_R) &= \left[\rho\,E\lvert_{\vec{v}_R} - 
   \rho\,e\lvert_{\vec{v}_R}\fll\right] - \left[\rho\,E\lvert_{\vec{v}_L} - \rho\,e\lvert_{\vec{v}_R}\fll\right] = \\
          &= \rho\,(e_\up{k} + e_\up{s} + e_\up{t})\lvert_{\vec{v}_R} - \rho\,(e_\up{k} + e_\up{s} + e_\up{t})\lvert_{\vec{v}_L},
\end{aligned}
\end{equation}
where we recall $e_\up{k} = \norm{\vec{u}}^2/2$, $e_\up{s} = c_\up{s}^2\,\tr{\left(\dev{\vec{G}}\,\dev{\vec{G}}\right)}/4$,
$e_\up{t} = \sigma\,\norm{\vec{b}}/\rho$.

Nonconservative products are discretised within the path-conservative formalism \cite{castro2006, pares2006}.
This means that, at each cell interface, generically indexed with ${i+1/2,j}$, in addition 
to numerical fluxes, two so-called fluctuations, denoted
${\vec{D}_{i+1/2,j}^-\left(\vec{v}_L,\ \vec{v}_R\right)}$ and 
${\vec{D}_{i+1/2,j}^+\left(\vec{v}_L,\ \vec{v}_R\right)}$, have to be computed.
Our choice of for the discrete jump terms (fluctuations) is such that, 
the left and right fluctuations have the same value
and may be denoted as $\vec{D}_{i+1/2,j}\left(\vec{v}_L,\ \vec{v}_R\right)$.
The same holds in the $y$-direction for $\vec{D}_{i,j+1/2}^-\left(\vec{v}_L,\ \vec{v}_R\right)$ and 
$\vec{D}_{i,j+1/2}^+\left(\vec{v}_L,\ \vec{v}_R\right)$.
The fluctuations are computed with a three-point point quadrature rule over a segment path 
$\vec{\Psi}(\vec{v}_L,\ \vec{v}_R,\ s) = (1 - s)\,\vec{v}_L + s\,\vec{v}_R$ in the primitive state space. 
Their discrete expressions read
\begin{equation}
\begin{aligned}
    &\vec{D}_{i+1/2,j}\left(\vec{v}_L,\ \vec{v}_R\right) =
    \dfrac{1}{2}\,\sum_{k = 1}^3 \omega_k\,\vec{B}_1^\up{p}\left[\fll\vec{\Psi}\left(\vec{v}_{L},\ \vec{v}_{R},\ s_k\right)\right]\,\left(\vec{v}_R - \vec{v}_L\right),\\
    &\vec{D}_{i,j+1/2}\left(\vec{v}_L,\ \vec{v}_R\right) =
    \dfrac{1}{2}\,\sum_{k = 1}^3 \omega_k\,\vec{B}_2^\up{p}\left[\fll\vec{\Psi}\left(\vec{v}_{L},\ \vec{v}_{R},\ s_k\right)\right]\,\left(\vec{v}_R - \vec{v}_L\right),
\end{aligned}
\end{equation}
where the nonconservative products for the convective subsystem are given by
\begin{equation}
\begin{aligned}
&\vec{B}_{1}^\up{p}(\vec{v})\,\Delta\vec{V}_i = \left(0,\ 0,\ \vec{0},\ 0,\ u_1(\vec{v})\,(\Delta \alpha_1)_i - K(\vec{v})\,(\Delta u_1)_i,\ \vec{0},\ \vec{0}\right)^\transpose,\\
&\vec{B}_{2}^\up{p}(\vec{v})\,\Delta\vec{V}_j = \left(0,\ 0,\ \vec{0},\ 0,\ u_2(\vec{v})\,(\Delta \alpha_1)_j - K(\vec{v})\,(\Delta u_2)_j,\ \vec{0},\ \vec{0}\right)^\transpose,
\end{aligned}
\end{equation}
with $K = \left(\rho_2\,a_2^2 - \rho_1\,a_1^2\right)\,\alpha_1\,\alpha_2/(\alpha_1\,\rho_2\,a_2^2 + \alpha_2\,\rho_1\,a_1^2)$.
For clarity, we explicitly give also the expressions for the fluxes of the convective subsystem
\begin{equation} 
\label{eqn.split.def2} 
\mathbf{F}_1 = \left( \begin{array}{c} \alpha_1\,\rho_1\,u_1 \\ \alpha_2\,\rho_2\,u_1  \\ \rho\,\vec{u}\,u_1 \\ \left(\rho\,E - \rho\,e\right)\,u_1\\ 0 \\ \vec{0} \\ \vec{0} \end{array} \right),  \quad 
\mathbf{F}_2 = \left( \begin{array}{c} \alpha_1\,\rho_1\,u_2 \\ \alpha_2\,\rho_2\,u_2  \\ \rho\,\vec{u}\,u_2 \\ \left(\rho\,E - \rho\,e\right)\,u_2 \\ 0 \\ \vec{0} \\ \vec{0} \end{array} \right),  \quad 
\end{equation}
while the source term is simply $\vec{S} = \vec{0}$.


At each cell interface, of generic index $i+\frac{1}{2}, j$ in the $x$-direction or $i, j+\frac{1}{2}$ in the $y$-direction, 
the boundary-extrapolated states $\vec{v}_L$ and $\vec{v}_R$ are taken from a cell-local space-time predictor 
solution $\vec{v}_{ij}(t,\ x,\ y)$. In particular, the space-time midpoint values for each face are
\begin{equation}
   \begin{aligned} 
&(\vec{v}_L)_{i+1/2,j} = \vec{v}_{i,j}(t^{n+1/2},\ x_{i+1/2},\ y_{j}),\\
&(\vec{v}_R)_{i+1/2,j} = \vec{v}_{i+1,j}(t^{n+1/2},\ x_{i+1/2},\ y_{j}),\\
&(\vec{v}_L)_{i,j+1/2} = \vec{v}_{i,j}(t^{n+1/2},\ x_{i},\ y_{j+1/2}),\\
&(\vec{v}_R)_{i,j+1/2} = \vec{v}_{i,j+1}(t^{n+1/2},\ x_{i},\ y_{j+1/2}), 
   \end{aligned}
\end{equation}
and they are explicitly computed as
\begin{equation}\label{eq.recp2c}
   \begin{aligned}
       \vec{v}_{ij}(t^{n+1/2},\ x,\ y)  & = \mathcal{P}\left[\vec{q}_{ij}(t^{n+1/2},\ x,\ y)  \right] \\
                                        & = \mathcal{P}\left\{\mathcal{C}\left[\vrec_{ij}(x,\ y)\right] + \Delta \vec{Q}_{ij}\right\},
   \end{aligned}
\end{equation}
where
\begin{equation}
\begin{aligned}
   \Delta\vec{Q}_{ij} = -&\frac{\Delta t}{2\,\Delta x}\,\left\{\flll\vec{F}_1\left[\vrec_{ij}(x_{i+1/2},\ y_{j})\right] - \vec{F}_1\left[\vrec_{ij}(x_{i-1/2},\ y_{j})\right]\right\} +\\
   - &\frac{\Delta t}{2\,\Delta y}\,\left\{\flll\vec{F}_2\left[\vrec_{ij}(x_i,\ y_{j+1/2})\right] - \vec{F}_2\left[\vrec_{ij}(x_i,\ y_{j-1/2})\right]\right\} +\\
   - &\frac{\Delta t}{2\,\Delta x}\,\vec{B}^{\up{p}}_1\left[\vrec_{ij}(x_{i},\ y_{j})\right]\,\Delta\vec{V}_i
   - \frac{\Delta t}{2\,\Delta y}\,\vec{B}^{\up{p}}_2\left[\vrec_{ij}(x_{i},\ y_{j})\right]\,\Delta\vec{V}_j + \\
   + &\frac{\Delta t}{2}\,\vec{S}\left[\vrec_{ij}(x_{i},\ y_{j})\right].\\
\end{aligned}
\end{equation}

For the sake of clarity, it should be pointed out that the primitive-to-conserved and conserved-to-primitive conversion operators in Eq. \eqref{eq.recp2c} are to be read
as \emph{pointwise} operations, or equivalently the formula can be taken as a projection between two different
polynomial spaces, one in which the conserved variables are polynomials but the primitive ones are not, and viceversa, 
but it is not \emph{strictly} satisfied in any point except those where the conversion of state 
variables has taken place, i.e. the space-time barycenters of each face.
paspso

\subsection{Staggered mesh and discrete divergence, curl and gradient operators}
\begin{figure}[!bp]
   \begin{center}
      \includegraphics[width=\textwidth]{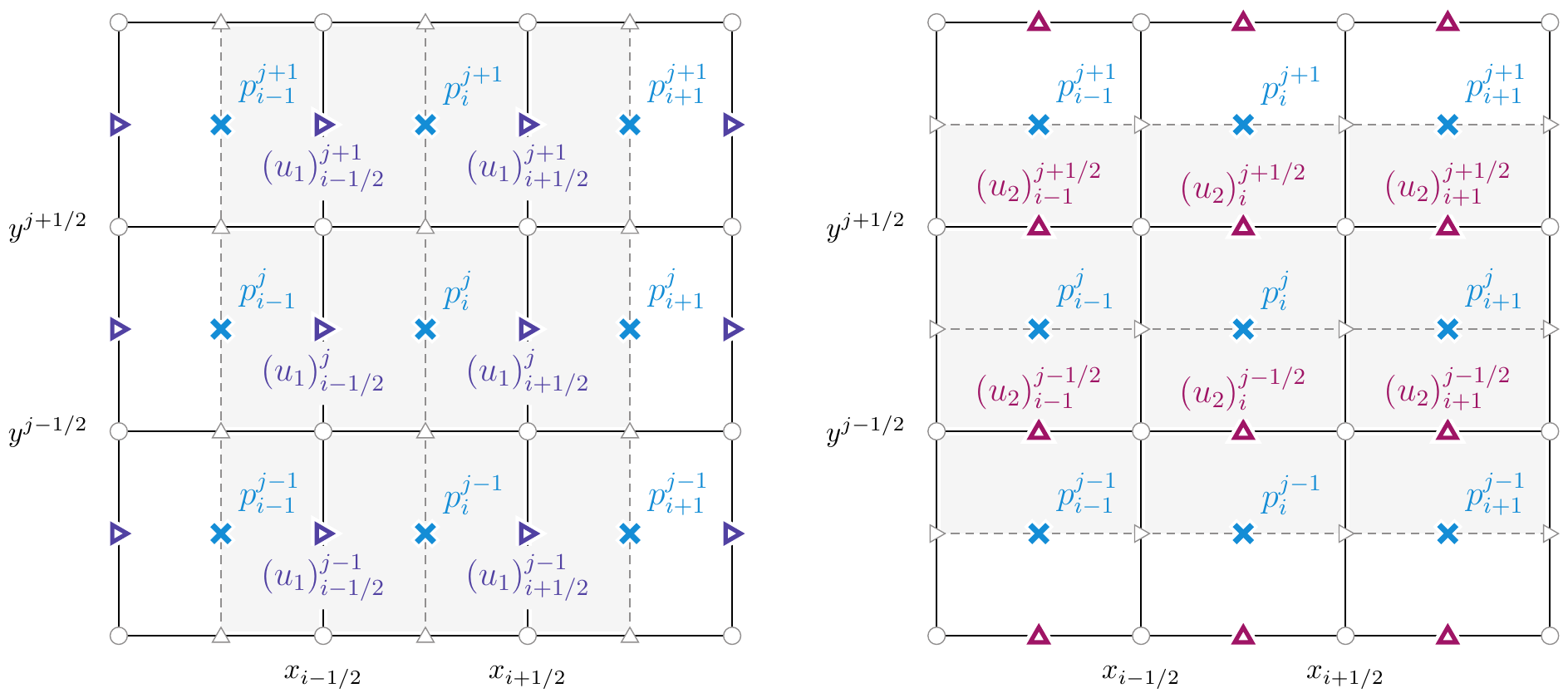}  
      \caption{Staggered mesh configuration with the pressure field $p_{i}^{j}$ defined in the 
       cell barycenters and the velocity field components $(u_1)_{i+1/2}^{j}$ 
       and $(u_2)_{i}^{j+1/2}$ defined 
       on the edge-based staggered dual grids.}  
      \label{fig.staggerededge}
   \end{center}
\end{figure}
%
The numerical scheme is presented in a two-dimensional context, however it is necessary and beneficial to retain
all components of three-dimensional vectors to simplify the treatment of the relaxation source term, which
acts on all components of the distortion matrix even regardless of whether derivatives in any direction vanish or not.
Again, we consider a 
two-dimensional physical domain $ \Omega $
covered by a set of uniformly sized and non-overlapping 
Cartesian control volumes $\Omega_{ij} = [x_{i-1/2},\ x_{i+1/2}] \times 
[y_{j-1/2},\ y_{j+1/2}]$ with 
mesh spacings $\Delta x = x_{i+1/2} - x_{i-1/2}$ and $\Delta y = y_{j+1/2} - 
y_{j-1/2}$ 
in $x$ and $y$ direction respectively. With $x_{i \pm 1/2}=x_i \pm \Delta x / 2$ and 
$y_{j \pm 1/2} = y_j \pm \Delta y /2$, we denote the barycenter coordinates of the control volumes. 
We will furthermore use the notation $\uvec{e}_x = 
(1,\ 0,\ 0)$, 
$\uvec{e}_y = (0,\ 1,\ 0)$ and $\uvec{e}_z = (0,\ 0,\ 1)$ for the unit vectors pointing into the 
directions 
of the Cartesian coordinate axes. 

\begin{figure}[!bp]
   \begin{center}
      \includegraphics[width=\textwidth]{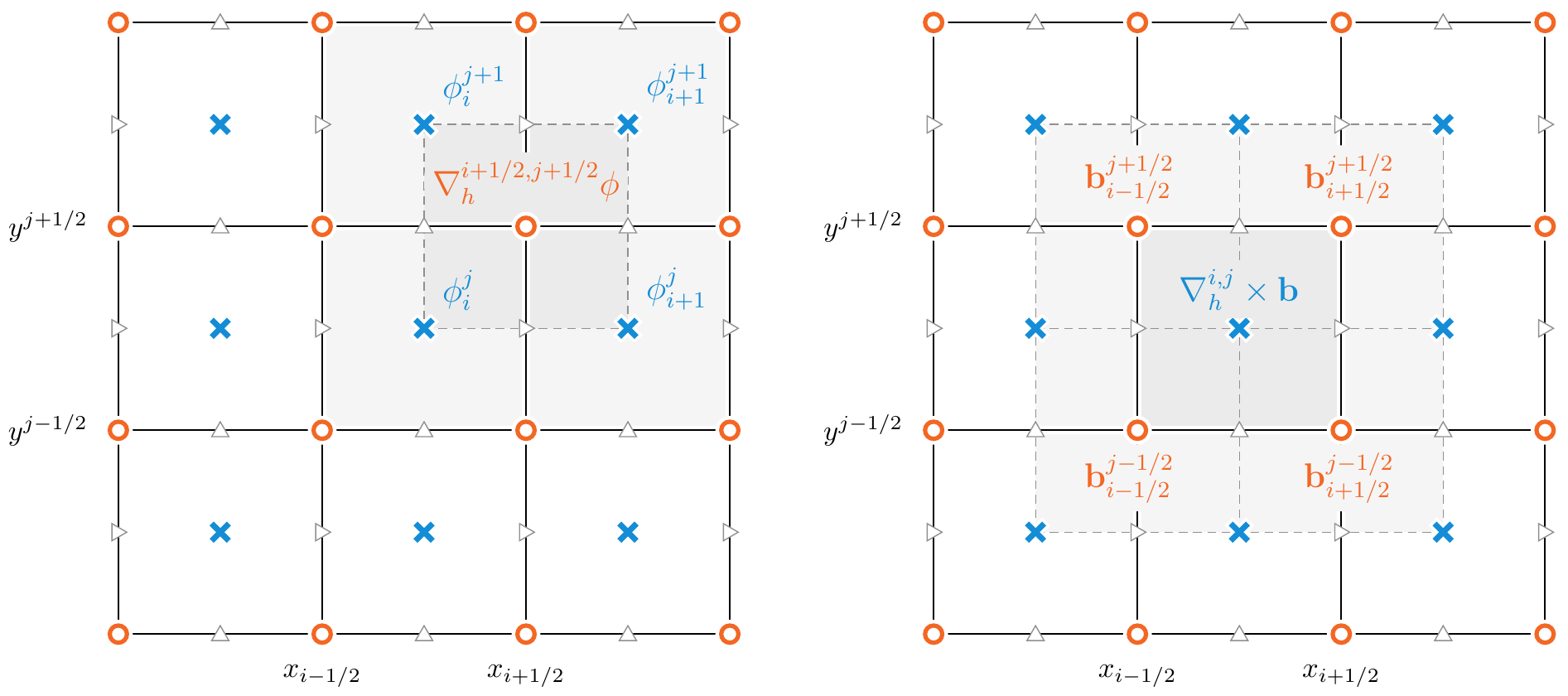}  
      \caption{Staggered mesh configuration with a scalar field field $\phi_{i}^{j}$ 
       defined in the 
       cell barycenters, and the interface 
       field field $\vec{b}_{i+1/2}^{j+1/2}$ defined on the vertices of the main grid.
       The shaded control volumes indicate the stencil for the computation of the
       corner gradients $\nabla_h^{i+1/2,j+1/2}\phi$ and for the cell-centred 
       curl operator $\nabla_h^{i,j}\times\vec{b}$.} 
      \label{fig.staggeredcorner}
   \end{center}
\end{figure}

The set of discrete 
times will be denoted by $t^n$. 
For a sketch of the employed staggered grid arrangement of the main quantities, 
see Fig. \ref{fig.staggeredcorner} and Fig. \ref{fig.staggerededge}. 

The main ingredients of the 
structure-preserving staggered 
semi-implicit scheme proposed 
in this Chapter are the definitions of
appropriate discrete divergence, gradient and curl operators acting on quantities that are 
arranged in different and judiciously 
chosen locations on the staggered mesh. The discrete pressure field at time
$t^n$ is denoted by $p_{ij}$ and its degrees of freedom are located in the center of each control 
volume as $p_{ij}^n=p(t^n, x_i,\ y_j)$. 

The discrete velocities $u_1^{n}$ and $u_2^{n}$ are arranged in an edge-based staggered fashion, i.e.  
$(u^n_{1})_{i+1/2}^{j} = u_1(t^n,\ x_{i+1/2},\ y_j)$ and 
$(u^n_{2})_{i}^{j+1/2} = u_2(t^n,\ x_{i},\ y_{j+1/2})$. 
The discrete vector field $\vec{b}^{n}$ is defined on 
the \textit{vertices} of each spatial control volume as  
$\vec{b}^n_{i+1/2,j+1/2} = \vec{b}(t^n,\ x_{i+1/2},\ y_{j+1/2})$. 
For clarity, see again Fig. \ref{fig.staggeredcorner}.  

The \textit{discrete divergence operator}, $\dgrad \cdot$, acting on a discrete vector field 
$\vec{u}^{n}$ is abbreviated 
by $\dgrad \cdot \vec{u}^{n}$ and its degrees of freedom are given by  
\begin{equation}
\dgrad^{i,j} \cdot \vec{u}^{n}  =  \left(\dgrad \cdot \vec{u}^{n}\right)_{i,j}  =
\frac{(u_1^n)_{i + 1/2}^{j} - (u_1^n)_{i - 1/2}^{j}}{\Delta x} + 
\frac{(u_2^n)_{i}^{j+1/2} - (u_2^n)_{i}^{j-1/2}}{\Delta y}, 
\label{eqn.div} 
\end{equation}
i.e. it is based on the \textit{edge-based} staggered values of the field $\vec{u}^{n}$. It defines  
a discrete divergence on the control volume $\Omega_{ij}$ via the Gauss theorem, 
\begin{equation}
\dgrad \cdot \vec{u}^{n} = \frac{1}{\Delta x\,\Delta y} \int \limits_{\Omega_{ij}} \nabla \cdot \vec{u} \,\de{\vec{x}} 
= \frac{1}{\Delta x\,\Delta y} \int \limits_{\partial \Omega_{ij}} \vec{u} \cdot \uvec{n} \, \de{S},  
\label{eqn.gauss}
\end{equation}
based on the mid-point rule for the computation of the integrals along each edge of $\Omega_{ij}$. In \eqref{eqn.gauss} the
outward pointing unit normal vector to the boundary $\partial \Omega_{ij}$ of $\Omega_{ij}$ is denoted by $\uvec{n}$. 

In a similar manner, the $z$ component of the \textit{discrete curl}, $\dgrad \times $, of a 
discrete vector field 
$\vec{b}^{n}$ (or $\vec{a}_1^n = (1,\ 0,\ 0)\,\vec{A}^n$ for example)
is denoted by $\left( \dgrad \times \mathbf{b}^{n} \right) \cdot \uvec{e}_z$ and its degrees of freedom are naturally defined as
\begin{equation}
\begin{aligned}
\left( \dgrad^{i,j} \times \mathbf{b}^{n} \right) \cdot \uvec{e}_z = 
{\left[\flll\left( \dgrad \times \mathbf{b}^{n} \right) \cdot \uvec{e}_z\right]}_{i,j}\,=&\\
 \frac{(b_2^n)_{i+1/2}^{j+1/2} - (b_2^n)_{i-1/2}^{j+1/2} + 
(b_2^n)_{i+1/2}^{j-1/2} - (b_2^n)_{i-1/2}^{j-1/2}}{2\,\Delta x} \,+& \\
\frac{(b_1)_{i+1/2}^{j+1/2} - (b_1^n)_{i+1/2}^{j-1/2} + 
(b_1^n)_{i-1/2}^{j+1/2} - (b_1^n)_{i-1/2}^{j-1/2}}{2\,\Delta y}\,\phantom{+}&
\label{eqn.rot} 
\end{aligned}
\end{equation}
making use of the \textit{vertex based} staggered values of the field $\vec{b}^{n}$, 
see Fig. \ref{fig.staggeredcorner}. 
In the present two-dimensional description of the scheme the first and second components of 
the discrete curl $\dgrad^{i,j}\times\vec{b}^n$ vanish identically.
Eqn. \eqref{eqn.rot} defines a discrete
curl on the control volume $\Omega_{ij}$ via the Stokes theorem
\begin{equation}
\left( \dgrad \times \vec{b}^{n} \right) \cdot \uvec{e}_z = \frac{1}{\Delta x\,\Delta y} \int \limits_{\Omega_{ij}} 
\left( \nabla \times  \vec{b} \right) \cdot \uvec{e}_z \, \de\mathbf{x} = \frac{1}{\Delta x\, \Delta y} \int \limits_{\partial \Omega_{ij}} \vec{b} \cdot \uvec{t} \, \de{S},  
\label{eqn.stokes}
\end{equation}
based on the trapezoidal 
rule for the computation of the integrals along each edge of $\Omega_{ij}$. 

Last but not least, we need to define a discrete gradient operator that is compatible with the discrete curl,
so that the continuous identity
\begin{equation}
\nabla \times \nabla \phi = \vec{0}
\label{eqn.rotgrad} 
\end{equation}
also holds on the discrete level. If we define a scalar field in the barycenters of the control volumes $\Omega_{ij}$ as
$\phi_{ij}^{n}=\phi(t^n,\ x_i,\ y_j)$ then the corner gradient generates a natural discrete gradient operator $\dgrad$ 
of the discrete scalar field $\phi^{n}$ that defines a discrete gradient in all vertices of the mesh. 
The corresponding degrees of freedom generated by $\dgrad \phi^{n}$ read (see 
Fig. \ref{fig.staggeredcorner})
\begin{equation}
\label{discr.grad}
\dgrad^{i+1/2,j+1/2}  \phi^{n} = \left(\dgrad  \phi^{n}\right)_{i+1/2}^{j+1/2} =  \left( 
\begin{array}{c}  
  \dfrac{\phi^n_{i + 1, j + 1} - \phi^n_{i, j + 1}
  + \phi^n_{i + 1, j}     - \phi^n_{i, j} }{2\,\Delta x}  \\[4mm]
  \dfrac{\phi^n_{i + 1, j + 1} - \phi^n_{i + 1, j} + 
  \phi^n_{i, j + 1} - \phi^n_{i,j} }{2\,\Delta y}   \\[3mm] 
 0\vphantom{\dfrac{0^n}{\Delta y}}
\end{array} \right).
\end{equation}

It is then straightforward to verify that an immediate consequence of Equations \eqref{eqn.rot} and 
\eqref{discr.grad} is 
\begin{equation}
\dgrad \times \dgrad \phi^{n} = \vec{0}, 
\label{eqn.discrete.curl} 
\end{equation}
i.e. one obtains a discrete analogue of \eqref{eqn.rotgrad}. This can be easily seen by computing 
\begin{equation}
\label{eqn.rot.grad} 
\begin{aligned}
\left( \dgrad^{i,j} \times \nabla^{i+1/2,j+1/2}  \phi^{n} \right) \cdot \uvec{e}_z  & =   \\[2mm] 
\frac{\left( \phi_{i + 1}^{ j + 1} - \phi_{i + 1}^{ j}   + \phi_{i}^{ j + 1}   - \phi_{i}^{j} \right)    + 
      \left( \phi_{i + 1}^{ j }    - \phi_{i + 1}^{ j-1} + \phi_{i}^{ j }      - \phi_{i}^{j-1} \right)  }{4\,\Delta x\, \Delta y\fll}\,  &- \\[2mm]  
\frac{\left( \phi_{i }^{ j + 1}    - \phi_{i }^{ j}      + \phi_{i-1}^{ j + 1} - \phi_{i-1}^{j} \right)  + 
      \left( \phi_{i }^{ j }       - \phi_{i }^{ j-1}    + \phi_{i-1}^{ j }    - \phi_{i-1}^{j-1} \right) }{4\,\Delta x\, \Delta y\fll}\,  &- \\[2mm] 
\frac{\left( \phi_{i + 1}^{ j + 1} - \phi_{i}^{ j + 1}   + \phi_{i + 1}^{ j}   - \phi_{i}^{j} \right)   + 
      \left( \phi_{i }^{ j + 1}    - \phi_{i-1}^{ j + 1} + \phi_{i }^{ j}      - \phi_{i-1}^{j} \right)    }{4\,\Delta y\, \Delta x\fll}\,  &+ \\[2mm]  
\frac{\left( \phi_{i + 1}^{ j }    - \phi_{i}^{ j }      + \phi_{i + 1}^{ j-1} - \phi_{i}^{j-1} \right) + 
      \left( \phi_{i }^{ j }       - \phi_{i-1}^{ j }    + \phi_{i }^{ j-1}    - \phi_{i-1}^{j-1} \right)    }{4\,\Delta y\, \Delta x\fll}\, &= 0.
\end{aligned}
\end{equation}
We furthermore define the following averaging operators from the edge-based staggered meshes to the cell barycenter $(x_i,\ y_j)$ and viceversa 
\begin{equation}
\begin{aligned}
(u_1^n)_{i,j}     =& \frac{1}{2} \left[ \left(u_1^n\right)_{i-1/2,j} + \left(u_1^n\right)_{i+1/2,j}\fll \right],\\                                                                   
(u_2^n)_{i,j}     =& \frac{1}{2} \left[ \left(u_2^n\right)_{i,j-1/2} + \left(u_2^n\right)_{i,j+1/2}\fll \right],\\     
(u_1^n)_{i+1/2,j} =& \frac{1}{2} \left[ \left(u_1^n\right)_{i,j} + \left(u_1^n\right)_{i+1,j}\fll \right],\\                                                           
(u_2^n)_{i,j+1/2} =& \frac{1}{2} \left[ \left(u_2^n\right)_{i,j} + \left(u_2^n\right)_{i,j+1}\fll \right].                                                                   
\end{aligned}
\end{equation}
Interpolation from the cell centers to the the cell-edge locations is required for the construction of a compact stencil symmetric positive definite
system for the pressure wave equation, while interpolation to the cell centers can be avoided by using face-interpolated pressure values for a direct update
of the momentum on the barycenter grid.

Finally we also introduce an interpolation operator to compute cell center approximations of the corner quantities like $\vec{b}$ and $\vec{a}_{sk}$
which we write
\begin{equation}
   \vec{b}_{i,j}^n = \dfrac{1}{4}\,\left(\vec{b}^n_{i-1/2,j-1/2}+\vec{b}^n_{i+1/2,j-1/2}+\vec{b}^n_{i-1/2,j+1/2}+\vec{b}^n_{i-1/2,j+1/2}\right)
\end{equation}
which represents a simple (linear) arithmetic averaging operator and introduces minimal numerical dissipation. 
Whenever flow convection is particularly strong, we found beneficial to
apply a partial upwinding to the interpolation operator for the interface field $\vec{b}$, 
so to add additional numerical stabilisiation to the scheme.
\begin{equation}
   \vec{b}_{i,j}^n = w_1\,(\vec{b}^n)_{i-1/2}^{j-1/2}+w_2\,(\vec{b}^n)_{i+1/2}^{j-1/2}+w_3\,(\vec{b}^n)_{i-1/2}^{j+1/2}+w_4\,(\vec{b}^n)_{i-1/2}^{j+1/2}.
\end{equation}
The coefficients $w_k$ are obtained by first constructing a set of preliminary weighs by two-dimensional upwinding.
\begin{equation}
   \begin{aligned}
      &w_1^\ast = \epsilon + u_1^+ + u_2^+,\qquad
       w_2^\ast = \epsilon + u_1^- + u_2^+,\\
      &w_3^\ast = \epsilon + u_1^+ + u_2^-,\qquad
       w_4^\ast = \epsilon + u_1^- + u_2^-,\\
   \end{aligned}
   \label{eq.uwindweights}
\end{equation}
with $\epsilon = 10^{-6}$, 
    $\vec{u}^+ = (u_1^+,\ u_2^+)^\transpose = \max\left(0,\ \vec{u}\right)/(\norm{\vec{u}} + \epsilon)$, 
and $\vec{u}^- = (u_1^-,\ u_2^-)^\transpose = \max\left(0,\ -\vec{u}\right)/(\norm{\vec{u}} + \epsilon)$.
Then the preliminary weights \eqref{eq.uwindweights} are normalized in such a way that
the upwind bias will be reduced for flows with weak convection.
The final weights thus are computed as
\begin{equation}
   w_k = \frac{w_k^\ast}{\sum_{k=1}^{4}{w_k^\ast}}\,\lambda + \frac{1 - \lambda}{4}, 
   \text{ with } \lambda = \min\left(1,\ 2\,\norm{\vec{u}}\,\sqrt{\Delta x\, \Delta y}\right).
\end{equation}


\subsection{Explicit discretization of involution constrained fields} 

The key ingredient of the numerical method proposed in this Chapter is the curl-compatible discretization of the terms
$\nabla \vec{G}_v(\vec{Q})$ and $\vec{B}_v(\vec{Q}) \nabla \vec{Q}$ present in \eqref{eqn.pde.split}. 
We propose the following compatible discretization for the interface field equation: 
\begin{equation}
\label{eqn.Jh}
\begin{aligned}
   \vec{b}^{n+1}_{i+1/2,j+1/2} = \vec{b}^n_{i+1/2,j+1/2} &- \Delta t\,\dgrad^{i+1/2,j+1/2}\left(\phi^n\right)+ \\
   & - \Delta t\,\left[\flll\left(\dgrad\times\vec{b}\right)\times\vec{u}\right]_{i+1/2,j+1/2}
\end{aligned}
\end{equation}
with the corner-averaged curl term is given by 
\begin{equation}
   \left[\flll\left(\dgrad\times\vec{b}\right)\times\vec{u}\right]_{i+1/2,j+1/2} = \dfrac{1}{4}\sum_{r=0}^1\sum_{s=0}^1\left(\dgrad^{i+r,j+s}\times\vec{b}^n\right)\times\vec{u}_{i+r,j+s}^n
\end{equation}

It is easy to check that for an initially curl-free vector field $\vec{b}^{n}$ 
which satisfies $\dgrad \times \vec{b}^{n} = 0$ also $\dgrad \times \vec{b}^{n+1} = 0$ holds. 
In order to see this, one needs to apply the discrete curl operator $\dgrad \times$ to Eqn. \eqref{eqn.Jh}.  One realizes 
that the second row of \eqref{eqn.Jh}, which contains the discrete curl of $\vec{b}^{n}$ vanishes immediately, 
due to $\dgrad \times \vec{b}^{n} = 0$. 
The curl of the first term on the right hand side in the first row of Eqn. \eqref{eqn.Jh} is zero because of 
$\dgrad \times \vec{b}^{n} = 0$ and the curl of the second term is zero because of 
$\dgrad \times \left(\dgrad \phi^{n}\right) = 0$, with the auxiliary scalar field $\phi^{n} = \vec{b}^{n} \cdot \vec{u}^{n}$, 
whose degrees of freedom are computed as 
$\phi^{n}_{i,j} = \vec{b}^{n}_{i,j}\cdot \vec{u}^{n}_{i,j}$ after interpolating the velocity vector and 
the gradient field $\vec{b}$ into the barycenters of the control volumes $\Omega_{i,j}$. 
The key ingredient of our compatible  discretization for the $\vec{b}$ equation is indeed 
the use of a discrete gradient operator that is compatible with the discrete curl operator, see Eq. \eqref{eqn.rot.grad}.

\subsection{Compatible numerical viscosity} 

The previous discretizations were all \textit{central} and thus without artificial numerical viscosity. In order to add a \textit{compatible numerical viscosity} operator, we need to recall the definition of the vector Laplacian at the continuous level, which reads: 
\begin{equation}
\nabla^2 \vec{b} = \nabla \left( \nabla \cdot \vec{b} \right) - \nabla \times \left(\nabla \times \vec{b}\right)
\label{vector.laplace} 
\end{equation}
In order to define a discrete analogue of \eqref{vector.laplace} we define another discrete divergence operator, 
which naturally follows from the definition of the discrete gradient \eqref{discr.grad}
\begin{equation}
   \begin{aligned}
      \dgrad^{i+1/2,j+1/2} \cdot \vec{b}^{n} &=
 \frac{(b_1^n)_{i+1}^{j+1} - (b_1^n)_{i}^{j+1} + 
(b_1^n)_{i+1}^{j} - (b_1^n)_{i}^{j}}{2\,\Delta x} + \\[2mm]
&+\frac{(b_2)_{i+1}^{j+1} - (b_2^n)_{i+1}^{j} + 
(b_2^n)_{i}^{j+1} - (b_2^n)_{i}^{j}}{2\,\Delta y}.
   \end{aligned}
\label{vector.divc} 
\end{equation}
and yields the degrees of freedom of the divergence of $\vec{b}^n$ at the cell corner locations, starting
from the cell center interpolated values of the vector field $\vec{b}^n$.
By shifting indices by a half step in both directions, the same operator can be used to obtain cell center values
for $\dgrad\cdot\vec{b}^n$ starting from the corner values of $\vec{b}^n$. In this case, the operator is
\begin{equation}
   \begin{aligned}
      \dgrad^{i,j} \cdot \vec{b}^{n} &=
 \frac{(b_1^n)_{i+1/2}^{j+1/2} - (b_1^n)_{i-1/2}^{j+1/2} + 
(b_1^n)_{i+1/2}^{j-1/2} - (b_1^n)_{i-1/2}^{j-1/2}}{2\,\Delta x} + \\[2mm]
&+\frac{(b_2)_{i+1/2}^{j+1/2} - (b_2^n)_{i+1/2}^{j-1/2} + 
(b_2^n)_{i-1/2}^{j+1/2} - (b_2^n)_{i-1/2}^{j-1/2}}{2\,\Delta y}.
   \end{aligned}
\label{vector.divcb} 
\end{equation}
The discrete vector Laplacian then simply reads 
\begin{equation}
\begin{aligned}
(\dgrad^2)^{i+1/2,j+1/2} \mathbf{b}^{n} &= \dgrad^{i+1/2,j+1/2} \cdot \left(\dgrad^{i,j} \mathbf{b}^{n}\right) = \\
&=\dgrad^{i+1/2,j+1/2} \left( \dgrad^{i,j} \cdot \mathbf{b}^{n} \right) - \dgrad^{i+1/2,j+1/2} \times \left(\dgrad^{i,j} \times \mathbf{b}^{n}\right),   
\label{disc.vector.laplace} 
\end{aligned}
\end{equation}
i.e. it is composed of a grad-div contribution minus a curl-curl term. 
Taking \eqref{disc.vector.laplace} into account, a compatible discretization of $\mathbf{b}$ \textit{with} numerical viscosity then reads 
\begin{equation}
\begin{aligned}
   \vec{b}^{n+1}_{i+1/2,j+1/2} = \vec{b}^n_{i+1/2,j+1/2} &- \Delta t\,\dgrad^{i+1/2,j+1/2}\left(\phi + h\,c_a\,\dgrad^{i,j}\cdot\vec{b}^n\right)+ \\
   & - \Delta t\,\left[\flll\left(\dgrad\times\vec{b}\right)\times\vec{u}\right]_{i+1/2,j+1/2} + \\
   & - \Delta t\,h\,c_a\,\dgrad^{i+1/2,j+1/2}\times\left(\dgrad^{i,j}\times\vec{b}^n\right)
\end{aligned}
\label{eqn.Jh.visc} 
\end{equation}
where $h = \max( \Delta x, \Delta y)$, is a characteristic mesh spacing and $c_a$ is a characteristic velocity
related to the artificial viscosity that one would like to add to the scheme. In practice we take
$c_a = k_L\,\lambda$, with $\lambda = \max_{\Omega}\left(\norm{\vec{u}}\right)$ for the evolution of 
the distortion field $\vec{A}$ and $\lambda = \max_\Omega\left(\norm{\vec{u}} + \sigma\norm{\vec{b}}/\rho\right)$
for the evolution of the interface field $\vec{b}$. Unless otherwise specified, we take $k_L = 0.1$.
For the sake of clarity, 
the additional numerical viscosity terms have been highlighted in red. It is obvious that 
also \eqref{eqn.Jh.visc} satisfies the curl-free property $\dgrad \times \mathbf{b}^{n+1} = \vec{0}$ if $\dgrad \times \mathbf{b}^{n} = \vec{0}$.   
In order to reduce the numerical dissipation, it is possible to employ a piecewise linear reconstruction 
and insert the barycenter extrapolated values into the discrete divergence operator under the discrete gradient.  
In two space dimensions, the curl-curl term in \eqref{eqn.Jh.visc} simplifies to 
\begin{equation}
   \dgrad^{i+1/2,j+1/2}\times\left(\dgrad^{i,j}\times\vec{b}^n\right) = 
   \left(
   \begin{array}{c}  
     -\dfrac{\omega^n_{i + 1, j + 1} - \omega^n_{i+1, j}
     + \omega^n_{i , j+1}     - \omega^n_{i, j} }{2\,\Delta y}  \\[4mm]
     \dfrac{\omega^n_{i + 1, j + 1} - \omega^n_{i, j+1} + 
     \omega^n_{i+1, j } - \omega^n_{i,j} }{2\,\Delta x}   \\[3mm] 
    0\vphantom{\dfrac{0^n}{\Delta x}}
   \end{array} \right).
\end{equation}
by denoting with $\omega^n_{i,j} = \left(\dgrad^{i,j}\times\vec{b}^n\right)\cdot\uvec{e}_z$ the third component 
of the discrete curl of $\vec{b}^n$.

\subsection{Implicit solution of the pressure equation} 

Up to now, the contribution of the pressure to the momentum and to the total energy conservation laws 
has been excluded, i.e. the terms contained in the pressure fluxes $\mathbf{F}_p$. The discrete momentum 
equations including the pressure terms read 
\begin{equation}
\label{eqn.rhou2d}
   \begin{aligned}
      &{(\rho\,u_1^{n+1})}_{i+1/2}^{j} = {(\rho\,u_1^{\ast})}_{i+1/2}^{j} + \Delta t \, (f_1^\ast)_{i+1/2}^{j} - \dfrac{\Delta t}{\Delta x}\,\left(p^{n+1}_{i+1,j} - p^{n+1}_{i,j}\right) \\
      &{(\rho\,u_2^{n+1})}_{i}^{j+1/2} = {(\rho\,u_2^{\ast})}_{i}^{j+1/2} + \Delta t \, (f_2^\ast)_{i}^{j+1/2} - \dfrac{\Delta t}{\Delta y}\,\left(p^{n+1}_{i,j+1} - p^{n+1}_{i,j}\right) 
   \end{aligned}
\end{equation}
where pressure is taken \textit{implicitly}, while nonlinear convective terms have already been discretized \textit{explicitly} 
via the operators ${(\rho\,u_1^{\ast})}_{i+1/2}^{j} $ and 
${(\rho\,u_2^{\ast})}_{i}^{j+1/2}$ given in \eqref{eqn.Qstar} and after averaging of the obtained 
quantities back to the 
edge-based staggered dual grid. 
The contribution to momentum due to the gravity source and the vertex fluxes (due to capillarity and viscosity), 
is computed, using the four point discrete divergence of the fluxes \eqref{vector.divcb}, as
\begin{equation}
   \begin{aligned}
      (\vec{f}^\ast)_{i}^j &=
 \frac{(\vec{\Omega}_{1k}^n)_{i+1/2}^{j+1/2} - (\vec{\Omega}_{1k}^n)_{i-1/2}^{j+1/2} + 
(\vec{\Omega}_{1k}^n)_{i+1/2}^{j-1/2} - (\vec{\Omega}_{1k}^n)_{i-1/2}^{j-1/2}}{2\,\Delta x} + \\[2mm]
&+\frac{(\vec{\Omega}_{2k})_{i+1/2}^{j+1/2} - (\vec{\Omega}_{2k}^n)_{i+1/2}^{j-1/2} + 
(\vec{\Omega}_{2k}^n)_{i-1/2}^{j+1/2} - (\vec{\Omega}_{2k}^n)_{i-1/2}^{j-1/2}}{2\,\Delta y} + \rho_{i,j}^{n+1}\,\vec{g}.
   \end{aligned}
\end{equation}
where $\vec{\Omega}_{1k}$ and $\vec{\Omega}_{2k}$ indicate the first
and the second row of the tensor $\vec{\Omega} = -\vec{\Sigma}_\up{t}(\vec{b}^{n+1}) - 
\vec{\Sigma}_\up{s}(\vec{A}^{n+1},\ \rho^{n+1})$  
collecting the effects of the stress tensors associated with corner quantities $\vec{b}$ and $\vec{A}$, i.e. 
capillarity and viscous forces respectively.
Both components of the flux divergence are then interpolated onto the corresponding cell edges and 
yielding
\begin{equation}
   {(f_1^\ast)}_{i+1/2}^{j} = \frac{{(\vec{f}^\ast)}_{i}^{j} + {(\vec{f}^\ast)}_{i+1}^{j}}{2}\cdot\uvec{e}_1, \quad
   {(f_2^\ast)}_{i}^{j+1/2} = \frac{{(\vec{f}^\ast)}_{i}^{j} + {(\vec{f}^\ast)}_{i}^{j+1}}{2}\cdot\uvec{e}_2
\label{eq.drhou}
\end{equation}
The first component of $\vec{f}^\ast$, ${(f_1^\ast)}_{i+1/2}^{j}$, will contribute to the momentum balance in the $x$-direction, 
and for this reason it is interpolated only at the $u_1$-velocity locations, while, 
the second component ${(f_2^\ast)}_{i}^{j+1/2}$ is part of the momentum balance in the $y$-direction and is interpolated
at the $u_2$-velocity locations.

   

A preliminary form of the discrete total energy equation reads 
\begin{equation}
\label{eqn.rhoE2d.prelim} 
 \begin{aligned}
      \rho\,e(p^{n+1}_{i,j}) + (\rho\,e_\up{t}^{n+1})_{i,j} + (\rho\,e_\up{s}^{n+1})_{i,j} + (\rho\,\tilde{e}_\up{k}^{n+1})_{i,j} = \rho\,E^\ast_{i,j} & +\\
      -\dfrac{\Delta t}{\Delta x}\,\left[
      \tilde{h}^{n+1}_{i+1/2,j}\,(\rho\,u_1^{n+1})_{i+1/2}^j - 
      \tilde{h}^{n+1}_{i-1/2,j}\,(\rho\,u_1^{n+1})_{i-1/2}^j                                  
      \right] & + \\ 
      -\dfrac{\Delta t}{\Delta y}\,\left[
      \tilde{h}^{n+1}_{i,j+1/2}\,(\rho\,u_2^{n+1})_{i}^{j+1/2} - 
      \tilde{h}^{n+1}_{i,j-1/2}\,(\rho\,u_2^{n+1})_{i}^{j-1/2}                                  
      \right] & + (\rho\,\tilde{w}_\up{g}^{n+1})_{i,j}
   \end{aligned}  
\end{equation}
with the term $(\rho\,\tilde{w}_\up{g}^{n+1})_{i,j} = \rho\,\vec{u}_{i,j}^{n+1}\cdot\vec{g}$ accounting for the work due to gravity forces.

Inserting the discrete momentum equations \eqref{eqn.rhou2d} into 
the discrete energy equation \eqref{eqn.rhoE2d.prelim} and making tilde symbols explicit via a simple Picard iteration 
(using the lower index $r$ in the following), as suggested in \cite{DumbserCasulli2016,SIMHD}, leads to the following  discrete wave equation for the unknown pressure:  
\begin{equation}
\label{eqn.p2d} 
   \begin{aligned}
      \rho\, e\left( p^{n+1}_{i,j}\right) 
      -\dfrac{\Delta t^2}{\Delta x^2}
      \left(\tilde{h}_{i+1/2,j}^{n+1}\right)_r\,\left(p^{n+1}_{i+1,j} - p^{n+1}_{i,j}\right)_{r+1}
      &+ \\
      +\dfrac{\Delta t^2}{\Delta x^2}
      \left(\tilde{h}_{i-1/2,j}^{n+1}\right)_r\,\left(p^{n+1}_{i,j} - p^{n+1}_{i-1,j}\right)_{r+1}
       &+ \\
      -\dfrac{\Delta t^2}{\Delta y^2}
      \left(\tilde{h}_{i,j+1/2}^{n+1}\right)_r\,\left(p^{n+1}_{i,j+1} - p^{n+1}_{i,j}\right)_{r+1}
      &+ \\
      +\dfrac{\Delta t^2}{\Delta y^2}\,\left(\tilde{h}_{i,j-1/2}^{n+1}\right)_r\,\left(p^{n+1}_{i,j} - p^{n+1}_{i,j-1}\right)_{r+1}
      & = \left(d_{i,j}\right)_r
   \end{aligned}
\end{equation}
with the known right hand side 
\begin{equation}
\begin{aligned}
\left(d_{i,j}\right)_r  =   \rho\, E_{i,j}^{*} &-  (\rho\,e_\up{t}^{n+1})_{i,j} - (\rho\,e_\up{s}^{n+1})_{i,j} - \left[\left(\rho\,\tilde{e}_\up{k}^{n+1}\right)_{i,j}\right]_r + \left[\left(\rho\,\tilde{w}_\up{g}^{n+1}\right)_{i,j}\right]_r + \\
&- \frac{\Delta t}{\Delta x} \left(\tilde{h}_{i+1/2,j}^{n+1}\right)_r\,\left[\left(\rho\, u_1^\ast\right)_{i+1/2}^{j} + \Delta t \, \left(f^\ast_1\right)_{i+1/2}^{j}\right] + \\
&+ \frac{\Delta t}{\Delta x} \left(\tilde{h}_{i-1/2,j}^{n+1}\right)_r\,\left[\left(\rho\, u_2^\ast\right)_{i-1/2}^{j} + \Delta t \, \left(f^\ast_1\right)_{i-1/2}^{j}\right] + \\ 
&- \frac{\Delta t}{\Delta y} \left(\tilde{h}_{i,j+1/2}^{n+1}\right)_r\,\left[\left(\rho\, u_2^\ast\right)_{i}^{j+1/2} + \Delta t \, \left(f^\ast_2\right)_{i}^{j+1/2}\right] + \\
&+ \frac{\Delta t}{\Delta y} \left(\tilde{h}_{i,j-1/2}^{n+1}\right)_r\,\left[\left(\rho\, u_2^\ast\right)_{i}^{j-1/2} + \Delta t \, \left(f^\ast_2\right)_{i}^{j-1/2}\right] \\
\label{eqn.rhs.2d} 
\end{aligned} 
\end{equation} 
The density at the new time $\rho_{i,j}^{n+1} = \rho_{i,j}^{\ast} $ is already known from \eqref{eqn.Qstar}, 
and so are the energy 
contribution $(\rho\,e_{\up{s}}^{n+1})_{i,j}$ of the distortion 
field $\vec{A}^{n+1}$ and the interface energy $(\rho\,e_{\up{t}}^{n+1})_{i,j}$ of the
field $\vec{b}^{n+1}$, 
after averaging onto the main grid the staggered 
field components of $\vec{b}$ and $\vec{A}$ that have been evolved in the vertices
via the compatible discretization \eqref{eqn.Jh}.  

Note that, the definitions given in Equation \eqref{eq.drhou}, are an important element of the scheme presented 
in this paper, aimed at improving its accuracy and robustess, with respect to simpler splitting techniques.

Concerning the kinetic energy contribution, it is updated explicitly via a fixed-point iteration like the enthalpy $\tilde{h}^{n+1}$
\begin{equation}
\begin{aligned}
   \left[\left(\rho\,\tilde{e}_\up{k}^{n+1}\right)_{i,j}\right]_r & = 
   \frac{1}{2}\,\rho_{i,j}^{n+1}\,{\left\{\frac{\left[\left(u_1^{n+1}\right)_{i-1/2}^{j}\right]_r + \left[\left(u_1^{n+1}\right)_{i+1/2}^{j}\right]_r}{2}\right\}}^2 + \\
   & + \frac{1}{2}\,\rho_{i,j}^{n+1}\,{\left\{\frac{\left[\left(u_2^{n+1}\right)_{i}^{j-1/2}\right]_r + \left[\left(u_2^{n+1}\right)_{i}^{j+1/2}\right]_r}{2}\right\}}^2,
\end{aligned}
\end{equation}
and the same update strategy is applied for the work due to gravity forces 
\begin{equation}
\begin{aligned}
   \left[(\rho\,\tilde{w}_\up{g}^{n+1})_{i,j}\right]_r 
   &= \frac{1}{2}\,\left[\left(\rho\,u_1^{n+1}\right)_{i-1/2}^{j} + \left(\rho\,u_1^{n+1}\right)_{i+1/2}^{j}\right]_r\,\vec{g}\cdot\uvec{e}_x + \\
   &+ \frac{1}{2}\,\left[\left(\rho\,u_2^{n+1}\right)_{i}^{j-1/2} + \left(\rho\,u_2^{n+1}\right)_{i}^{j+1/2}\right]_r\,\vec{g}\cdot\uvec{e}_y
\end{aligned}
\end{equation}

For a general equation of state, the final system for the pressure \eqref{eqn.p2d} forms a \textit{mildly nonlinear system} (see \cite{DumbserCasulli2016}) of the form 
\begin{equation}
\rho\,\mathbf{e} \left( \mathbf{p}_{r+1}^{n+1} \right) + \mathbf{M}_r \, \mathbf{p}_{r+1}^{n+1} = \mathbf{d}_r^{n} 
\label{eqn.nonlinear} 
\end{equation}
with a linear part contained in $\mathbf{M}$ that is symmetric and at least positive semi-definite. Hence, with the usual 
assumptions on the nonlinearity detailed in \cite{CasulliZanolli2012}, it can be efficiently solved with 
the nested Newton method of Casulli and Zanolli \cite{CasulliZanolli2010,CasulliZanolli2012}. 
For our particular choice of equation of state (stiffened gas), the system is \emph{linear} and 
thus we can employ an even simpler Jacobi-preconditioned matrix free conjugate gradient method for its solution.

Note that in the incompressible 
limit $\mathbb{M}\up{a} \to 0$, following the asymptotic analysis performed in \cite{KlaMaj,KlaMaj82,Klein2001,Munz2003,MunzDumbserRoller}, 
the pressure tends to a constant and the contribution of the kinetic energy $\rho\,\tilde{e}_\up{k}$  
can be neglected with respect to the internal energy $\rho\,e$. 
Therefore, in the incompressible limit the system \eqref{eqn.p2d} tends to the classic 
pressure Poisson equation used in incompressible flow solvers. 
In each Picard iteration, after the solution of the pressure system \eqref{eqn.p2d}, 
the enthalpies at the interfaces are recomputed and the momentum is updated by 
\begin{equation}
\label{eqn.rhou2d.pic} 
\begin{aligned}
\left[\left(\rho\,u_1^{n+1}\right)_{i+1/2}^{j}\right]_{r+1} & = \left(\rho\,u_1^{\ast}\right)_{i+1/2}^{j} + {\Delta t}\left(\rho_{i+1/2,j}^{n+1}\,\vec{g}\cdot\uvec{e}_x  - \frac{p^{n+1}_{i+1,j} - p^{n+1}_{i,j}}{\Delta x} \right)_{r+1} ,  \\
\left[\left(\rho\,u_2^{n+1}\right)_{i}^{j+1/2}\right]_{r+1} & = \left(\rho\,u_2^{\ast}\right)_{i}^{j+1/2} + {\Delta t}\left(\rho_{i,j+1/2}^{n+1}\,\vec{g}\cdot\uvec{e}_y  - \frac{p^{n+1}_{i,j+1} - p^{n+1}_{i,j}}{\Delta y} \right)_{r+1} ,  
\end{aligned} 
\end{equation} 
At the end of the Picard iterations, the total energy is updated as 
\begin{equation}
\label{eqn.rhoE2d} 
 \begin{aligned}
      \rho\,E^{n+1}_{i,j} = \rho\,E^{\ast}_{i,j} &
      -\dfrac{\Delta t}{\Delta x}\,\left[
      \tilde{h}^{n+1}_{i+1/2,j}\,(\rho\,u_1^{n+1})_{i+1/2}^j - 
      \tilde{h}^{n+1}_{i-1/2,j}\,(\rho\,u_1^{n+1})_{i-1/2}^j                                  
      \right] + \\ 
      &-\dfrac{\Delta t}{\Delta y}\,\left[
      \tilde{h}^{n+1}_{i,j+1/2}\,(\rho\,u_2^{n+1})_{i}^{j+1/2} - 
      \tilde{h}^{n+1}_{i,j-1/2}\,(\rho\,u_2^{n+1})_{i}^{j-1/2}                                  
      \right] + \\ 
      & + (\rho\,\tilde{w}_\up{g}^{n+1})_{i,j}.
   \end{aligned}  
\end{equation}
While for the final main-grid update of the momentum variables we compute a set of interpolated
cell-face values for the pressure field $p_{i+1/2,j}^{n+1} = (p_{i,j}^{n+1} + p_{i+1,j}^{n+1})/2$ 
and $p_{i,j+1/2}^{n+1} = (p_{i,j}^{n+1} + p_{i,j+1}^{n+1})/2$, and then update the momentum
\emph{directly one the main cell-center grid} with
\begin{equation}
\begin{aligned}
    \rho\,\vec{u}_{i,j}^{n+1} = \rho\,\vec{u}_{i,j}^{\ast} 
    & - \frac{\Delta t}{\Delta x}\,\left(p_{i+1/2,j}^{n+1} - p_{i-1/2,j}^{n+1}\right) + \\
    & - \frac{\Delta t}{\Delta y}\,\left(p_{i,j+1/2}^{n+1} - p_{i,j-1/2}^{n+1}\right) + \Delta t\,\rho_{i,j}^{n+1}\,\vec{g}.
\end{aligned}
\end{equation}
This approach further differentiates 
the method proposed in this paper from the one given in \cite{sigpr}, and 
is preferred in this work as opposed to averaging the momentum from the cell-face grid to the
cell center grid, in order to avoid the Lax--Friedrichs-type numerical diffusion associated with
incurred in the final momentum is averaged back onto the main grid. 


\subsection{Proof of the Abgrall compatibility condition}

We provide here a simple proof that the proposed scheme respects the so-called Abgrall condition \cite{abgrallcondition}, 
that is, it preserves exactly those flows characterised by a constant velocity and 
constant pressure. In absence of other
driving forces, such uniform flows should not be affected by spurious perturbations, 
regardless of the distribution of density or volume fraction, which do not affect 
the dynamics in these situations.

The starting point is showing that the velocity field produced by the convective step is 
kept uniform by the MUSCL--Hancock scheme applied to the convective subsystem.
In one space dimension, the mixture density $\rho$ obey the update formula
\begin{equation}
\label{eq.abgralldensity}
   \rho_i^{n+1}        = \rho_i^n - \frac{\Delta t}{\Delta x}\,\left(f_{i+1/2}^\rho - f_{i-1/2}^\rho\right),
\end{equation}
with the Rusanov flux yielding explicitly
\begin{equation}
   f_{i+1/2}^\rho(\rho_L,\ \rho_R) = \frac{1}{2}\,u_1\,\left(\rho_L + \rho_R\right) - \frac{1}{2}\,|u_1|\,\left(\rho_R - \rho_L\right)
\end{equation}
which is easily shown by direct sum of the equations for the phase 
densities $\alpha_1\,\rho_1$ and $\alpha_2\,\rho_2$. 
Since it is a fundamental assumption that the velocity field is constant at time level $t^n$,
in this proof we can denote $u_1 = (u_1)_i = (\rho\,u_1)_i^n/\rho_i^n$ for any cell $i$, 
The update formula for the mixture momentum $\rho\,u_1$ similarly reads
\begin{equation}
\label{eq.abgrallmomentum}
   (\rho\,u_1)_i^{n+1} = (\rho\,u_1)_i^n - \frac{\Delta t}{\Delta x}\,\left(f_{i+1/2}^{\rho\,u_1} - f_{i-1/2}^{\rho\,u_1}\right),
\end{equation}
and exploiting the constant velocity assumption,
the flux is 
\begin{equation}
   f_{i+1/2}^{\rho\,u_1} = \frac{1}{2}\,u_1\,u_1\,\left(\rho_L + \rho_R\right) - \frac{1}{2}\,|u_1|\,u_1\,\left(\rho_R - \rho_L\right) = u_1\,f_{i+1/2}^\rho,
\end{equation}
which means that \eqref{eq.abgrallmomentum}
can be simplified to
\begin{equation}
\label{eq.abgrallmomentum2}
   (\rho\,u_1)_i^{\ast} = (\rho\,u_1)_i^n - u_1\,\frac{\Delta t}{\Delta x}\,\left(f_{i+1/2}^{\rho} - f_{i-1/2}^{\rho}\right),
\end{equation}
from which, setting $(u_1)_i^{\ast} = (u_1)_i^n = u_1$ allows to factor out 
Eq. $\eqref{eq.abgralldensity}$, showing that the constant velocity field is preserved 
by the scheme. The same formal proportionality argument can be applied also to the
cell-local predictor of the MUSCL--Hancock scheme, showing that the velocity field generated
by the predictor step is unaltered in the same way.

It remains to be shown that a constant pressure field $p = p_i^n = p_i^{n+1}$ for any index $i$ is a solution 
of the discrete wave equation resulting from the energy balance
\begin{equation}
\begin{aligned}
   \rho\,e_{i}^{n+1} + (\rho\,e_\up{k})_i^{n+1} &= \rho\,e_{i}^{n} + (\rho\,e_\up{k})_i^{\ast} +\\
   &- \frac{\Delta t}{\Delta x}\,\left[
   \tilde{h}_{i+1/2}^{n+1}\,\left(\rho\,u_1\right)_{i+1/2}^{\ast} - 
   \tilde{h}_{i-1/2}^{n+1}\,\left(\rho\,u_1\right)_{i-1/2}^{\ast}\right], 
\end{aligned}
\end{equation}
together with the equivalences $(\rho\,e_\up{k})_{i}^{n+1} = (\rho\,e_\up{k})_{i}^{\ast}$ and 
$\left(\rho\,u_1\right)_{i+1/2}^{n+1} = \left(\rho\,u_1\right)_{i+1/2}^{\ast}$ resulting from the constant-pressure
assumption which means that momentum $\rho\,u_1$ and kinetic energy $\rho\,e_\up{k}$ at the new 
time level coincide with those resulting from the convective subsystem.
Collecting the constant velocity $u_1$ and plugging in a generic linear equation of state of the form
$\rho\,e = k_0 + k_1\,p$, which is the form of the stiffened gas EOS applied to our model
when $\alpha_1$ is constant throughout
the domain, gives 
\begin{equation}
\label{eq.abgrallrhoh}
\begin{aligned}
   k_0 + k_1\,p_i^{n+1} + (\rho\,e_\up{k})_i^{n+1} &= k_0 + k_1\,p_i^{n} + (\rho\,e_\up{k})_i^{\ast} +\\
   &- \frac{\Delta t}{\Delta x}\,u_1\,\left(
   \rho_{i+1/2}^{\ast}\,\tilde{h}_{i+1/2}^{n+1} - 
   \rho_{i-1/2}^{\ast}\,\tilde{h}_{i-1/2}^{n+1}\right), 
\end{aligned}
\end{equation}
which with the constant pressure assumption $p_i^{n+1} = p_i = p$ yields a condition
\begin{equation}
   \rho_{i+1/2}^{\ast}\,\tilde{h}_{i+1/2}^{n+1} - \rho_{i-1/2}^{\ast}\,\tilde{h}_{i-1/2}^{n+1} = 0,
\end{equation}
highlighting that the enthalpy estimates $\tilde{h}_{i+1/2}^{n+1}$
must be chosen as 
\begin{equation}
   \tilde{h}_{i+1/2}^{n+1} = \frac{\rho\,\tilde{h}_{i+1/2}^{n+1}}{\rho_{i+1/2}^\ast} = \frac{\rho\,e(p_{i+1/2}) + p_{i+1/2}}{\rho_{i+1/2}^\ast}, 
\end{equation}
meaning that
the density used for the computation of enthalpies must necessarily be the one produced by the convective 
step $\rho_{i+1/2}^{\ast}$.
Any interpolation scheme for computing $p_{i+1/2}$ will clearly work in a constant pressure field, and 
we use a simple arithmetic average $p_{i+1/2} = (p_i + p_{i+1})/2$, and the same average is employed for 
computing $\rho_{i+1/2}^\ast$, but in this case it is important that the interpolation operator be
the same applied for the computation of the momentum $(\rho\,u_1)_{i+1/2}^{\ast}$ from the cell-center quantities.
Also, note that in order to be able to set $(\rho\,e_\up{k})_{i}^{n+1} = (\rho\,e_\up{k})_{i}^{\ast}$, 
the kinetic energy computed from averaging the (constant) velocities from the cell center to the edges and vice-versa, 
must be compatible with that obtained by the MUSCL--Hancock scheme itself, which is verified thanks again to 
the fact that a compatible numerical dissipation was chosen for density, momentum, and kinetic energy.
Specifically, it is immediately apparent that, analogously to the momentum flux, we can write the kinetic
energy flux as $f_{i+1/2}^{\rho\,e_\up{k}} = e_\up{k}\,f_{i+1/2}^\rho$, thus the specific kinetic energy $e_\up{k} = u_1^2/2$
is kept constant after the convective step.

Given the conditions above, one can immediately see that a constant pressure field, 
with $\rho\,e_\up{k}^{n+1} = \rho\,e_\up{k}^{\ast}$ is 
solution to the discrete wave equation
\begin{equation}
\begin{aligned}
    \rho\,e_i^{n+1} + \rho\,e_\up{k}^{n+1} - \frac{\Delta t^2}{\Delta x^2}\,\left[
    \tilde{h}_{i+1/2}^{n+1}\,\left(p_{i+1}^{n+1} - p_i^{n+1}\right) 
    - \tilde{h}_{i-1/2}^{n+1}\,\left(p_{i}^{n+1} - p_{i-1}^{n+1}\right]
    \right) =& \\
    =  \rho\,e_i^{n} + \rho\,e_\up{k}^{\ast}   - \frac{\Delta t}{\Delta x}\,\left[
   \tilde{h}_{i+1/2}^{n+1}\,\left(\rho\,u_1\right)_{i+1/2}^{\ast} - 
   \tilde{h}_{i-1/2}^{n+1}\,\left(\rho\,u_1\right)_{i-1/2}^{\ast}\right]&, 
\end{aligned} 
\end{equation} 

A further condition on the scheme must be imposed whenever the volume fraction $\alpha_1$ is not 
constant in space.
In this case, the stiffened gas equation of state applied to each phase provides a more complex closure law 
of the type $\rho\,e = \alpha_1\,\rho_1\,e_1 + \alpha_2\,\rho_2\,e_2$ or
$\rho\,e = \alpha_1\,k_{01} + \alpha_1\,k_{11}\,p + \alpha_2\,k_{02} + \alpha_2\,k_{12}\,p$.
Applied to the discrete wave equation, the mixture equation of state gives
\begin{equation}
   (\alpha_1)_i^{n+1} = (\alpha_1)_i^n  - \frac{\Delta t}{\Delta x}\,\left[u_1\,(\alpha_1)_{i+1/2}^{\ast} - u_1\,(\alpha_1)_{i-1/2}^{\ast}\right],
\end{equation}
that is, at least when the velocity field $u_1$ is a constant, the scheme for the update of $\alpha_1$
must coincide with one in flux form, for some appropriate choice of $(\alpha_1)_{i+1/2}^\ast$.

For a constant velocity field, the nonconservative products not associated with pure 
convection of the volume fraction vanish and, the update scheme reads
\begin{equation}
   (\alpha_1)_{i}^{n+1} = (\alpha_1)_{i}^{n} - \frac{\Delta t}{\Delta x}\,\left(D_{i+1/2}^{\alpha_1} + D_{i-1/2}^{\alpha_1}\right)
\end{equation}
with the path-conservative fluctuations as well as the numerical dissipation from the Rusanov flux collected in 
\begin{equation}
   \begin{aligned}
      & D_{i+1/2}^{\alpha_1} = u_1\,\frac{1}{2}\,\left[\flll(\alpha_1)_R - (\alpha_1)_L\right]_{i+1/2} - |u_1|\,\frac{1}{2}\,\left[\flll(\alpha_1)_R - (\alpha_1)_L\right]_{i+1/2},\\
      & D_{i-1/2}^{\alpha_1} = u_1\,\frac{1}{2}\,\left[\flll(\alpha_1)_R - (\alpha_1)_L\right]_{i-1/2} + |u_1|\,\frac{1}{2}\,\left[\flll(\alpha_1)_R - (\alpha_1)_L\right]_{i-1/2}.\\
   \end{aligned}
\end{equation}
This automatically gives raise to an upwind discretisation that suggests the interpolated
values of $\alpha_{i+1/2}^\ast$ should be computed with the same upwinding rule
\begin{equation}
   (\alpha_1)_{i+1/2}^\ast\left(\flll(\alpha_1)_L,\ (\alpha_1)_R\right) = 
   \frac{1 + \sign{(u_1)}}{2}\,(\alpha_1)_L + \frac{1 - \sign{(u_1)}}{2}\,(\alpha_1)_R, 
\end{equation}
for any left and right states $(\alpha_1)_L$ and $(\alpha_1)_R$ obtained from the predictor step of
the MUSCL--Hancock scheme.
Then it can be verified that for any distribution of volume fraction $\alpha_1$ and density $\rho$
the discrete wave equation will indeed preserve constant-pressure, constant-velocity solutions exactly.

%% file: methodsource.tex
\section{Semi-analytic integration of strain relaxations sources}
\label{sec:semianalytical}

A necessary element for the successful solution of the unified model of continuum mechanics is the
accurate integration of the distortion matrix $\vec{A}$. 

The evolution dynamics of the distortion matrix $\vec{A}$ and of the metric tensor $\Ge = \vec{A}^\transpose\,\vec{A}$ 
take place on a wide span of timescales: given a fixed evolution speed of
the kinematics of distortion (due to flow convection and velocity gradients), 
one can find anything from infinitely slow strain relaxation 
in elastic solids, to infinitely fast shear dissipation in perfect fluids, with viscous fluids
also being a nontrivial example of fast-acting (stiff) strain relaxation.

From the mathematical standpoint, such timescales can be quantified by means of a relaxation time $\tau$ 
in the evolution equation of the distortion matrix
\begin{equation}
\label{eq.strainrel}
    \partial_t \vec{A} + \left(\nabla\vec{A}\right)\,\vec{u} + \left(\nabla\vec{u}\right)\,\vec{A} = \vec{Z} = 
    -\frac{3}{\tau}\,\left(\det\vec{A}\right)^{5/3}\,\vec{A}\,\dev\left(\vec{A}^\transpose\,\vec{A}\right) 
\end{equation}
and in the corresponding equation for the metric tensor
\begin{equation} 
\label{eq.strainrelg}
\partial_t\vec{G} + \left(\nabla\vec{G}\right)\,\vec{u} + \Ge\,\nabla\vec{u} + \left(\nabla\vec{u}\right)^\transpose \, \Ge = - \frac{6}{\tau}\,{(\det{\Ge})}^{5/6}\,\Ge\,\dev{\Ge}.
\end{equation}
The relaxation time $\tau$, in principle a function of the state variables, but often a fixed constant, 
is what defines the stiff nature of the algebraic source terms governing the relaxation towards an equilibrium 
state of the material strain. To highlight the connection between the distortion matrix $\vec{A}$, 
the metric tensor $\vec{G}$ and what we generically call \textit{strain}, it should be recalled that, 
in a purely elastic context, for small deformations, the linear strain $\vec{\epsilon}$ can be 
directly expressed as
$\vec{\epsilon} = \left(\vec{I} - \vec{G}\right)/2$, for this reason we refer to the above right hand side terms 
as \textit{strain relaxation sources}.

One of the major difficulties in the solution of the unified model of continuum mechanics is indeed
the presence of these nonlinear source terms. In the past, the locally implicit ADER treatment of 
source terms has proven effective, as well as the splitting or fractional step approach, in conjunction
with the implicit Euler scheme for stable time integration. However, we found that 
a new approach has to be adopted 
for certain choices
of the material parameters, for example for extremely fast relaxation times in complex flows, or 
for the nonlinearly stress-dependent timescales encountered in the application of 
the model to material failure dynamics.

A final major step forward in the development of a robust solver for the strain relaxation system~\eqref{eq.strainrel}, 
in particular allowing to accurately capture the Navier--Stokes limits regardless of the timestep size, is
based on three key observations:
\begin{enumerate}
\item The splitting approach is not always adequate
\item The structure of the problem can be significantly simplified by choosing the 
appropriate reference frames
\item Equilibrium states can be computed algebraically without time integration
\end{enumerate}
The details regarding the semi-analytical solution strategy adopted in this work are given in 
the following paragraphs.

\subsection{Limits of the splitting approach}
In previously discussed techniques \cite{chiocchettimueller, tavellicrack} 
for the solution of relaxation processes we have adopted the 
fractional step (or splitting) approach. The technique is very useful as it allows to separate
the solution of the relaxation source from all other dynamics, and attack the resulting 
ordinary differential equation system with ad hoc techniques.
However, the relaxation processes in the unified model of continuum mechanics, besides
complex nonlinear dynamics, also feature nontrivial
equilibrium states that \textit{must} be reliably captured and preserved. If not, important properties
of the continuum model, like the convergence to the Navier--Stokes--Fourier system \cite{GPRmodel}, 
may be lost in its discrete transposition.

To quantitatively argue the point, we 
highlight the problem with regards to a simpler example,
namely given by the thermal impulse equation \cite{boscherisigpr}, in the simple case of a vanishing velocity field, which is
\begin{equation}
   \label{eq.thermalimpulsesource}
   \partial_t\left(\vec{J}\right) + \nabla T = -\vec{J}/\taut,
\end{equation}
and assume available a generic numerical scheme
by which one can compute for each cell/degree of freedom an update
$\Pstar = (\vec{J^\ast} - \vec{J}^n)/\Delta t$ such that $\vec{J^\ast}$ is the solution to 
the update of the left hand side of \eqref{eq.thermalimpulsesource}, 
i.e. the homogeneous system that is solved by application of the splitting approach.
In this particular case $\Pstar$ can be seen as a discretisation of $-\nabla T$.
Then, a straightforward application of the splitting method would find the solution at time 
$t^{n+1} = t^n + \Delta t$
of the initial value problem
\begin{equation}
   \label{eq.jrelaxation}
   \frac{\de\vec{J}}{\de\vec{t}} = -\frac{1}{\taut}\,\vec{J}, \quad t\in[t^n,\ t^{n+1}], \quad \vec{J}(t^n) = \vec{J}^\ast.
\end{equation}
If Equation \eqref{eq.jrelaxation} is integrated via the implicit Euler method, then one can prove
that the final value of $\vec{J}^{n+1}$ will indeed yield an asymptotic preserving discretization of
the PDE (provided obviously that the discretisation of the left hand side is compatible).
However, it is easy to see that, isolated from the left hand side of \eqref{eq.thermalimpulsesource},
the ordinary differential problem \eqref{eq.jrelaxation} asymptotically relaxes $\vec{J}$ to zero if the relaxation 
time is sufficiently small with respect to the timestep size, regardless of the value of the initial condition $\vec{J}^\ast$.
This implies that if one were to integrate \eqref{eq.jrelaxation} \emph{exactly}, then the updated value 
of the thermal impulse $\vec{J}$
would be $\vec{J}^{n+1} = \vec{0}$ instead of $\vec{J}^{n+1} = -\taut\,\nabla T$, highlighting a very clear issue in a naive application
of the splitting approach.

In order to overcome this issue, a simple modification to the ordinary differential problem \eqref{eq.jrelaxation}
allows to account for the left hand side of \eqref{eq.thermalimpulsesource} and thus 
converge to the correct asymptotic state $\vec{J} = -\taut\,\nabla T$ in the stiff limit $\tau_h \to 0$.


An alternative ordinary differential problem to be solved is then 
   \begin{equation}
   \label{eq.goodjode}
      \dfrac{\de{\vec{J}}}{\de{t}} = \Pstar - \dfrac{1}{\taut}\,\vec{J}, \quad t\in[t^n,\ t^{n+1}],  \quad \vec{J}(t^n) = \vec{J}^\ast, 
   \end{equation}
   where, as stated, $\Pstar$ accounts for the discrete update from the left hand side of \eqref{eq.thermalimpulsesource}.

   Again, \eqref{eq.goodjode}
   can be seen as a 
   system of three uncoupled first order linear ordinary differential equations (ODEs) and 
   an exact solution is indeed found thanks to the linearity and independence of
   the three equations.
   Explicitly, the solution is 
   \begin{equation} \label{eqn.heatsolution}
      \vec{J}^{n+1} = (\vec{J}^n - \taut\,\Pstar)\,\exp(- \Delta t/\taut) + \taut\,\Pstar.
   \end{equation}
   The only degenerate case to be considered is that
   if $\Delta t/\taut$ is very small (of the order of $10^{-8}$),
   i.e. if the source term is not stiff at all, then \eqref{eqn.heatsolution} might yield 
   inaccurate results, due to floating point representation issues.
   In this case, one may simply opt to to switch to explicit Euler integration, 
   which for such mild (vanishing) sources yields perfectly valid solutions.
   A similar application

\subsection{Simplification of the problem by polar decomposition and principal axes coordinates.}

The nine components of the distortion matrix/basis triad $\vec{A}$ encode two different kinds of information:
six degrees of freedom are 
directly linked to the stress tensor $\vec{\sigma} = -\rho\,\cshear^2\,\vec{G}\,\dev\vec{G}$, while 
the remaining three are associated with an angular orientation which does not influence 
stresses or energies but are nonetheless part of the structure of the governing equations.
This can be formalized by means of the polar decomposition of $\vec{A}$, by which we can highlight
the six stress-inducing components of $\vec{A}$, identifying them as the square root $\vec{G}^{1/2}$ of
the metric tensor $\vec{G}$. 
Moreover, one can easily see that, if an appropriate \textit{fixed} transformation of the reference 
frame is applied to \eqref{eq.strainrel} (a polar decomposition followed by a spectral decomposition), such that at time
$t = t^n$ one has $\vec{A}$ in 
diagonal form, and if the convection/production term on the left hand side of \eqref{eq.strainrel} is 
null (if the flow field is uniform, i.e. $\nabla\vec{u} = \vec{0}$, or if formally we want 
to study the invariance properties of the strain relaxation source), 
then the diagonality of $\vec{A}$ is maintained for all $t \ge t^n$. This means that 
the relaxation source on the right hand side of \eqref{eq.strainrel} does not alter the rotational component
of $\vec{A}$. In the following we establish the notation for the polar decomposition procedure enabling 
separate treatment of rotational degrees of freedom of $\vec{A}$ and volumetric/shear/relaxation effects
and provide some formal justification of the validity of the approach.

\subsubsection{Polar decomposition of the distortion matrix}
Given the definition of the metric tensor $\vec{G} = \vec{A}^\transpose\,\vec{A}$, the distortion
matrix $\vec{A}$ can always be expressed as 
\begin{equation}
\vec{A} = \vec{R}\,\vec{G}^{1/2}, \quad\text{with}\quad \vec{G}^{1/2} = \vec{E}\,\hat{\vec{G}}^{1/2}\,\vec{E}^{-1},    
\end{equation}
with $\vec{R}$ an orthogonal transformation with positive unitary determinant, i.e. a rotation matrix.
Numerically, the matrix square root $\vec{G}^{1/2}$ can be evaluated by means of the Denman--Beavers algorithm, 
or alternatively, thanks to the symmetry of $\vec{G}$, one can reliably and accurately 
compute the eigenvectors $\vec{E}$ and diagonal form $\hat{\vec{G}}$ 
from the eigen-decomposition of the metric tensor
\begin{equation}
    \vec{G} =\vec{E}\,\hat{\vec{G}}\,\vec{E}^{-1}
\end{equation}
with the Jacobi eigenvalue algorithm. In the work reported in the present paper, 
the latter is indeed the method of choice for the task.
This allows, for any given state $\vec{A}$, to compute a rotation matrix
\begin{equation}
    \vec{R} = \vec{A}\,\vec{G}^{-1/2},\quad\text{with}\quad\vec{G}=\vec{A}^\transpose\,\vec{A}
\end{equation}
which allows to apply any operation to $\vec{G}=\vec{A}^\transpose\,\vec{A}$, 
or its square root $\vec{G}^{-1/2}$, exploiting their symmetric and positive definite properties,
or, temporarily in the solution procedure and locally in space, 
adopting a simpler form of the governing equations, written in terms of $\vec{G}$ instead
of $\vec{A}$.
Then the effects of such operations can be mapped onto $\vec{A}$ directly via the 
rotation $\vec{A} = \vec{R}_{\ast}\,\vec{G}^{-1/2}$, with $\vec{R}_\ast$ a rotation matrix that can 
be computed independently from the nonlinear source terms as $\vec{R}_\ast = \vec{A}_\ast\,\vec{G}^{-1/2}_\ast$, 
having defined $\vec{A}_\ast$ the distortion matrix obtained from the left hand side of the evolution
equation as customarily done per the fractional step method.
To summarize, one may first compute $\vec{R}_n = \vec{A}_n\,\vec{G}^{-1/2}_n$, use it to map to 
an auxiliary frame in which one can easily integrate the source term as applied to the symmetric 
positive definite metric tensor $\vec{G}$, and then map back to $\vec{A}$ by means of a \textit{different} rotation
matrix $\vec{R}_\ast = \vec{A}_\ast\,\vec{G}^{-1/2}_\ast$, already obtained as a function of the left hand side only.

\subsubsection{Invariance under strain relaxation of the rotational component of the distortion matrix}
We can study the effects of the relaxation source 
on the distortion matrix $\vec{A}$ 
and the metric tensor $\vec{G}$ in isolation from those of flow gradients
by formally setting $\nabla\vec{u} = \vec{0}$
obtaining the uniform flow equations
\begin{align}
   &\partial_t\vec{A} + \left(\nabla\,\vec{A}\right)\,\vec{u} = - k\,\vec{A}\,\dev\vec{G},\\[1mm]
   &\partial_t\vec{G} + \left(\nabla\,\vec{G}\right)\,\vec{u} = - 2\,k\,\vec{G}\,\dev\vec{G},
\end{align}
with, $k = 3\,\tau^{-1}\,\left(\det\vec{G}\right)^{5/6}$, or, 
in 
the co-moving reference, equivalently
\begin{align}
   &\frac{\de \vec{A}}{\de t} = - k\,\vec{A}\,\dev\vec{G},\label{eq.lagrA}\\[1mm]
   &\frac{\de \vec{G}}{\de t} = - 2\,k\,\vec{G}\,\dev\vec{G},\label{eq.lagrG}
\end{align}
where $\de/\de t$ is the customary notation for the total/convective/Lagrangian derivative.
In the following paragraphs we derive a set of evolution equations, valid in uniform flow, 
where only strain relaxation effects can be observed, which are provided as an argument for the
separate computation of rotations and volumetric/shear effects in the proposed semi-analytical solver.

\paragraph{Evolution equation for the square root of the metric tensor.}
The metric tensor equation \eqref{eq.lagrG} can be mapped to a principal reference frame, by setting
$\vec{G} = \ege$, with $\vec{E}$ the matrix of eigenvectors of $\vec{G}$ computed
at a fixed time $t = t_\ast$ and $\Gh$ a diagonal matrix defined at any time $t$ by 
$\Gh(t) = \vec{E}^{-1}\,\vec{G}(t)\,\vec{E}$.
Then \eqref{eq.lagrG} can be rewritten as
\begin{equation}
   \frac{\de}{\de t}\left(\ege\right) = - 2\,k\,\egdge,\\
\end{equation}
and, since $\vec{E}$ is fixed at time $t = t_\ast$, we obtain
\begin{equation}
   \frac{\de\Gh}{\de t} = -2\,k\,\Gh\dev\Gh, 
\end{equation}
which implies that for a fixed orthonormal transformation $\vec{E}$ 
such that at a given time $t_\ast$ one has a diagonal 
$\Gh(t_\ast) = \vec{E}^{-1}\,\vec{G}(t_\ast)\,\vec{E}$, 
the matrix $\Gh$ will remain diagonal for all times, implying that 
the strain relaxation source does not affect the eigenvectors of $\vec{G}$ but only its 
eigenvalues.
Moreover, since in the fixed principal axes reference frame $\Gh$ and its square root
are guaranteed to be diagonal at any time $t$, then one can write
\begin{equation}
\begin{aligned}
   \frac{\de \Gh}{\de t} &= \frac{\de}{\de t}\left(\rGh\,\rGh\right) = \rGh\,\ddt\left(\rGh\right)  + \ddt\left(\rGh\right) \,\rGh = \\[1mm]
    &= 2\,\rGh\,\ddt\left(\rGh\right),
\end{aligned}
\end{equation}
since, in this specific reference frame, all matrix involved are diagonal and products are commutative.
Therefore the evolution equation for the square root of $\Gh$, under the effects of strain relaxation only,
and in the principal frame reads
\begin{equation}
\label{eq.rtg}
   \ddt\left(\rGh\right) = \frac{1}{2}\,\Gh^{-1/2}\,\frac{\de \Gh}{\de t}  = -k\,\rGh\,\dev\Gh.
\end{equation}

\paragraph{Evolution equation for the rotational component of the distortion matrix.}
A simple governing equation for the rotational component $\vec{R}$ 
of the distortion matrix $\vec{A}$ can be derived, 
in uniform flow $\nabla\vec{u} = \vec{0}$, or equivalently under the effect of the strain relaxation
source only,
by substituting the polar decomposition 
\begin{equation}
   \vec{A} = \vec{R}\,\vec{G}^{1/2} = \vec{R}\,\erge
\end{equation}
in \eqref{eq.lagrA}, obtaining
\begin{equation}
   \ddt\left(\vec{R}\,\erge\right) = -k\,\vec{R}\,\ergdge.
\end{equation}
Then, expanding the derivative, we have
\begin{equation}
\label{eq.ralmost}
   \vec{R}\,\vec{E}\,\ddt\left(\rGh\right)\,\vec{E}^{-1} + \frac{\de \vec{R}}{\de t}\,\erge = -k\,\vec{R}\,\ergdge,
\end{equation}
and by substituting \eqref{eq.rtg} in \eqref{eq.ralmost} we recover
\begin{equation}
   \vec{R}\,\vec{E}\,\left(-k\,\rGh\,\dev\Gh\right)\,\vec{E}^{-1} + \frac{\de \vec{R}}{\de t}\,\erge
   = -k\,\vec{R}\,\ergdge,
\end{equation}
and hence the simple (non)-evolution law for the rotations of $\vec{A}$ reads
\begin{equation}
   \frac{\de \vec{R}}{\de t} = \partial_t\vec{R} + \left(\nabla\vec{R}\right)\vec{u} = \vec{0},
\end{equation}
which means that the rotational component $\vec{R}$ is not affected by the relaxation source and can only change
due to non-uniformity of the flow, unlike the principal components of the metric tensor $\vec{G}$, accounting for volumetric
and shear effects.

As early as 2010, in \cite{Pesh2010}, S. Godunov and I. Peshkov, recognised 
that the structure of the governing equations of the cobasis $\vec{A}$ could be exploited to simplify
the relaxation procedure. 
Similarly, in \cite{Jackson2019, Jackson2019a}, Jackson and Nikiforakis
develop and employ an approximate semi-analytical integration technique for 
the relaxation source of the UMCM, which however is tied to a periodic reset procedure for 
the cobasis $\vec{A}$ and for strain energy, restoring a new local \textit{small deformations} 
configuration whenever necessary. Hence the motivation for developing a novel 
integration methodology capable of handling arbitrary deformations and extreme parameter choices 
potentially varying between
$\tau = 10^{20}$ and $\tau = 10^{-14}$ as a function of stress as in \cite{boscheriulaggpr}, 
without modification to the structure of the PDE system.

\subsection{Direct solution of strain equilibrium states}

A final welcome fact about the particular structure of the 
strain production/relaxation equations \eqref{eq.strainrel} and 
\eqref{eq.strainrelg} is that, if the source term is very stiff, one may easily compute 
the asymptotic solution without necessarily resolving the stiff dynamics: the asymptotic state
does not depend on the initial conditions and can be computed by means of a simple and quickly convergent
fixed point iteration scheme.

Moreover, from the numerical standpoint, time-dependent closed form solutions of the strain relaxation equations
can be less accurate than expected for certain initial conditions, parameter choices, or flow configurations.
In such cases, the availability of the solution for the sought equilibrium state is not only a 
matter of efficiency, but also of accuracy and robustness of the computations.

Additionally, the fixed point iteration allows to prove that the integration 
scheme for the relaxation source will indeed provide, in the Navier--Stokes limit of the model, 
the correct unique asymptotic stress tensor, reflecting the convergence limit of the continuum model 
to the Navier--Stokes equations.

\section{Semi-analytic solution of the strain relaxation equations with finite relaxation time}
\label{sec.gprsource}

The solver for the strain relaxation source is based on the exponential integrator developed in \cite{tavellicrack} for the computation
of diffuse interface fractures and material failure, but exploits in a deeper manner the particular structure 
of the equation being solved, following the technique presented in \cite{chiocchettimueller} for 
finite-rate pressure relaxation. 

We recall that the solver employed in \cite{tavellicrack} required, in general, the 
solution of a sequence of a nonhomogeneous nine-by-nine systems of linear ordinary differential 
equations for the nine independent components of the distortion matrix $\mathbf{A}$, 
which involves the numerical computation of matrix exponentials and the inversion of 
the Jacobian matrix of the ODE system. Both these operations constitute delicate tasks in linear algebra
that require special care to be carried out in an efficient and accurate manner.

The approach used in this work entirely foregoes the solution of such nine-by-nine systems (six-by-six, in 
the case of the symmetric tensor $\Ge$) and the associated linear algebra intricacies. Instead, we 
compute the analytical solution to one of several different linearised equations that 
approximate the nonlinear ordinary differential equation
\begin{equation} \label{eqn.strainode}
\od{\Ge}{t} =  \Lstar - \frac{6}{\tau}\,{(\det{\Ge)}}^{5/6}\,\Ge\,\dev{\Ge},
\end{equation}
while admitting simple solutions that can be computed in a robust fashion. Here
with $\Lstar$ we denote a \textit{constant convective/productive forcing term} to be 
given in the following paragraphs, in analogy to the previously defined $\Pstar$
discrete time derivative of the thermal impulse.


An important aspect of the scheme is that it avoids fractional-step-type splitting, so that
the Navier--Stokes stress tensor and the Fourier heat flux can be recovered regardless 
of the ratio between the computational timestep 
size and the relaxation timescales. This means that the global timestep size need not be adjusted to 
accommodate for the fast dynamics of the relaxation sources. This is achieved by first computing, cell by cell, 
the update to $\Ge$ (or to $\vec{J}$) associated with the left hand side of \eqref{eq.strainrelg}, or of \eqref{eq.goodjode} for
heat conduction, and then including its effects 
in \eqref{eqn.strainode},  in the form of the constant forcing term $\Lstar$.
Formally, the first step amounts to computing the 
solution $\Ge_\ast = \Ge(t^{n+1})$ to the initial value problem

\begin{equation} \label{eqn.ivplhs}
\left\{
   \begin{aligned}
   & \dfrac{\de\Ge}{\de t} + \left(\nabla\vec{G}\right)\,\vec{u} + \Ge\,\nabla\vec{u} + \left(\nabla\vec{u}\right)^\transpose \, \Ge = \vec{0},\\[2mm]
   & \Ge(t^n) = \Ge_n.
   \end{aligned}
\right.
\end{equation}
In our case, instead of Eq. \eqref{eqn.ivplhs}, we solve the more general equation for the distortion matrix
\begin{equation} \label{eqn.ivplhsa}
\left\{
   \begin{aligned}
   & \dfrac{\de\vec{A}}{\de t} + \left(\nabla\vec{A}\right)\,\vec{u} + \left(\nabla\vec{u}\right) \, \vec{A} = \vec{0},\\[2mm]
   & \vec{A}(t^n) = \vec{A}_n.
   \end{aligned}
\right.
\end{equation}
from which we obtain a pointwise update
to the cell averages or degrees of freedom, which is in turn used to define the constant convective/productive forcing term
\begin{equation}
 \Lstar = \dfrac{\Ge_\ast - \Ge_n}{\Delta t}, \text{ with } \vec{G}_\ast = \vec{A}_\ast^\transpose\,\vec{A}_\ast,\ \vec{G}_n = \vec{A}_n^\transpose\,\vec{A}_n.
\end{equation}

Then a sub-timestepping loop with adaptive step size $\delta t^m$ is entered in order to approximate the solution of \eqref{eqn.strainode} 
with a sequence of solutions of linearised ODEs. Such a sub-timestepping loop is 
useful for ensuring 
the robustness and accuracy of the solver in complex flow configurations, 
but generally in our computations
the solver achieves convergence in one single sub-timestep $\delta t^m = \Delta t$.
We refer to \cite{chiocchettimueller,tavellicrack} for more details on the sub-timestepping approach, 
and we carry on our presentation of the method by listing three possible approximation choices for the solution 
of \eqref{eqn.strainode}.

\subsection{Approximate analytical solution for strain, approach 1}

When dealing with fluid flows (that is, when the source term acts on fast timescales), 
rather often one may assume that $\Ge$ be a perturbation of a spherical 
tensor, that is, $\dev{\Ge}$ can be assumed small.
Then, it is advantageous to rewrite \eqref{eqn.strainode} as 
\begin{equation} \label{eqn.smalldevg}
\od{\Ge}{t} = \Lstar - k\,\dev{\Ge}\,\dev{\Ge} + 
   k\,{\left(\dfrac{\tr \Ge}{3}\right)}^2\,\mathbf{I} - k\,\dfrac{\tr \Ge}{3}\,\Ge,
\end{equation}
with $k = {6}\,{\det(\Ge)}^{5/6}/\tau$ taken constant for the sub-timestep. This splits the source in four
pieces. The first is the constant $\Lstar$, associated with convection which, by definition, cannot be stiff as its size
is limited by the CFL constraint of the global timestepping scheme. The second is a (small by hypothesis) quadratic term
in $\dev{\Ge}$ which can be safely approximated as constant. 
The third is a function of $\tr \Ge$ only, again formally taken constant. This assumption 
can be justified by writing the evolution equation for the trace of $\Ge$
\begin{equation}
   \dfrac{\de}{\de{t}}\left(\tr{\Ge}\right) = \Lstar - k\, \tr{\left(\dev{\Ge}\,\dev{\Ge}\right)},
\end{equation}
which shows that either $\tr{\Ge}$ varies on a timescale associated with 
convection (by definition, slow), or as a quadratic function of $\dev{\Ge}$ (small by assumption).
The approximate equation \eqref{eqn.smalldevg} then admits the simple exact solution

\begin{equation}
   \Ge_{m+1} = \Ge(t^m + \delta t^m) = \exp\left(-k\,\dfrac{\tr \Ge}{3}\,\delta t^m\right)\,\left(\Ge_m + \mathbf{F}_0\right) - \mathbf{F}_0,
\end{equation}
with 
\begin{equation}
    \mathbf{F}_0 = -\dfrac{3}{k\,\tr{\Ge}}\,\left[\Lstar + 
        k\,{\left(\dfrac{\tr{\Ge}}{3}\right)}^2\,\mathbf{I} - k\,\dev{\Ge}\,\dev{\Ge}\right].
\end{equation}

We should remark that \textit{nowhere} in this approximate solution we neglected the contributions due to 
$\dev{\Ge}$, they only have taken to be constant for a sub-timestep. Specifically, 
our constant approximations 
are initially set to $\tr\Ge = \tr\Ge_m$ and $\dev \Ge = \dev\Ge_m$ and then updated within a fixed point iteration
as $\tr\Ge = \tr(\Ge_m + \Ge_{m+1})/2$ and $\dev \Ge = \dev(\Ge_m + \Ge_{m+1})/2$.

\subsection{Approximate analytical solution for strain, approach 2}

Whenever the deviatoric part of $\Ge$ cannot be assumed small, i.e. in practice when 
\begin{equation}
   \sqrt{\tr{\left(\dev{\Ge_m}\,\dev{\Ge_m}\right)}} > {(\det{\Ge_m})}^{1/3}/5, 
\end{equation}
better accuracy in the approximation of \eqref{eqn.strainode} can be obtained by observing 
that it is possible to switch the order of the operands of the matrix product $\Ge\,\dev{(\Ge)}$ 
appearing in \eqref{eqn.strainode}. Thus we rewrite \eqref{eqn.strainode} as
\begin{equation} \label{eqn.largedevg}
\od{\Ge}{t} = \Lstar - k\,\dev{(\Ge)}\,\Ge, 
\end{equation}
where we will take $\Lstar$, $k = {6}\,{(\det{\Ge})}^{5/6}/\tau$, and $\dev{(\Ge)}$ to be constant at each sub-timestep.
In order to simplify the solution of \eqref{eqn.largedevg}, we 
work in the principal reference frame which diagonalizes $\Ge$ and $\dev{(\Ge)}$, i.e. we compute the orthonormal matrix
$\mathbf{E}$ such that $\hat{\Ge} = \mathbf{E}^{-1}\,{\Ge}\,\mathbf{E}$
is a diagonal matrix and apply the associated change of basis to all vectors and tensors in our equation. 
This way the exact solution to Eq. \eqref{eqn.largedevg} is 
\begin{equation}
\begin{aligned}
   &\Ge_{m+1} = \Ge(t^m + \delta t^m) = \mathbf{E}\,\left\{\exp\left(-k\,\dev{\hat{\Ge}}\,\delta t^m\vphantom{l^{l^{l^l}}}\right)\cdot\right.&\\
   &\left.\cdot\left[
   \mathbf{E}^{-1}\,\Ge_n + \dfrac{1}{k}\,{\left(\dev{\hat{\Ge}}\right)}^{-1}\,\mathbf{E}^{-1}\,\Lstar
   \right]                 - \dfrac{1}{k}\,{\left(\dev{\hat{\Ge}}\right)}^{-1}\,\mathbf{E}^{-1}\,\Lstar\right\}.
\end{aligned}
\end{equation}
The three-by-three matrix $\mathbf{E}$ having for columns the eigenvectors of $\Ge$ can be quickly and robustly computed 
to arbitrary precision by means of Jacobi's method for the eigenstructure of symmetric matrices, 
and its inverse is simply given by $\mathbf{E}^{-1} = \mathbf{E}^{\transpose}$.
Furthermore, $\dev{\hat{\Ge}}$ can be inverted trivially in the principal reference frame by just taking 
the reciprocal of each diagonal entry. Like for the previous solution we iteratively update the estimate $\dev \hat{\Ge} = \dev(\hat{\Ge}_m + \hat{\Ge}_{m+1})/2$.


\paragraph{Determinant constraint.}
In the solution of the equation for the metric tensor $\Ge$, specifically when the computation
involves fluid-type behaviour, special care must be paid to preserve the nonlinear algebraic constraint
$\det{\Ge}(t,\ \vec{x}) = \left[\rho(t,\ \vec{x})/\rho_0\right]^2 = D(t,\ \vec{x})$. 
The
constraint must be actively enforced since the discretisation scheme may in principle introduce 
significant errors that over time could let the solution drift away from a compatible state in which
$\det\vec{A} \neq \rho/\rho_0$. A simple approach to the problem
constitutes in uniformly multiplying all components of $\Ge$ by $(D/\det\Ge)^{1/3}$ so that the 
resulting determinant is $D$.

The specific numerical value of the determinant $D$ is clearly known (from density) 
at the time levels $t^n$ and $t^{n+1}$, 
however it must be somehow approximated for all
the in-between times during which we operate our sub-timestepping procedure. 
In this work we impose that for a given sub-timestep indexed by $m$, connecting $t^m$ and $t^{m+1}$ the 
determinant $D$ be computed as 
\begin{equation} \label{eqn.determinant}
   D = \beta_s\,D_s + (1 - \beta_s)\,D_f,
\end{equation}
where we define $D_s = \det\left[\Ge + (t^m + \delta t^m - t^n)\,\Lstar\right]$ to be the value that
the determinant would have following a linear segment path connecting the two 
states $\Ge$ and $\Ge_\ast$, that is, the value that would allow preserve
exact integration of the (zero) source term in the solid case. 
Instead, $D_f = \det{\Ge} + (t^m + \delta t^m - t^n)\,(\det{\Ge_\ast} - \det\Ge_{n})$ 
is a second order approximation of the determinant in the fluid limit.
The mixing ratio $\beta_s$ for the two approximations $D_s$ and $D_f$ is 
an heuristic measure of how close to a solid can the material be considered and its
expression is 
\begin{equation}
   \beta_s = \min{\left[1,\ \dfrac{||\Lstar||_2^2}{||6/\tau\,(\det{\Ge_m})^{5/6}\,\Ge_m\,\dev{\Ge_m}||_2^2 + 10^{-14}}\right]}^{4}, 
\end{equation}
with $||\mathbf{A}||_2^2$ denoting the square of the Frobenius norm of a given tensor $\mathbf{A}$.

\subsection{Approximate analytical solution for strain, approach 3: Fixed point iteration for the Navier--Stokes equilibrium state}
Oftentimes the timescale $\tau$ of strain relaxation is so fast that
one may decide to just compute the strain state for which the forcing term due to
convection $\Lstar$ and the relaxation source are balanced yielding an 
equilibrium state corresponding to the Navier--Stokes limit of the GPR model.
Such an equilibrium state can be easily computed by means of a fixed point iteration in the form
\begin{equation} \label{eqn.fpi}
       \Ge_{i+1} = \widetilde{\Ge}_i{\left(\dfrac{D}{\det\widetilde{\Ge}_i}\right)}^{1/3}, \text{ with } \widetilde{\Ge}_i = 
       \dfrac{\tau\,\dev\left(\Ge^{-1}_i\,\Lstar\right)}{6\,(\det\Ge_i)^{5/6}} + 
       \dfrac{\tr\vec{G}_i}{3}\,\mathbf{I},
 \end{equation} 
 with $D$ the target determinant as defined in \eqref{eqn.determinant} and $i$ the iteration index in the fixed point procedure.
 We found that the fixed point iteration \eqref{eqn.fpi} is always convergent regardless of the initial guess, 
 but nonetheless we care to provide a simple and efficient choice in the form
 \begin{equation}
    \Ge_{1} = \widetilde{\Ge}_0{\left(\dfrac{D}{\det\widetilde{\Ge}_0}\right)}^{1/3} , 
        \text{ with } \widetilde{\Ge}_0 = \mathbf{I} + \dfrac{\tau}{6\,\det(\Ge_m)^{5/6}}\dev{\Lstar}.
 \end{equation}
Details on the derivation of this fixed point iteration scheme, as well as a proof of convergence are given in \cite{chiocchettithesis}.

 \subsection{Summary of the selection procedure for the approximation method}
 At each sub-timestep between $t^m$ and $t^{m+1}$, our solver for the equation of the elastic metric tensor $\Ge$ has to select
 the optimal approximation method for the specific distortion configuration at hand.
 The selection procedure is carried out as follows:
 \begin{enumerate}
    \item If the source is not stiff, i.e. if $\beta_s > 1 - 10^{-14}$, then we use explicit Euler integration and compute
    the solution at the next time sub-level $\Ge_{m+1} = \Ge_m + \Delta t\,\left(\Lstar - \dfrac{6}{\tau\,\det{(\Ge_m)}^{5/6}\,\Ge\,\dev{\Ge}}\right)$.
   \item Else define the indicator matrix $\mathbf{\Lambda} = \abs{\left(\Ge_m^{-1}\,\Lstar - k\,\dev{\Ge_m}\right)}$ and 
       if the sum of the off-diagonal components of $\mathbf{\Lambda}$ is less than $\tr{\mathbf{\Lambda}}/5$ and $\delta t^m > \tau$ then
       the scheme selects the fixed point iteration \eqref{eqn.fpi}.
   \item Else if $\sqrt{\fls\tr{(\dev{\Ge_m}\,\dev{\Ge_m})}} < (\det{\Ge_m})^{1/3}/5$ or if any of the diagonal 
       entries of $\dev{\hat{\Ge}_m}$ has magnitude smaller than $\tr{\hat{\Ge}_m}/1000$ then the scheme uses \eqref{eqn.smalldevg}.
   \item If none of the above, then we apply approximation \eqref{eqn.largedevg}.
\end{enumerate}
Regardless of the chosen approximation method, at the end of each sub-timestep, the result $\Ge_{m+1} = \Ge(t^{m+1})$
must be multiplied by ${(D/\det{\Ge_{m+1}})}^{1/3}$ so that the determinant constraint is satisfied.

%% file: results.tex
\subsection{Experimental verification of the Abgrall condition}
\begin{figure}[!bp]
   \includegraphics[width=0.5\textwidth]{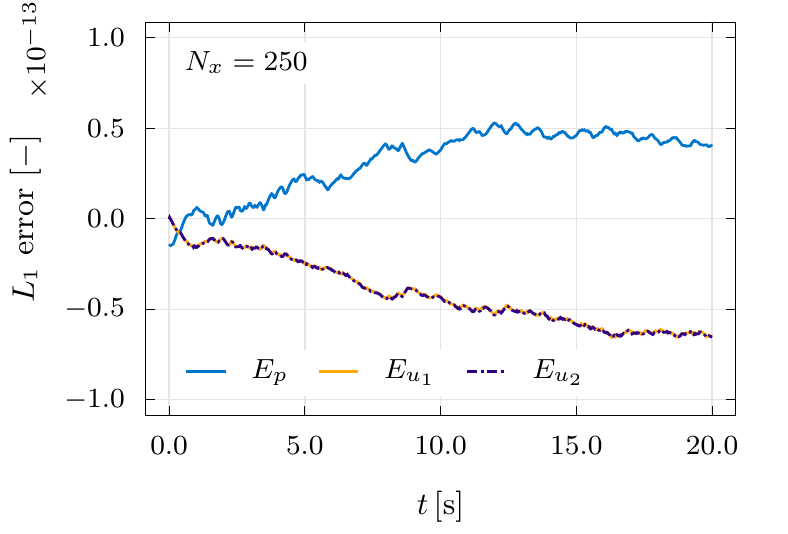}%
   \includegraphics[width=0.5\textwidth]{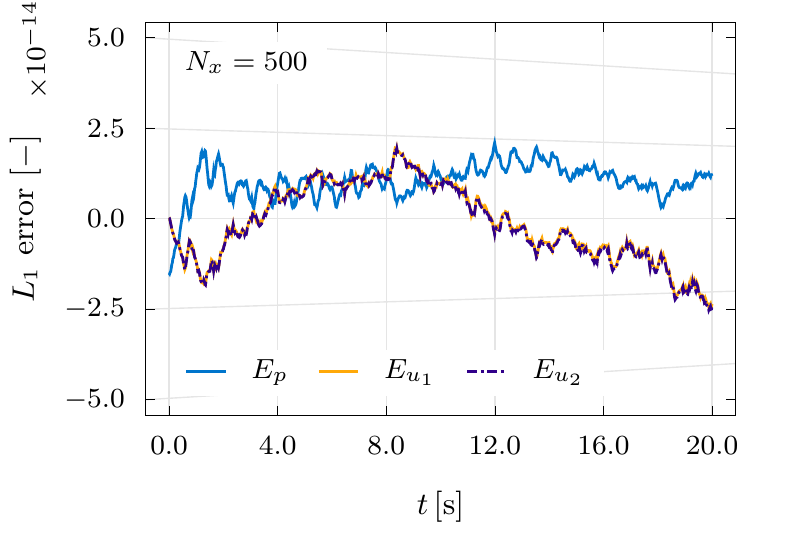}
   \caption{Timeseries of errors for the Abgrall condition verification test for the pressure field $p$, and for the velocity components $u_1$, $u_2$.
   No linear or exponential growth is observed, showing that indeed the implementation of the scheme does satisfy the Abgrall condition.} 
   \label{fig:abgrall}
\end{figure}
We begin the validation of the proposed numerical method by showing that the implementation
of the numerical scheme does indeed 
satisfy the Abgrall condition \cite{abgrallcondition}, i.e. it preserves uniformity of the velocity and pressure fields
regardless of the distribution of density or volume fraction.
The test is carried out on a domain $\Omega = [-1,\ 1]\times[-1,\ 1]$, where we set up a circular jump in $\rho_1,$ $\rho_2$, and $\alpha_1$.
Precisely, if $r = \norm{\vec{x}} < 1/2$ we set $\rho_1 = 1$, $\rho_2 = 1/2$, and $\alpha_1 = 0.3$ and otherwise if $r \geq 1/2$
we have $\rho_1 = 1/2$, $\rho_2 = 1$, and $\alpha_1 = 0.7$. The uniform pressure is $p = 1$ and a 
constant velocity field $\vec{u} = (1,\ 1,\ 0)^\transpose$ is initially imposed, so that the solution will consist in 
simple advection (and numerical diffusion) of the initial circular density and volume fraction jump.
Surface tension, shear, and gravity effects are not present, hence we set $\sigma = 0$ 
and $\vec{A} = \vec{0}$, $\vec{b} = \vec{0}$, $\vec{g} = \vec{0}$, $\tau = 10^{-14}$, $\cshear=0$, and the governing equations
reduce to Kapila's model. 
The material parameters for this test are $\gamma_1 = 4.0$, $\gamma_2 = 1.4$, $\Pi_1 = 2.0$, $\Pi_2 = 0$.
Throughout this section, unless explicitly noted, we adopt SI units for all quantities.

In Figure~\ref{fig:abgrall} we report the time evolution of errors (in the $L_1$ norm) for both nonzero components of 
the velocity field and for pressure with regard to simulations carried out on two different uniform Cartesian grids of 
$250^2$ and $500^2$ elements. For both meshes, the errors are of order $10^{-14}$ to $10^{-13}$ and present no exponential or
linear growth, thus they can perfectly be explained as accumulation of roundoff errors due to floating point arithmetic
and small errors due to the iterative solution of the discrete wave equation for the pressure field.

\begin{table}[!t]
    \caption{Numerical convergence results for the curl-preserving semi-implicit scheme applied to 
    a droplet in equilibrium under surface tension forces at rest. With $N_x$ we indicate the number of cells in one row of the Cartesian computational
    grid.}
    \label{tab:convergencetablesteady}
    \begin{tabularx}{\textwidth}{rcrrrRRR}
    \toprule
      & $N_x$  & $\mathcal{E}_{L_1}$ & $\mathcal{E}_{L_2}$ & $\mathcal{E}_{L_\infty}$  &  $\mathcal{O}_{L_1}$ & $\mathcal{O}_{L_2}$ & $\mathcal{O}_{L_\infty}$ \\
    \midrule
$\alpha_1\,\rho_1$ & $128$  &  $8.72\!\times\!10^{-3}$  &  $4.01\!\times\!10^{-3}$  &  $5.48\!\times\!10^{-3}$  &   --       &   --       &   --     \\
                   & $192$  &  $3.48\!\times\!10^{-3}$  &  $1.68\!\times\!10^{-3}$  &  $2.89\!\times\!10^{-3}$  &  $2.27$ & $2.15$ & $1.58$ \\
                   & $256$  &  $1.81\!\times\!10^{-3}$  &  $9.17\!\times\!10^{-4}$  &  $2.08\!\times\!10^{-3}$  &  $2.28$ & $2.10$ & $1.16$ \\
                   & $384$  &  $7.08\!\times\!10^{-4}$  &  $3.87\!\times\!10^{-4}$  &  $1.04\!\times\!10^{-3}$  &  $2.31$ & $2.13$ & $1.72$ \\
                   & $512$  &  $3.80\!\times\!10^{-4}$  &  $2.16\!\times\!10^{-4}$  &  $6.64\!\times\!10^{-4}$  &  $2.16$ & $2.04$ & $1.54$ \\[2mm]
$\rho\,u$ & $128$  &  $5.94\!\times\!10^{-4}$  &  $3.06\!\times\!10^{-4}$  &  $4.88\!\times\!10^{-4}$  &   --       &   --       &   --     \\
          & $192$  &  $3.72\!\times\!10^{-4}$  &  $2.04\!\times\!10^{-4}$  &  $3.59\!\times\!10^{-4}$  &  $1.16$ & $1.00$ & $0.76$ \\
          & $256$  &  $2.40\!\times\!10^{-4}$  &  $1.36\!\times\!10^{-4}$  &  $2.49\!\times\!10^{-4}$  &  $1.52$ & $1.42$ & $1.28$ \\
          & $384$  &  $1.12\!\times\!10^{-4}$  &  $6.94\!\times\!10^{-5}$  &  $1.38\!\times\!10^{-4}$  &  $1.72$ & $1.66$ & $1.46$ \\
          & $512$  &  $6.88\!\times\!10^{-5}$  &  $4.06\!\times\!10^{-5}$  &  $8.73\!\times\!10^{-5}$  &  $1.93$ & $1.87$ & $1.59$ \\[2mm]
$\rho\,E$ & $128$  &  $2.56\!\times\!10^{-3}$  &  $1.19\!\times\!10^{-3}$  &  $1.66\!\times\!10^{-3}$  &   --       &   --       &   --     \\
          & $192$  &  $1.03\!\times\!10^{-3}$  &  $4.98\!\times\!10^{-4}$  &  $8.71\!\times\!10^{-4}$  &  $2.27$ & $2.15$ & $1.58$ \\
          & $256$  &  $5.34\!\times\!10^{-4}$  &  $2.73\!\times\!10^{-4}$  &  $6.26\!\times\!10^{-4}$  &  $2.28$ & $2.10$ & $1.15$ \\
          & $384$  &  $2.10\!\times\!10^{-4}$  &  $1.16\!\times\!10^{-4}$  &  $3.13\!\times\!10^{-4}$  &  $2.30$ & $2.11$ & $1.71$ \\
          & $512$  &  $1.13\!\times\!10^{-4}$  &  $6.46\!\times\!10^{-5}$  &  $2.01\!\times\!10^{-4}$  &  $2.15$ & $2.03$ & $1.54$ \\[2mm]
$\alpha_1$  & $128$  &  $8.34\!\times\!10^{-3}$  &  $3.88\!\times\!10^{-3}$  &  $5.38\!\times\!10^{-3}$  &   --       &   --       &   --     \\
            & $192$  &  $3.33\!\times\!10^{-3}$  &  $1.62\!\times\!10^{-3}$  &  $2.83\!\times\!10^{-3}$  &  $2.26$ & $2.15$ & $1.58$ \\
            & $256$  &  $1.73\!\times\!10^{-3}$  &  $8.88\!\times\!10^{-4}$  &  $2.04\!\times\!10^{-3}$  &  $2.28$ & $2.10$ & $1.15$ \\
            & $384$  &  $6.81\!\times\!10^{-4}$  &  $3.77\!\times\!10^{-4}$  &  $1.02\!\times\!10^{-3}$  &  $2.30$ & $2.11$ & $1.71$ \\
            & $512$  &  $3.67\!\times\!10^{-4}$  &  $2.10\!\times\!10^{-4}$  &  $6.52\!\times\!10^{-4}$  &  $2.15$ & $2.03$ & $1.54$ \\[2mm]
$b_1$ & $128$  &  $7.67\!\times\!10^{-2}$  &  $3.88\!\times\!10^{-2}$  &  $5.34\!\times\!10^{-2}$  &   --       &   --       &   --     \\
      & $192$  &  $3.59\!\times\!10^{-2}$  &  $2.04\!\times\!10^{-2}$  &  $3.11\!\times\!10^{-2}$  &  $1.87$ & $1.59$ & $1.33$ \\
      & $256$  &  $2.19\!\times\!10^{-2}$  &  $1.32\!\times\!10^{-2}$  &  $2.21\!\times\!10^{-2}$  &  $1.72$ & $1.52$ & $1.19$ \\
      & $384$  &  $1.11\!\times\!10^{-2}$  &  $7.05\!\times\!10^{-3}$  &  $1.39\!\times\!10^{-2}$  &  $1.68$ & $1.54$ & $1.15$ \\
      & $512$  &  $6.51\!\times\!10^{-3}$  &  $4.38\!\times\!10^{-3}$  &  $9.96\!\times\!10^{-3}$  &  $1.86$ & $1.65$ & $1.15$ \\
    \bottomrule
    \end{tabularx}
\end{table}

\begin{table}[!t]
    \caption{Numerical convergence results for the curl-preserving semi-implicit scheme applied to 
    a droplet in equilibrium under surface tension forces moving in a uniform flow. With $N_x$ we indicate the number of cells in one row of the Cartesian computational
    grid.}
    \label{tab:convergencetableunsteady}
    \begin{tabularx}{\textwidth}{rcrrrRRR}
    \toprule
      & $N_x$  & $\mathcal{E}_{L_1}$ & $\mathcal{E}_{L_2}$ & $\mathcal{E}_{L_\infty}$  &  $\mathcal{O}_{L_1}$ & $\mathcal{O}_{L_2}$ & $\mathcal{O}_{L_\infty}$ \\
    \midrule
$\alpha_1\,\rho_1$ & $128$  &  $3.45\!\times\!10^{-2}$  &  $1.67\!\times\!10^{-2}$  &  $1.70\!\times\!10^{-2}$  &   --       &   --       &   --     \\
                   & $192$  &  $2.27\!\times\!10^{-2}$  &  $1.08\!\times\!10^{-2}$  &  $1.10\!\times\!10^{-2}$  &  $1.04$ & $1.07$ & $1.09$ \\
                   & $256$  &  $1.71\!\times\!10^{-2}$  &  $8.17\!\times\!10^{-3}$  &  $8.02\!\times\!10^{-3}$  &  $0.98$ & $0.97$ & $1.08$ \\
                   & $384$  &  $1.16\!\times\!10^{-2}$  &  $5.55\!\times\!10^{-3}$  &  $5.27\!\times\!10^{-3}$  &  $0.96$ & $0.95$ & $1.04$ \\
                   & $512$  &  $8.83\!\times\!10^{-3}$  &  $4.22\!\times\!10^{-3}$  &  $4.02\!\times\!10^{-3}$  &  $0.95$ & $0.95$ & $0.94$ \\[2mm]
$\rho\,u$ & $128$  &  $3.69\!\times\!10^{-2}$  &  $1.69\!\times\!10^{-2}$  &  $1.75\!\times\!10^{-2}$  &   --       &   --       &   --     \\
          & $192$  &  $2.44\!\times\!10^{-2}$  &  $1.11\!\times\!10^{-2}$  &  $1.13\!\times\!10^{-2}$  &  $1.02$ & $1.03$ & $1.08$ \\
          & $256$  &  $1.87\!\times\!10^{-2}$  &  $8.48\!\times\!10^{-3}$  &  $8.31\!\times\!10^{-3}$  &  $0.92$ & $0.93$ & $1.07$ \\
          & $384$  &  $1.29\!\times\!10^{-2}$  &  $5.82\!\times\!10^{-3}$  &  $5.47\!\times\!10^{-3}$  &  $0.92$ & $0.93$ & $1.03$ \\
          & $512$  &  $9.83\!\times\!10^{-3}$  &  $4.44\!\times\!10^{-3}$  &  $4.20\!\times\!10^{-3}$  &  $0.94$ & $0.94$ & $0.92$ \\[2mm]
$\rho\,E$ & $128$  &  $1.06\!\times\!10^{-2}$  &  $5.12\!\times\!10^{-3}$  &  $5.24\!\times\!10^{-3}$  &   --       &   --       &   --     \\
          & $192$  &  $6.96\!\times\!10^{-3}$  &  $3.31\!\times\!10^{-3}$  &  $3.37\!\times\!10^{-3}$  &  $1.03$ & $1.07$ & $1.09$ \\
          & $256$  &  $5.26\!\times\!10^{-3}$  &  $2.51\!\times\!10^{-3}$  &  $2.47\!\times\!10^{-3}$  &  $0.98$ & $0.97$ & $1.08$ \\
          & $384$  &  $3.57\!\times\!10^{-3}$  &  $1.70\!\times\!10^{-3}$  &  $1.62\!\times\!10^{-3}$  &  $0.96$ & $0.95$ & $1.04$ \\
          & $512$  &  $2.71\!\times\!10^{-3}$  &  $1.30\!\times\!10^{-3}$  &  $1.24\!\times\!10^{-3}$  &  $0.96$ & $0.95$ & $0.93$ \\[2mm]
$\alpha_1$ & $128$  &  $3.44\!\times\!10^{-2}$  &  $1.66\!\times\!10^{-2}$  &  $1.70\!\times\!10^{-2}$  &   --       &   --       &   --     \\
           & $192$  &  $2.26\!\times\!10^{-2}$  &  $1.08\!\times\!10^{-2}$  &  $1.09\!\times\!10^{-2}$  &  $1.03$ & $1.07$ & $1.09$ \\
           & $256$  &  $1.71\!\times\!10^{-2}$  &  $8.15\!\times\!10^{-3}$  &  $8.02\!\times\!10^{-3}$  &  $0.98$ & $0.97$ & $1.08$ \\
           & $384$  &  $1.16\!\times\!10^{-2}$  &  $5.54\!\times\!10^{-3}$  &  $5.26\!\times\!10^{-3}$  &  $0.96$ & $0.95$ & $1.04$ \\
           & $512$  &  $8.81\!\times\!10^{-3}$  &  $4.21\!\times\!10^{-3}$  &  $4.02\!\times\!10^{-3}$  &  $0.96$ & $0.95$ & $0.93$ \\[2mm]
$b_1$ & $128$  &  $3.29\!\times\!10^{-1}$  &  $1.44\!\times\!10^{-1}$  &  $1.43\!\times\!10^{-1}$  &   --       &   --       &   --     \\
      & $192$  &  $2.17\!\times\!10^{-1}$  &  $9.79\!\times\!10^{-2}$  &  $1.03\!\times\!10^{-1}$  &  $1.03$ & $0.95$ & $0.81$ \\
      & $256$  &  $1.60\!\times\!10^{-1}$  &  $7.33\!\times\!10^{-2}$  &  $7.98\!\times\!10^{-2}$  &  $1.06$ & $1.01$ & $0.88$ \\
      & $384$  &  $1.05\!\times\!10^{-1}$  &  $4.90\!\times\!10^{-2}$  &  $5.55\!\times\!10^{-2}$  &  $1.04$ & $0.99$ & $0.90$ \\
      & $512$  &  $7.70\!\times\!10^{-2}$  &  $3.65\!\times\!10^{-2}$  &  $4.34\!\times\!10^{-2}$  &  $1.08$ & $1.03$ & $0.86$ \\
    \bottomrule
    \end{tabularx}
\end{table}

\subsection{Numerical convergence study for a steady droplet in equilibrium}
A numerical convergence study is carried out in order to assess the order of accuracy of the proposed semi-implicit curl-preserving method.
The problem setup consists 
of the uniform convective transport of a diffuse droplet initialised according to the exact solution derived in \cite{hstglm}, 
to which we remand the interested reader for a detailed description of the initialisation procedure. 



Here we carry out an additional second test with a stationary droplet (the initial velocity field is 
set to $\vec{u} = \vec{0}$) in order to 
assess the convergence rates both in steady and unsteady problems.
The order is theoretically expected to be 2 for steady solutions, since all fluxes are integrated separately 
with second order accurate discretisations. In the unsteady case, due to 
the first order splitting of convection, pressure, and capillarity effects, we expect that the scheme be 
first order accurate.
Error norms and convergence rates, with respect to the variables $\alpha_1\,\rho_1$, $\rho\,\vec{u}$, $\rho\,E$, $\alpha_1$, and $b_1$, 
are given in Table~\ref{tab:convergencetablesteady} for the steady case and in Table~\ref{tab:convergencetableunsteady} for the unsteady one.
Both sets of simulations experimentally confirm the expected order of accuracy.

\subsection{Validation of the viscosity model and algorithms}
\begin{figure}[!bp]
   \includegraphics[width=0.5\textwidth]{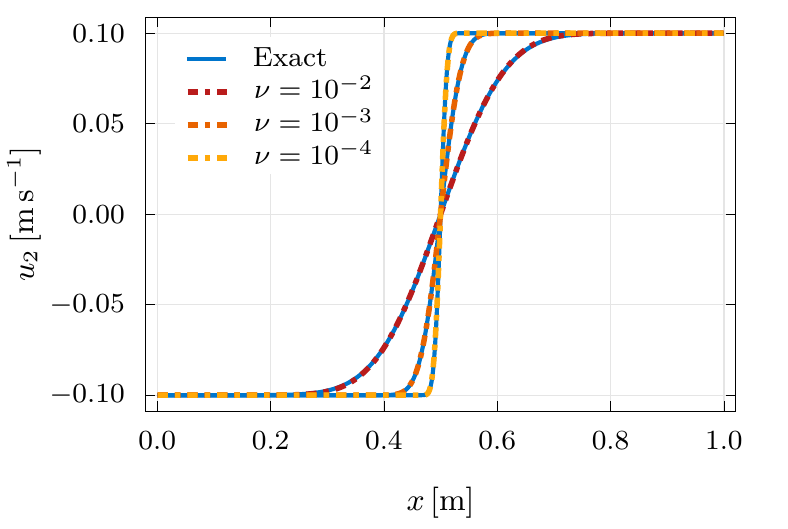}%
   \includegraphics[width=0.5\textwidth]{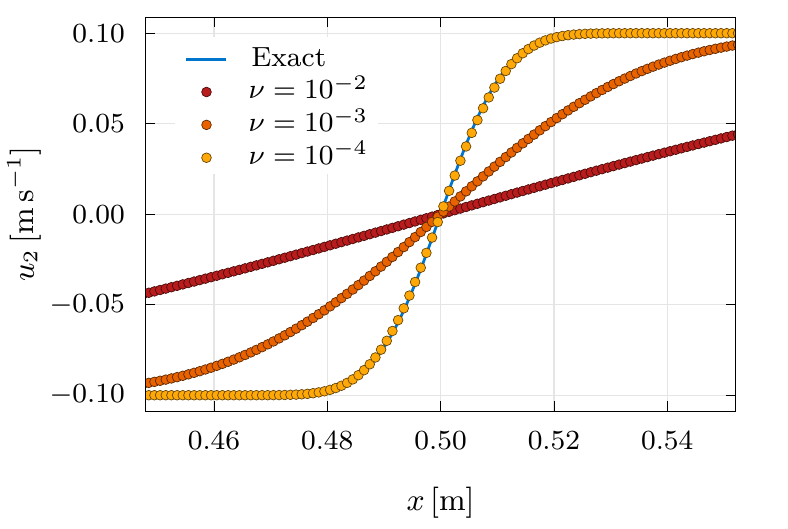}
   \caption{Numerical solutions of the first problem of Stokes, for three different values of the kinematic viscosity $\nu$. 
   The simulations are carried out with the semi-implicit structure-preserving Finite Volume scheme employing 
   the semi-analytical integrator for strain relaxation.
    The right panel is a zoomed-in view about the location of the 
    shear interface, showing perfect agreement of the GPR model with the analytical solution of the Navier--Stokes equations.}   
   \label{fig:stokes}
\end{figure}

\paragraph{First problem of Stokes.}
In the context of this paper, the test serves the dual purpose of showing that the unified model of continuum mechanics
 does indeed include the Navier--Stokes equations as a special case, and that the semi-analytical integration scheme 
 can capture the same limit.
An important benchmark for viscous flow solvers is the first problem of Stokes: 
for this test the computational domain is $\Omega = [0,\ 1]\times[0,\ 0.1]$, with periodic boundary conditions in $y$ direction and 
wall boundaries in the $x$ direction. The initial condition of the problem is given by  
uniform density $\rho=1$ and pressure $p = 1/\gamma$.
The $x$-component of the velocity field is initialized as $u_1=0$, the distortion field is initially set to 
$\vec{A}=\vec{I}$, while the $y$-velocity component $u_2$ is $u_2=-(u_2)_0$ for $x<1/2$
and $u_2=(u_2)_0$ for $x\geq 1/2$. The parametric quantities for this test are given as $(u_2)_0 = 0.1$, 
$\gamma_1 = \gamma_2=1.4$, $\Pi_1 = \Pi_2 = 0$, $\rho_0=1$, $\cshear=1$. 

The simulations are performed with the structure preserving semi-implicit finite volume
scheme on a grid composed of $500$ by $50$ square control volumes up to a final time of $t=0.4$. 
The exact solution of the incompressible Navier--Stokes equations for this test case is expressed in terms of the $y$-velocity component $u_2$ 
as  
\begin{equation}
\label{eqn.fpstokes} 
u_2(x,\ t) = (u_2)_0 \, \up{erf}\left(\frac{x}{2\,\sqrt{\nu\,t}} \right),  
\end{equation} 
with $\nu = \mu / \rho_0$.  
The test is repeated for three different values of kinematic viscosity $\nu=10^{-2}$, $\nu=10^{-3}$, 
 $\nu=10^{-4}$. 
The comparison between the Navier--Stokes reference solution \eqref{eqn.fpstokes} and the numerical results obtained with the new 
scheme for the unified model of continuum mechanics are presented in Fig. \ref{fig:stokes}, where one can observe an excellent agreement between the 
two for various kinematic viscosities $\nu$. This proves that the proposed numerical algorithm can accurately capture the 
Navier--Stokes regime of the governing equations.

\paragraph{Double shear layer problem.}
Here we solve the double shear layer test problem \cite{Bell1989,Tavelli2015,GPRmodel,sarayhtc}. 
The computational domain is $\Omega=[0,\ 1]^2$ with periodic boundary conditions everywhere. The $x$-component of the velocity 
field is initialised as  
\begin{equation}
u_1=\left\{
\begin{array}{l}
\tanh\left[\fls (y-0.25)\,\tilde{\rho} \right], \qquad \textnormal{ if } y \leq 0.5, \\[1.5mm]
\tanh\left[\fls (0.75-y)\,\tilde{\rho} \right], \qquad \textnormal{ if } y > 0.5,
\end{array}
\right.
\label{eq:DSL0} 
\end{equation} 
and the $y$-components is $u_2= \delta\,\sin(2\,\pi\,x)$, 
the density is initially uniform $\rho = \rho_0 = 1$ and the pressure is $p = 100/\gamma_1$.
Since
the test problem is usually adopted in a single-phase context, we set
the volume fraction function is $\alpha = 0.5$ throughout the computational domain to emulate the single-phase equations: both phases compute the exact same solution, 
and one simply obtains the mixture (single-phase) density by direct sum of the two partial densities.

\begin{figure}[!bp]
   \centering
   \includegraphics[width=0.845\textwidth]{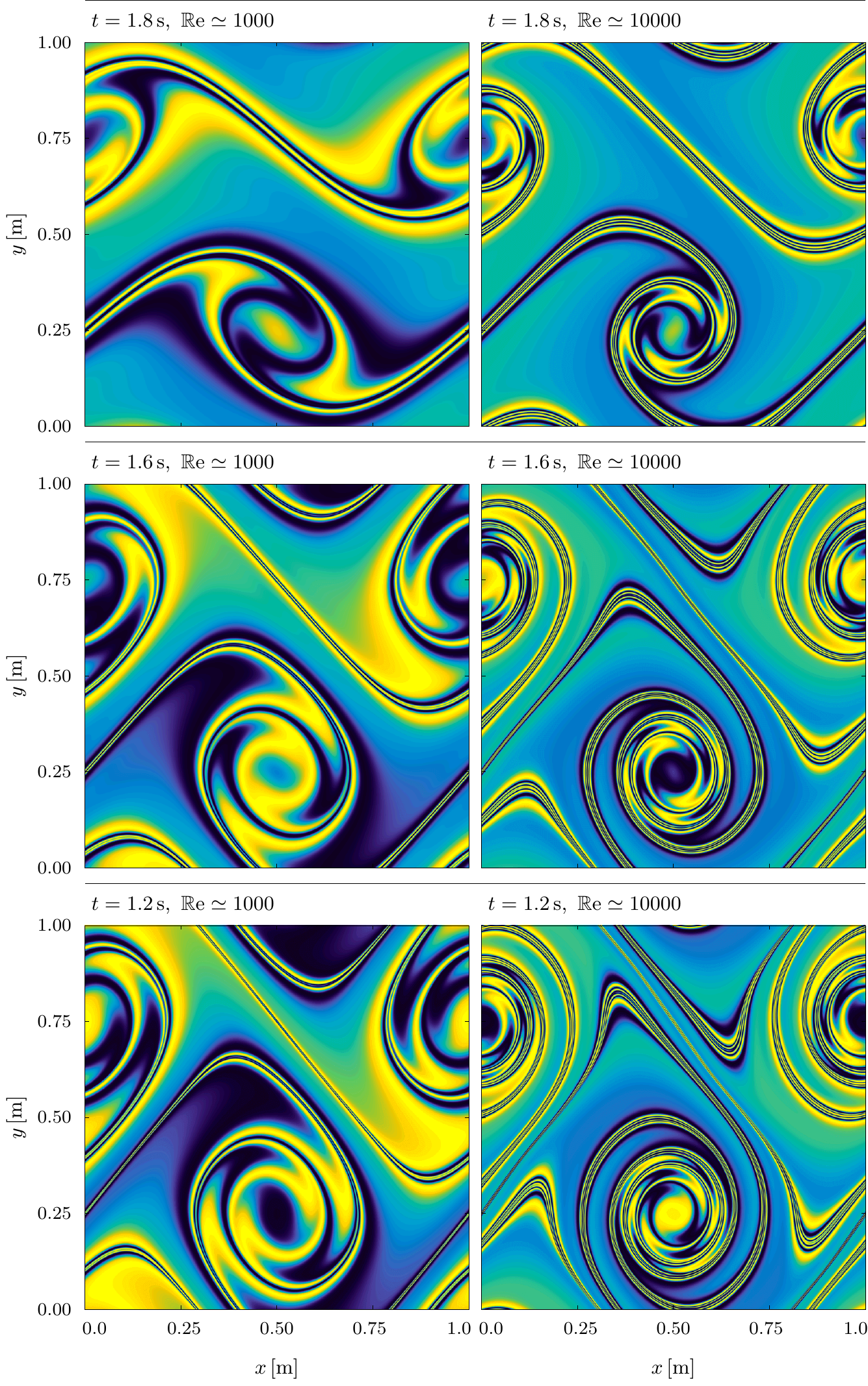}
   \caption{Filled contours of the $A_{12}$ component of the distortion field $\vec{A}$ in the double shear layer problem
   for two values of kinematic viscosity $\nu=2\!\times\!10^{-3}\,\up{m^2\,s^{-1}}$ ($\reynolds \simeq 1000$) and $\nu=2\!\times\!10^{-4}\,\up{m^2\,s^{-1}}$ 
   ($\reynolds \simeq 10000$).}   
   \label{fig:dsl}
\end{figure}
%

For this test case we set the parameters that determine the shape of the velocity field to 
$\delta=0.05$ and  $\tilde{\rho}=30$. The viscosity coefficient is set to $\nu= \mu / 
\rho_0 = 2 \times 10^{-4}$ and $\nu= \mu / 
\rho_0 = 2 \times 10^{-3}$ in two separate runs of the test problem. 
The other parameters of the model are $\gamma_1 = \gamma_2 = 1.4$, $\Pi_1 = \Pi_2 = 0$,
$\rho_0=1$, $\cshear=8$. The initial condition for the distortion field is 
$\vec{A}=\mathbf{I}$ and surface tension effects are not to be accounted for in this test, which means
that we set $\vec{b} = \vec{0}$.
Simulations are carried out with the new structure-preserving semi-implicit 
finite volume scheme up to a final time of $t=1.8$. The computational mesh is composed of 
$5120 \times 5120$ control volumes.
In Figure~\ref{fig:dsl} we show the temporal evolution of the distortion field component $A_{12}$.
The results, highlighting the incredibly rich structure found in the rotational component $\vec{R}$ of
the distortion field $\vec{A}$, are in excellent agreement with those 
obtained in \cite{sarayhtc} using a thermodynamically compatible scheme on the same fine uniform grid.

\subsection{Riemann problems and circular explosion problem} 
We continue the validation of the semi-implict scheme with two one-dimensional Riemann Problems
showing that the semi-implicit numerical method can reproduce the correct wave structure of the Kapila model \cite{kapila2001}
of two-phase flow: in these tests,
surface tension, shear, and gravity effects are not present, hence we set $\sigma = 0$ 
and $\vec{A} = \vec{I}$, $\vec{b} = \vec{0}$, $\vec{g} = \vec{0}$, $\tau = 10^{-14}$, $\cshear=0$, and the governing equations
reduce to Kapila's model exactly.
\begin{figure}[!bp]
   \includegraphics[width=0.5\textwidth]{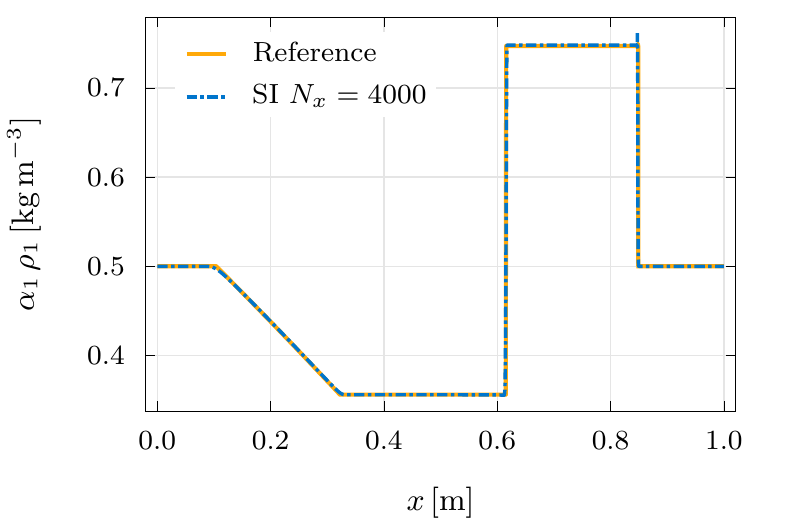}%
   \includegraphics[width=0.5\textwidth]{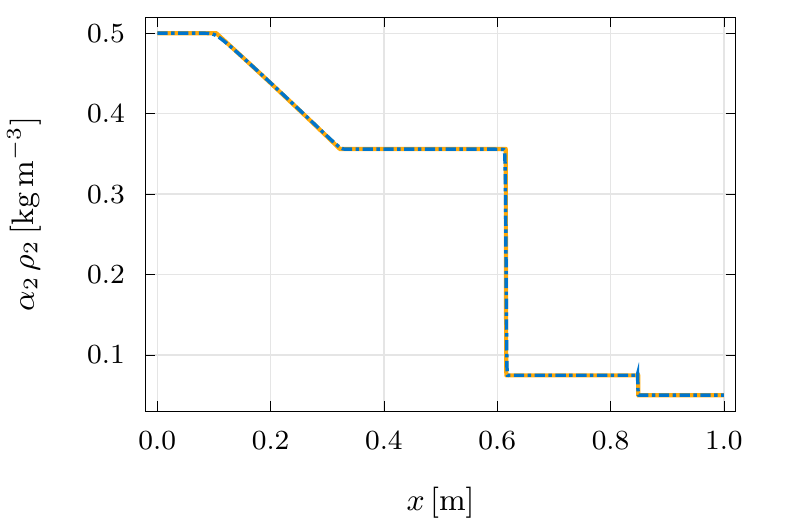}\\
   \includegraphics[width=0.5\textwidth]{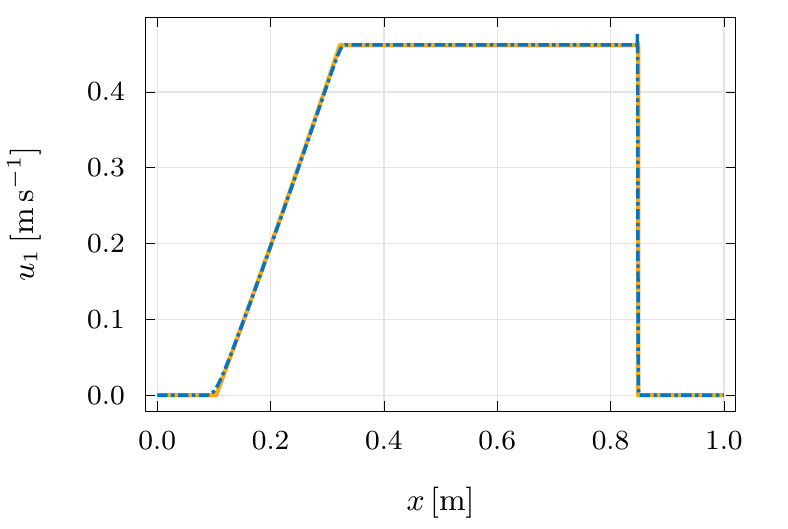}%
   \includegraphics[width=0.5\textwidth]{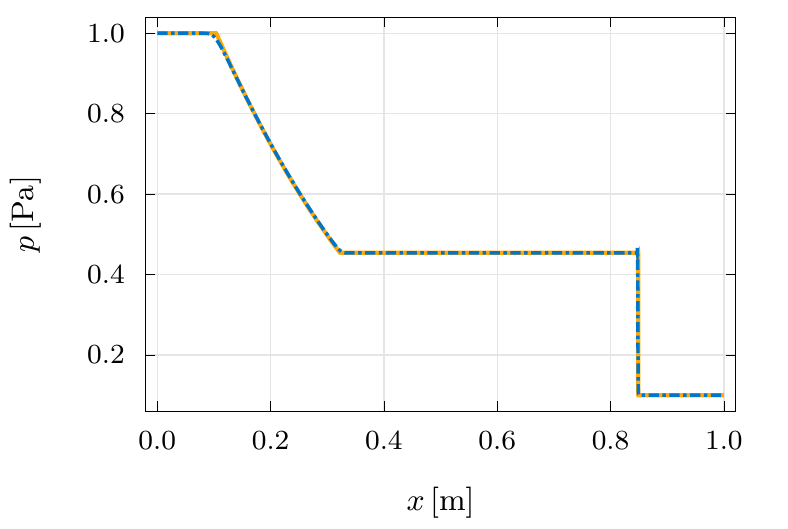}
   \caption{Numerical solution of the multiphase Riemann problem RP1 obtained with the
   semi-implict curl-preserving scheme on a uniform Cartesian grid with mesh size $\Delta x = 1/4000$, 
   compared with a reference solution obtained with a standard explicit path-conservative second order TVD MUSCL--Hancock scheme
   on a mesh of $40\,000$ cells.}   
   \label{fig:rp1}
\end{figure}
\begin{figure}[!bp]
   \includegraphics[width=0.5\textwidth]{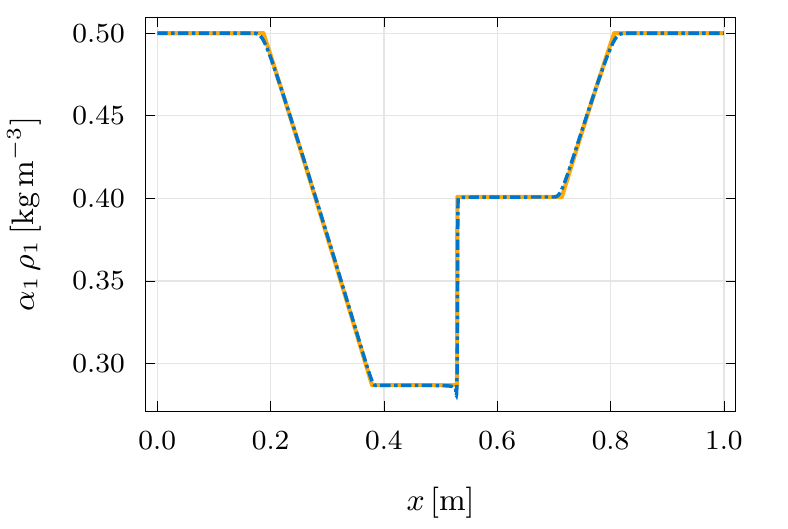}%
   \includegraphics[width=0.5\textwidth]{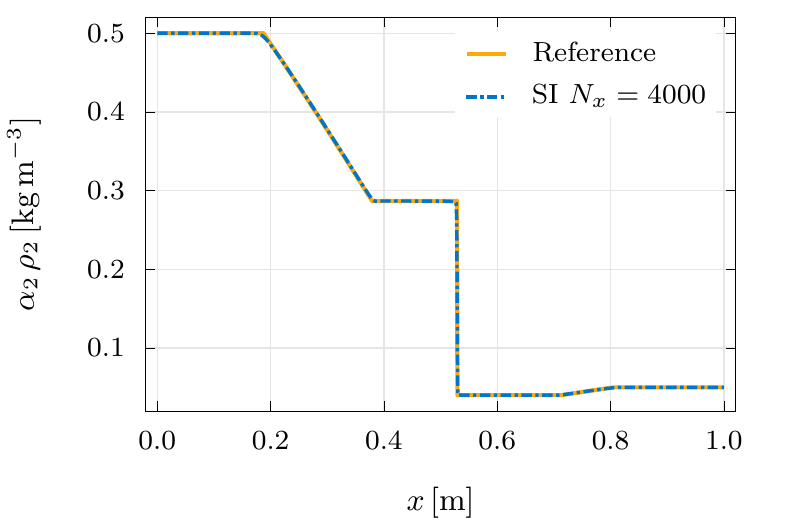}\\
   \includegraphics[width=0.5\textwidth]{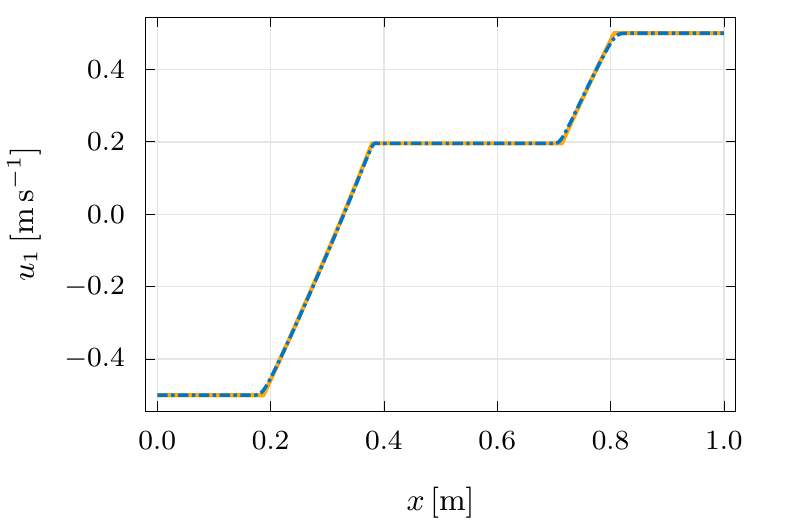}%
   \includegraphics[width=0.5\textwidth]{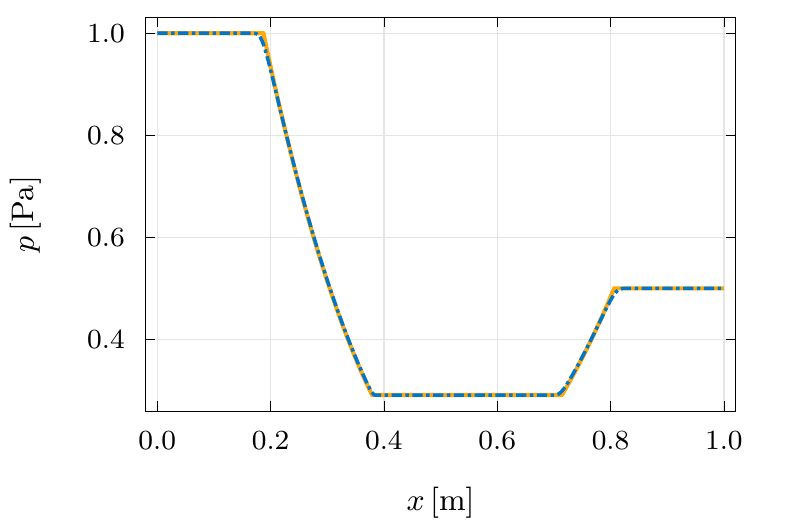}
   \caption{Numerical solution of the multiphase Riemann problem RP2 obtained with the
   semi-implicit curl-preserving scheme on a uniform Cartesian grid with mesh size $\Delta x = 1/4000$, 
   compared with a reference solution obtained with a standard explicit path-conservative second order TVD MUSCL--Hancock scheme
   on a mesh of $40\,000$ cells.}   
   \label{fig:rp2}
\end{figure}
\begin{figure}[!bp]
   \includegraphics[width=0.5\textwidth]{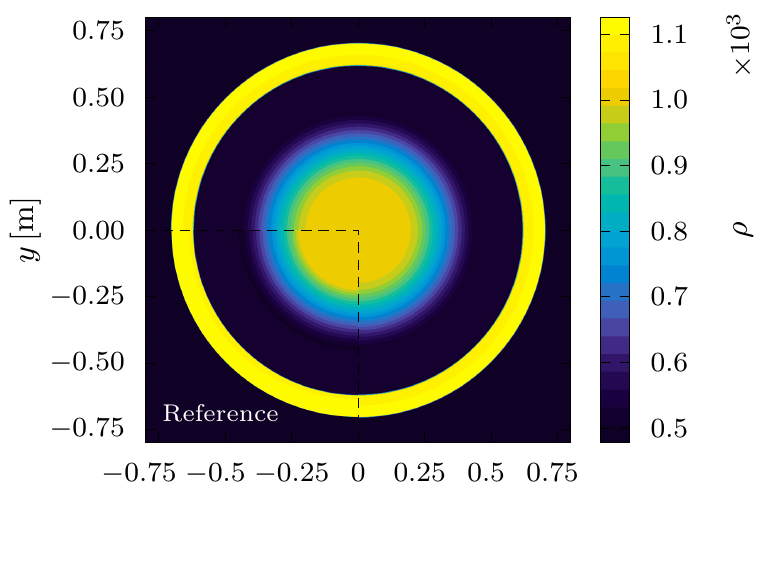}%
   \includegraphics[width=0.5\textwidth]{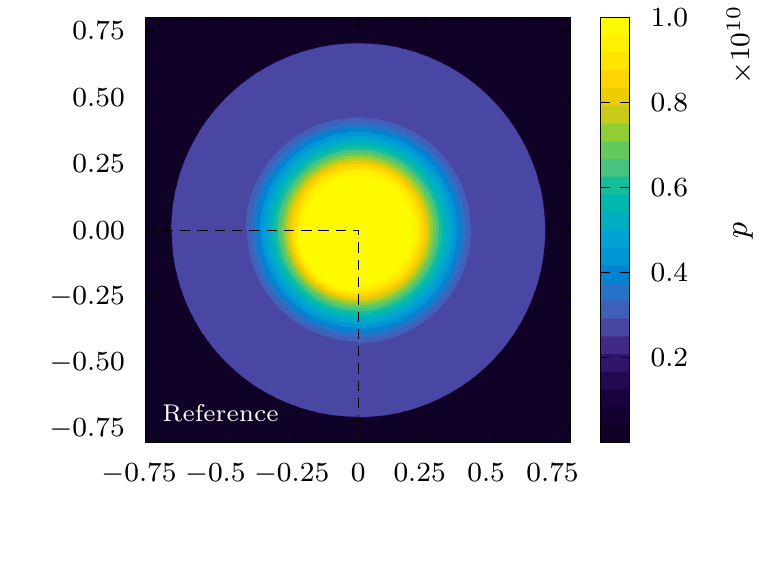}\\
   \includegraphics[width=0.5\textwidth]{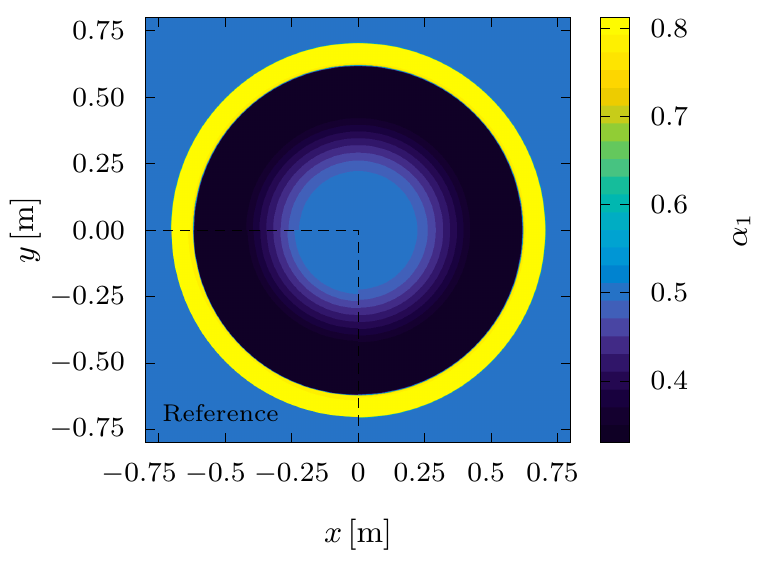}%
   \includegraphics[width=0.5\textwidth]{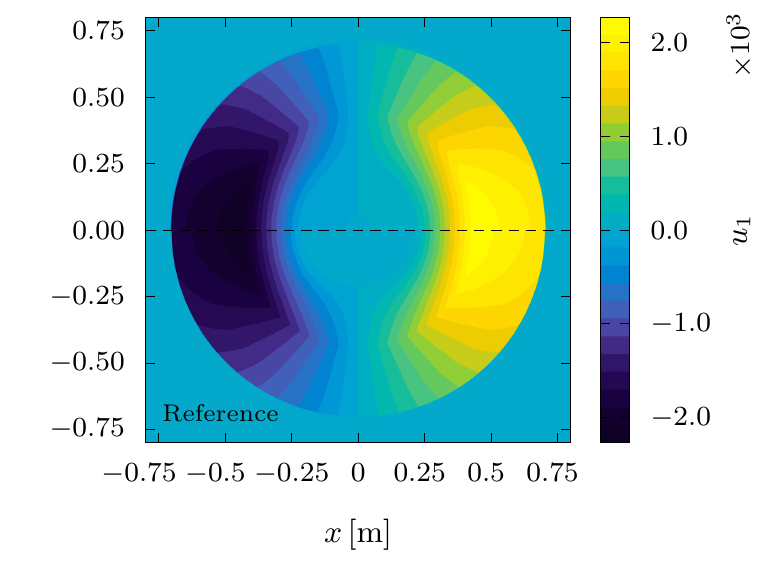}
   \caption{Filled contour plots for the two-phase inviscid circular explosion problem 
   with the curl-preserving semi-implicit scheme for viscous two-phase flow on a fine uniform 
   Cartesian grid counting $4096^2$ cells.
   In the lower half
   of each panel, we show a reference solution computed with an explicit path-conservative second order 
   method on the same fine uniform mesh of $4096^2$ elements. 
   }   
   \label{fig:explosion}
\end{figure}
First, we set up two simple planar Riemann problems RP1 and RP2 by partitioning the computational domain in 
two regions with constant state separated by a discontinuity normal to the $x$-direction, i.e. for RP1 we set 
$rho_1 = 1$, $\rho_2 = 1$, $\vec{u} = 0$, $p = 1$, $\alpha_1 = 0.5$, if $x \leq 0.5$, and 
$rho_1 = 1$, $\rho_2 = 0.1$, $\vec{u} = 0$, $p = 0.1$, $\alpha_1 = 0.5$ otherwise.

In Figures~\ref{fig:rp1} and \ref{fig:rp2}, we show the one-dimensional profiles of the partial densities $\alpha_1\,\rho_1$ and $\alpha_2\,\rho_2$, 
of the $x$-component of the velocity field $u_1$, and of the pressure $p$, obtained by applying 
the semi-implicit scheme on a mesh of $4000$ by $400$ square cells over the domain $\Omega=[0,\ 1]\times[0,\ 0.1]$.
The results are compared with a reference solutions computed by a second order explicit path-conservative TVD MUSCL--Hancock method on 
a fine one-dimensional grid compsed of $40\,000$ uniform control volumes.
The results match up to some minor artifacts observed in the solution obtained from the semi-implicit
method, in this regard, it is important to note that such a semi-implicit scheme is designed 
specifically for low Mach number flows, rather than for the treatment of shockwaves. The 
observed artifacts are not a hinderance for large scale simulations lying within the design regime of the method.

Then, a two-dimensional explosion problem is set up in a similar fashion, with the discontinuity 
now representing an inner and outer state rather than a left and right one.
The computational domain is the square $\Omega = [-0.8,\ 0.8]\times[-0.8,\ 0.8]$ and is disctretised with 
a mesh of $4096^2$ uniform Cartesian control volumes.
The inner state is given by $\rho_1^\up{i} = 1000$, $\rho_2^\up{i} = 1000$, $\vec{u}^\up{i} = \vec{0}$, $p^\up{i} = 10^{10}$, $\alpha_1^\up{i} = 0.5$.
The outer state is $\rho_1^\up{o} = 1000$, $\rho_2^\up{o} = 1$, $\vec{u}^\up{o} = \vec{0}$, $p^\up{o} = 10^{5}$, $\alpha_1^\up{0} = 0.5$.
The discontinuity is initially located at $r = \norm{\vec{x}} = 1/2$, and the final simulation time is $t_\up{end} = 6\times10^{-5}$.
The parameters of the stiffened gas equation of state are $\gamma_1 = 4.4$, $\gamma_2 = 1.4$, $\Pi_1 = 6\times10^8$, $\Pi_2 = 0$.

In Figure~\ref{fig:explosion} we show the filled contour plots for the mixture density $\rho$, for the pressure $p$, for the volume fraction $\alpha_1$, 
and for the $x$-component of the velocity field $u_1$. The results correctly preserve symmetries and are in agreement with the reference 
solutions obtained by means of an explicit path-conservative second order TVD MUSCL--Hancock method on the same mesh of $4096^2$ cells.

\subsection{Long-time evolution of an oscillating droplet at low Mach number}
In this Section, we reproduce the numerical experiments 
shown in \cite{hstglm} with regard to high order ADER Discontinuous Galerkin \pnpm{N}{N}
schemes with \aposteriori\ Finite Volume subcell limiting, now using instead the novel semi-implicit
curl-preserving scheme introduced in this work.
\begin{figure}[!bp]
   \includegraphics[width=0.5\textwidth]{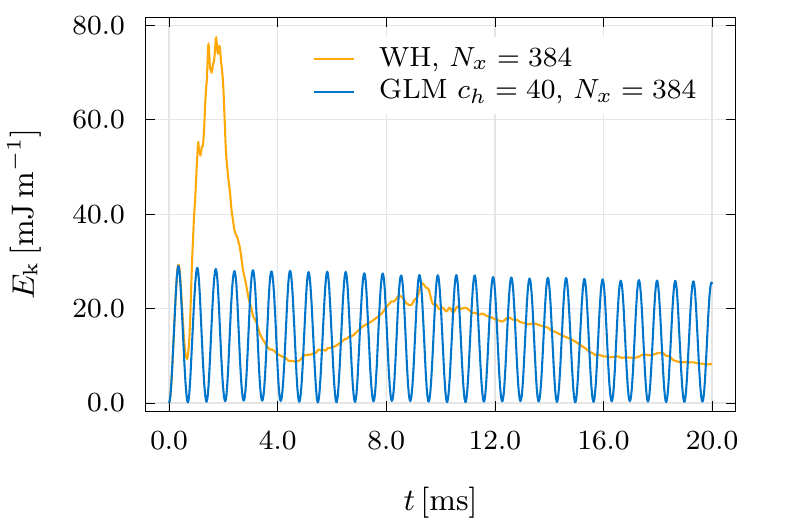}%
   \includegraphics[width=0.5\textwidth]{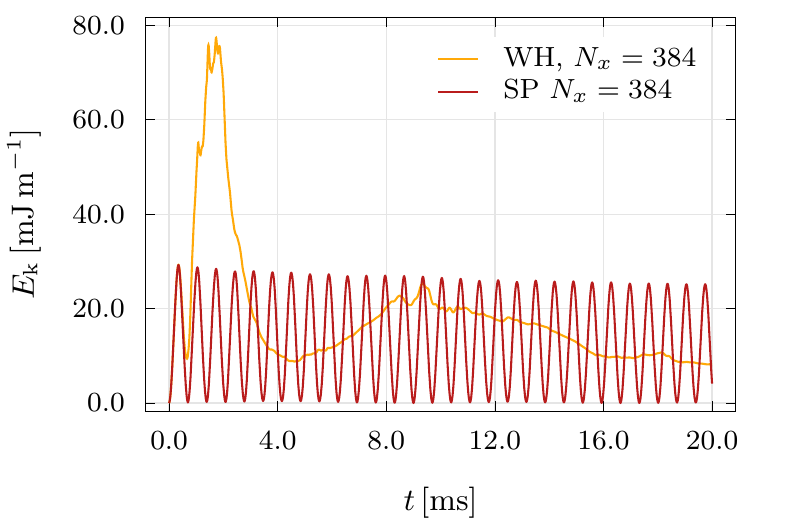}
   \caption{Time evolution of the total kinetic energy $E_\up{k}$ for the elliptical droplet oscillation problem.
   On the left: solution obtained with the a second order semi-implicit Finite Volume scheme with GLM curl cleaning 
   on a uniform Cartesian mesh composed of $384\times384$ elements. 
   On the right: solution obtained on the same mesh with the curl-preserving variant of the same semi-implicit scheme.
   In both cases it is apparent that the numerical diffusion is much lower than that of standard second
   order TVD schemes (see \cite{hstglm}): the oscillation amplitude shows only mild decay for the full 15 oscillation 
   cycles simulated in the numerical experiment. For comparison, without involution enforcement, the weakly hyperbolic model (WH)
   produces chaotic results before the end of the first oscillation cycle.}   
   \label{fig:dropletellipse.ek}
\end{figure}
\begin{figure}[!bp]
   \includegraphics[width=0.5\textwidth]{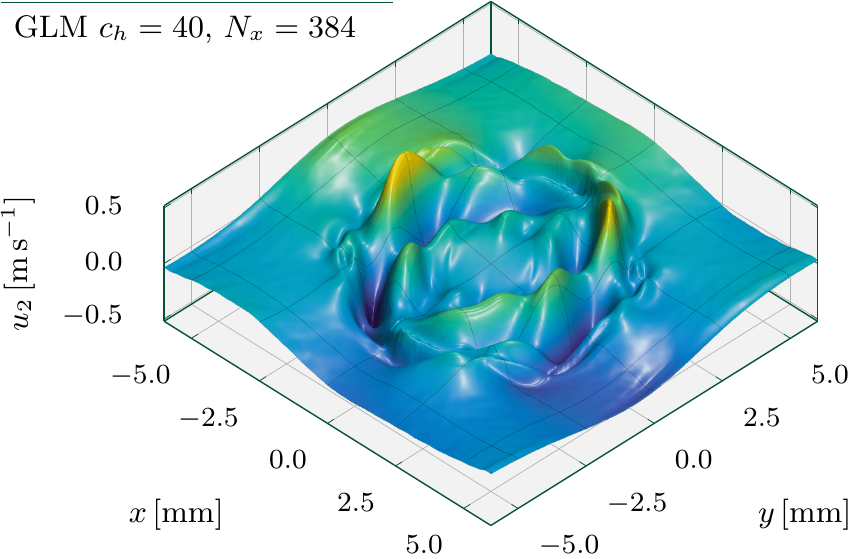}%
   \includegraphics[width=0.5\textwidth]{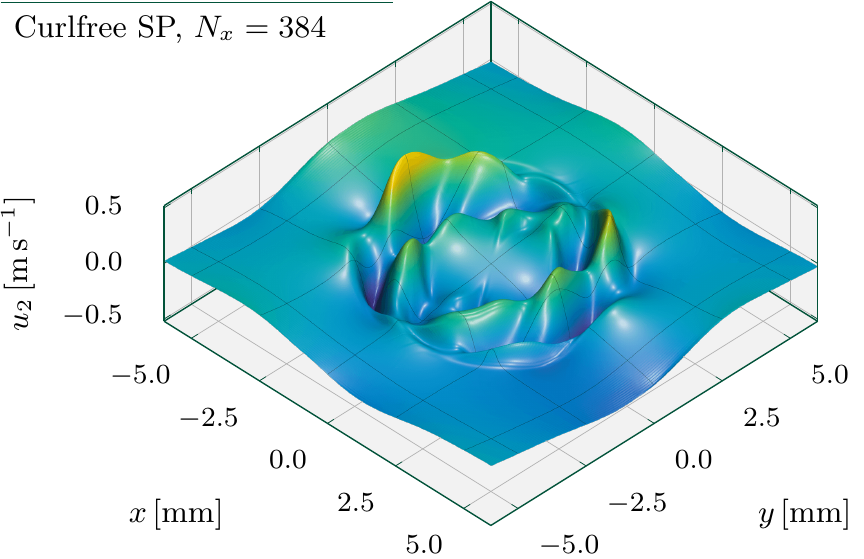}\\[6mm]
   \includegraphics[width=0.5\textwidth]{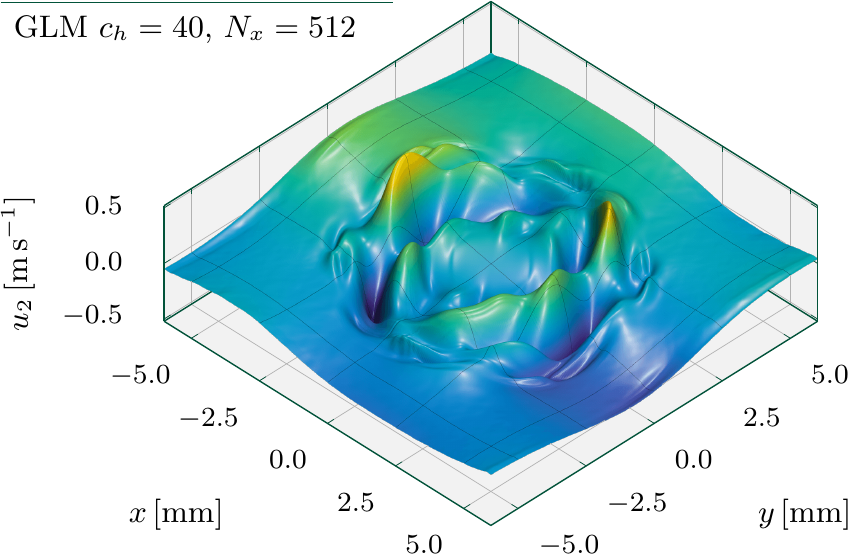}%
   \includegraphics[width=0.5\textwidth]{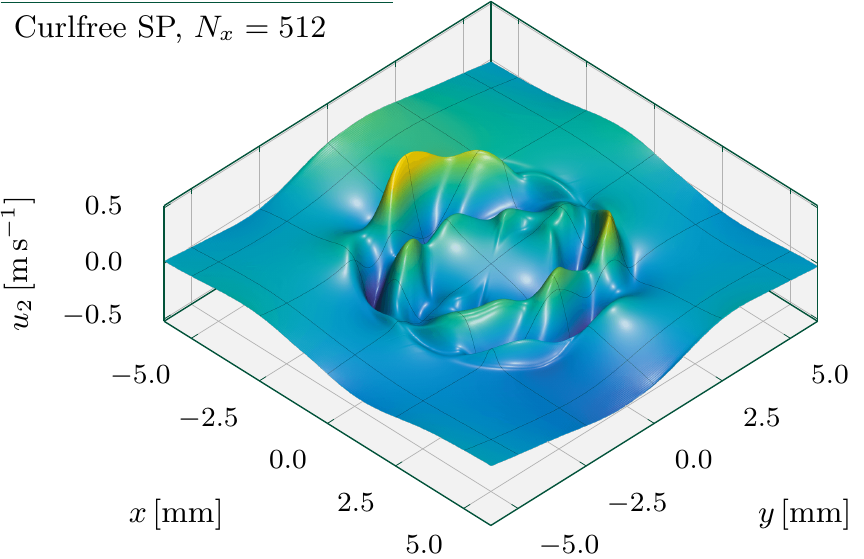}
   \caption{Three-dimensional view of the second component of the velocity field $u_2$ for the 
   elliptical droplet oscillation problem at time $t = 6.91\times\!10^{-4}\,\up{s}$ (half of the first oscillation cycle).
   In the top row the uniform Cartesian mesh counts $384^2$ elements, in the bottom row the mesh is refined to $512^2$ elements.
   On the left, curl constraints are enforced with GLM curl-cleaning and the equations solved 
   with the semi-implicit Finite Volume scheme presented in this work (without curlfree discretisation of the interface field), 
   while on the right the results are computed with the curlfree structure-preserving semi-implicit method.
   The difference in the quality of the results obtained with the two different approaches is immediately apparent by
   comparing the \textit{specular highlights} in the plots on the left column (from GLM curl cleaning, rough surface), with 
   those on the right column (from the exactly curlfree scheme, smooth surface). Moreover, note that the smoothness of the curlfree solution
   is not due to numerical diffusion (peak heights are the same) and that the observations are the same for different meshes.}   
   \label{fig:3d}
\end{figure}
\begin{figure}[!bp]
   \includegraphics[width=0.5\textwidth]{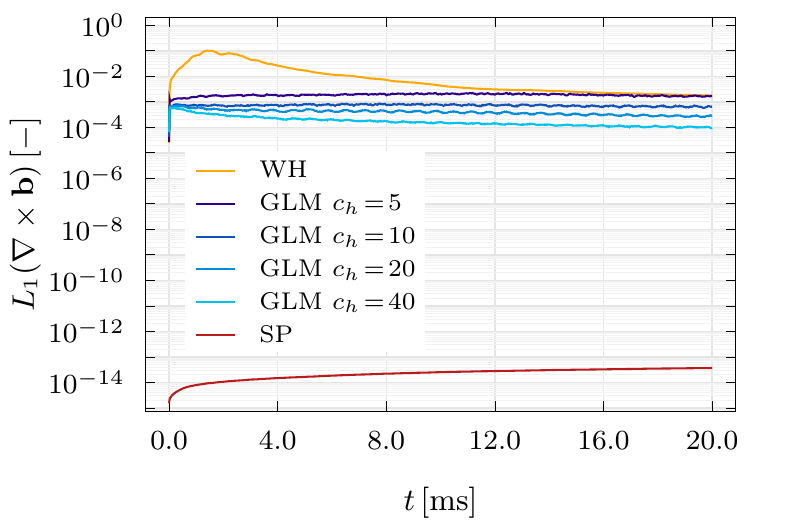}%
   \includegraphics[width=0.5\textwidth]{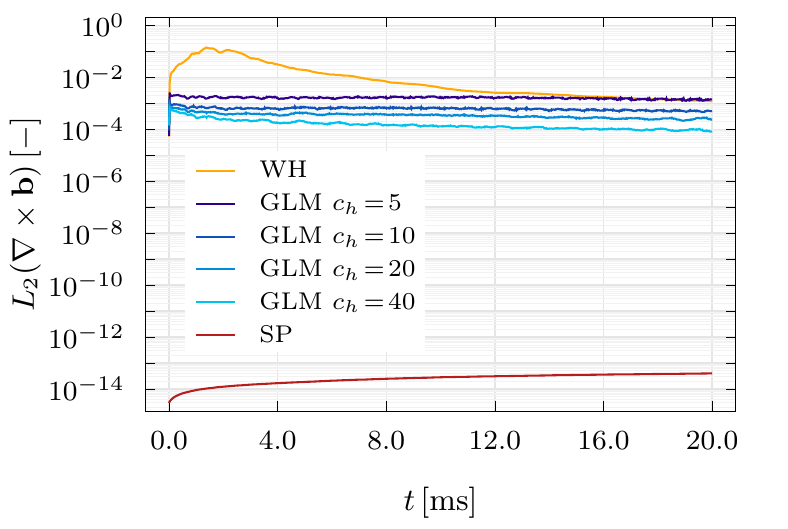}\\
   \includegraphics[width=0.5\textwidth]{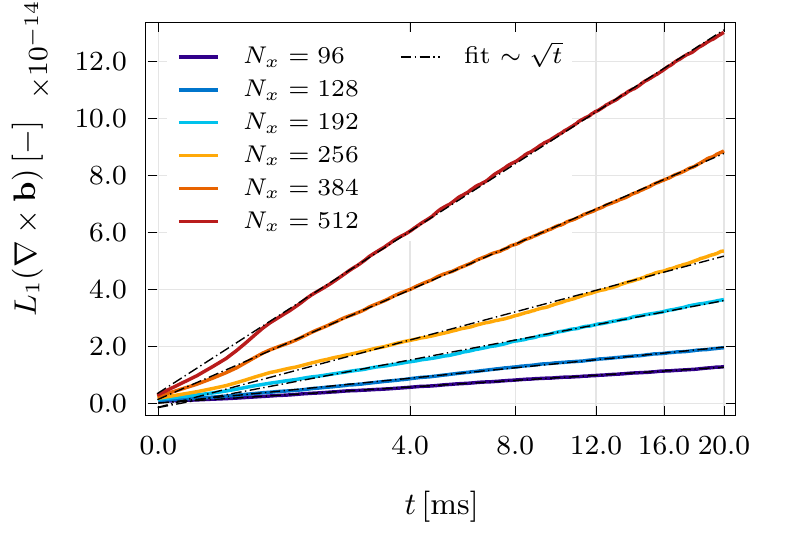}%
   \includegraphics[width=0.5\textwidth]{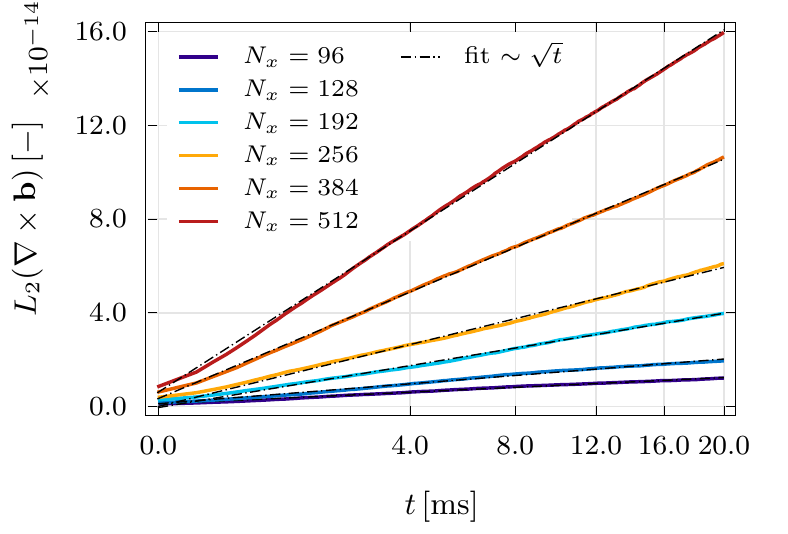}
   \caption{Time evolution of curl errors for the oscillating droplet test problem.
   In each panel of the top row, we compare the curl errors given by the GLM curl cleaning variant of the model
   with cleaning speed $c_h$ ranging from $5\,\up{m\,s^{-1}}$ to $40\,\up{m\,s^{-1}}$, and those given by
   the curl-preserving semi-implicit method.
   The errors are computed in the $L_1$ (left panel) and the $L_2$ (right panel) norms, and are 
   relative to a common uniform Cartesian mesh counting $192^2$ elements.
   The curl-preserving semi-implicit method produces errors more than nine orders of magnitude smaller
   than those obtained with GLM curl cleaning in the same semi-implicit framework.
   If involutions are not enforced (WH) curl errors grow uncontrollably and appear to decrease but the
   effect is only due to the complete breakdown of the simulation.
   In the bottom row, we show, for a variety of uniform Cartesian meshes, that the errors of the 
   curl-preserving semi-implicit Finite Volume scheme are due to statistical accumulation of roundoff
   errors: as expected from random-walk/Brownian \cite{browneinstein,feynmanbook} processes, they grow in time proportionally to $\sqrt{t}$, 
   and they grow faster for finer meshes which require more floating point operations in order to 
   integrate the solution up to a given time.
   }   
   \label{fig:dropletellipse.curl}
\end{figure}
\begin{figure}[!bp]
   \centering
   \includegraphics[width=0.890\textwidth]{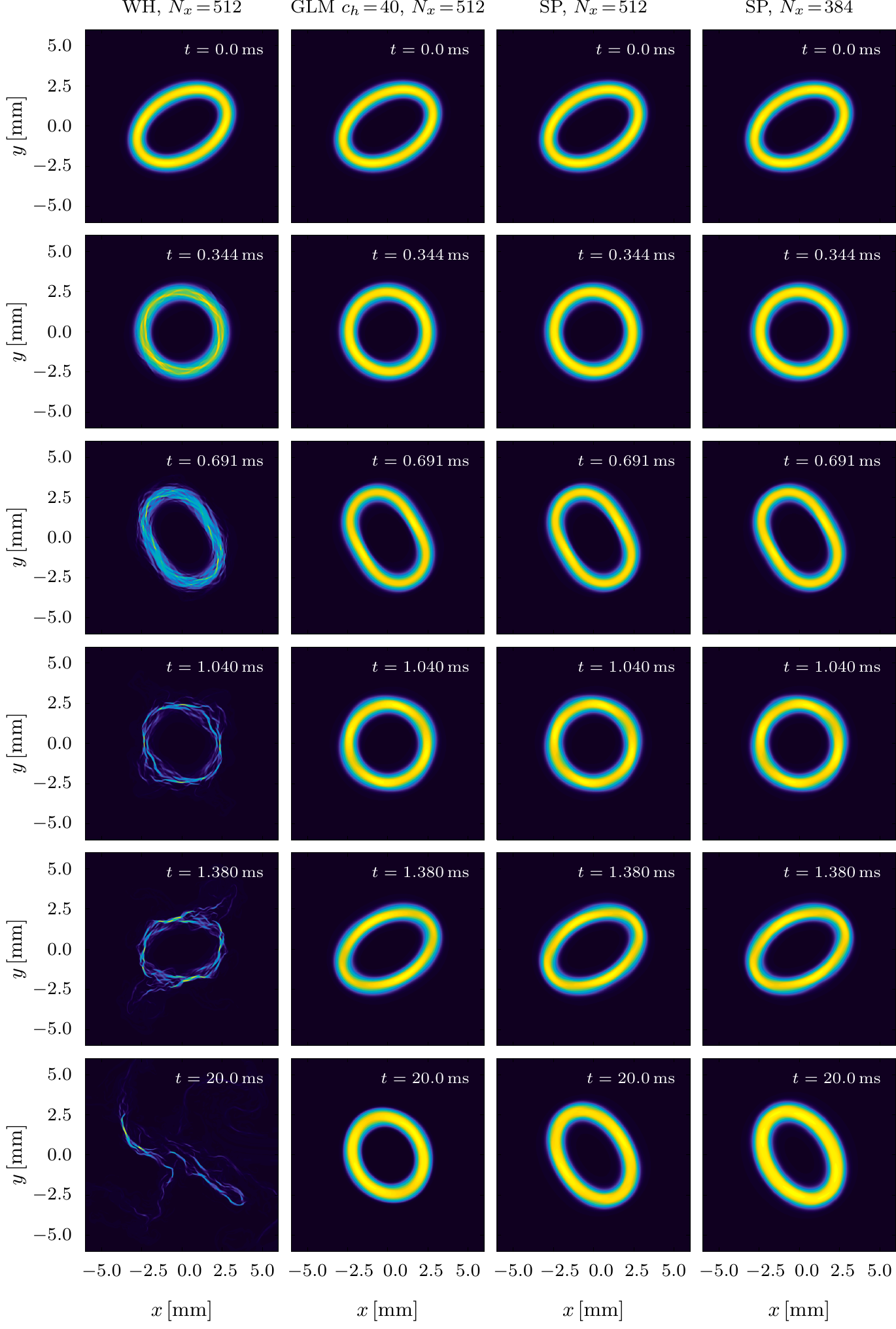}%
   \caption{Snapshots of the interface energy $\sigma\,\norm{\vec{b}}$ for the droplet oscillation problem with different 
   mesh resolutions and for different schemes. WH is the unconstrained weakly hyperbolic model, GLM refers to the
   curl cleaning method, SP indicates the curl-free semi-implicit scheme.}   
   \label{fig:dropletellipse.film}
\end{figure}

While we refer to \cite{hstglm} for a detailed exposition of 
the computational setup, we list here the parameters of the problem:
the density fields are set to the uniform values $\rho_1^0$ and $\rho_2^0$ throughout the
computational domain, as is the velocity field for which we set $\vec{u} = {(0,\ 0,\ 0)}^\transpose$.
The numerical values employed for this test problem are: {$\rho_1^0 = 1000\,\up{kg\,m^{-3}}$},
{$\rho_2^0 = 1\,\up{kg\,m^{-3}}$}, {$p_\up{atm} = 100\,\up{kPa}$}, {$R_x = 3\,\up{mm}$}, {$R_y =
2\,\up{mm}$}, {$\alpha_\up{min} = 0.01$}, {$\alpha_\up{max} = 0.99$}, {$\sigma =
60\,\up{N\,m^{-1}}$}. The parameters for the stiffened gas equation of state are: {$\Pi_1 =
1\,\up{MPa}$}, {$\Pi_2 = 0$}, {$\gamma_1 = 4$}, {$\gamma_2 = 1.4$}. The domain is the square $\Omega
= [-6\,\up{mm},\ 6\,\up{mm}]\times[-6\,\up{mm},\ 6\,\up{mm}]$ and additionally, the initial
condition is rotated counter-clockwise by 30 degrees, in order to avoid mesh alignment. In a first
batch of tests, we set {$\epsilon = 0.5\,\up{mm}$} and discretise the computational domain with
$384^2$ square cells, which yields the same number of degrees of freedom previously employed for the same
test with the ADER-DG \pnpm{5}{5} scheme in \cite{hstglm}.
These simulations are intended to test the capability of the method in a dynamical setting
where the interface deforms significantly under the effect of strong surface tension, and show
that the curl-preserving semi-implicit scheme, when applied to such low Mach problems, 
can be competitive even with very high order methods and in particular that 
that the GLM curl cleaning approach can deal with the violations of involution constraints that such
deformations generate.

In Figure~\ref{fig:dropletellipse.ek}, we show the time evolution of the total kinetic energy and compare
the results obtained without any enforcement of curl involutions (which we recall, as shown in \cite{hstglm}, are 
totally unstable and quickly lead to the breakdown of computations), with those from 
GLM curl cleaning formulation of the governing equations (this time implemented in the semi-implicit code) 
with the results from the structure preserving methodology. Both curl-control solutions yield stable results, but higher
accuracy is achieved with the exactly curlfree scheme. 

The same can be observed in Figure~\ref{fig:3d}, 
where we plot a 
three-dimensional view of the second component of the velocity field $u_2$ for the 
at time $t = 6.91\times\!10^{-4}\,\up{s}$ (half of the first oscillation cycle), computed on two different meshes of size $384^2$ and $512^2$, 
and
for both the GLM variant of the method and for the exactly curlfree method.
The difference in the quality of the results obtained with the two different approaches is immediately apparent by
comparing the \textit{specular highlights} in the plots on the left column, obtained with GLM curl cleaning (the surface is rough), with 
those on the right column, obtained with the exactly curlfree scheme (the surface is perfectly smooth). Moreover, note that the smoothness of the curlfree solution
is not due to numerical diffusion (peak heights are the same) and that the observations are the same for different meshes.

In Figure~\ref{fig:dropletellipse.curl} we show the time evolution of curl errors.
Specifically, we compare the curl errors given by the GLM curl cleaning variant of the model
with cleaning speed $c_h$ ranging from $5\,\up{m\,s^{-1}}$ to $40\,\up{m\,s^{-1}}$, and those given by
the curl-preserving semi-implicit method.
The errors are computed in the $L_1$ and the $L_2$ norms and are 
relative to a common uniform Cartesian mesh counting $192^2$ elements.
The curl-preserving semi-implicit method produces errors more than nine orders of magnitude smaller (less than $\sim10^{-13}$ instead of more than $\sim10^{-4}$)
than those obtained with GLM curl cleaning in the same semi-implicit framework.
If involutions are not enforced (WH) curl errors grow uncontrollably and appear to decrease but the
effect is only due to the complete breakdown of the simulation and the disappearance of \textit{all} physical flow features.
Additionally, in order to argue that the curl errors produced by the exactly curlfree method originate
from statistical accumulation of roundoff error, we plot their time evolution as a function of the square root 
of time $\sqrt{t}$ (essentially proportional to the square root of the number of timesteps as well).
This scaling law shows that curl errors are growing in a perfectly linear fashion as a function of $\sqrt{t}$, 
as expected from a process of random-walk, or Brownian origin \cite{browneinstein,feynmanbook}.
Moreover the growth rate is higher on finer meshes, since more operations are required to integrate
the solution up to a given time $t$. Nonetheless, this does not constitute an issue since the order of magnitude of the errors 
is that of machine-epsilon accumulated roundoff.

Finally, in Figure~\ref{fig:dropletellipse.film} we collect several snapshots of the global dynamics of the droplet oscillation problem, 
which clearly show that both the exactly curlfree methodology, and the GLM curl-cleaning method yield stable computations, 
counter to the unconstrained (weakly hyperbolic) variant of the model for hyperbolic surface tension, which catastrophically deteriorates
within short times.

\subsection{Binary droplet collision with high density ratio}
\begin{figure}[!bp]
   \centering
   \includegraphics[width=0.75\textwidth]{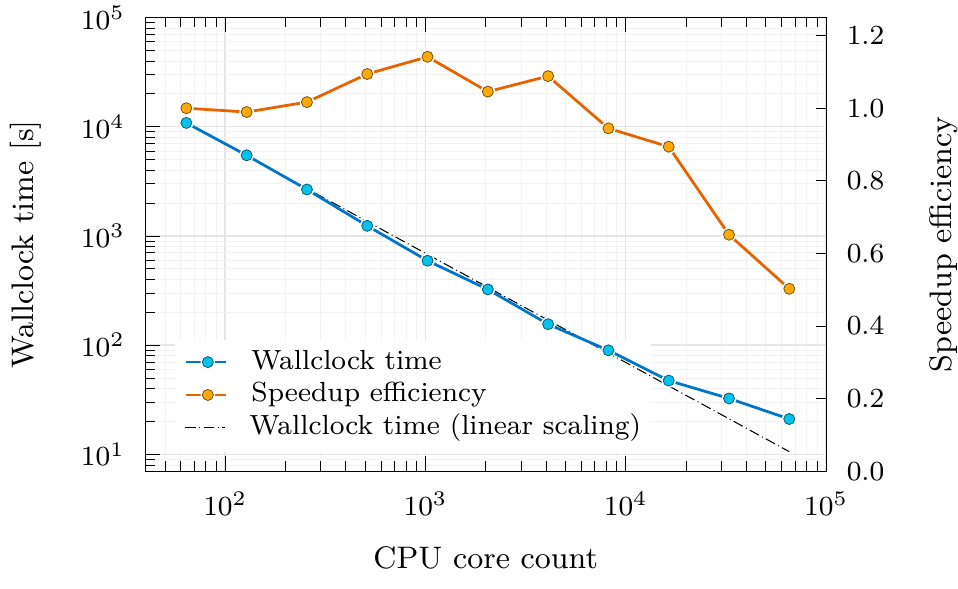}
   \caption{Strong scaling results (computational times and speedup efficiency) from 
   64 CPU cores to 65\,536 CPU cores (512 nodes) of the
   HPE-Hawk supercomputer at HLRS in Stuttgart. The semi-implicit structure-preserving scheme for 
   hyperbolic viscous two-phase flow with surface tension achieves excellent scaling performance 
   ($>\sim95\%$)
   up to about 16k CPU cores, also thanks to cache effects mitigating the main bottleneck of the
   scheme which is the solution of the pressure system via matrix-free conjugate gradient method. 
   On 65k CPU cores, the strong scaling efficiency drops to about $50\%$.
   The computational grid is composed of 
   33.5 million Cartesian cells.}   
   \label{fig:scaling}
\end{figure}
\begin{figure}[!bp]
   \includegraphics[width=\textwidth]{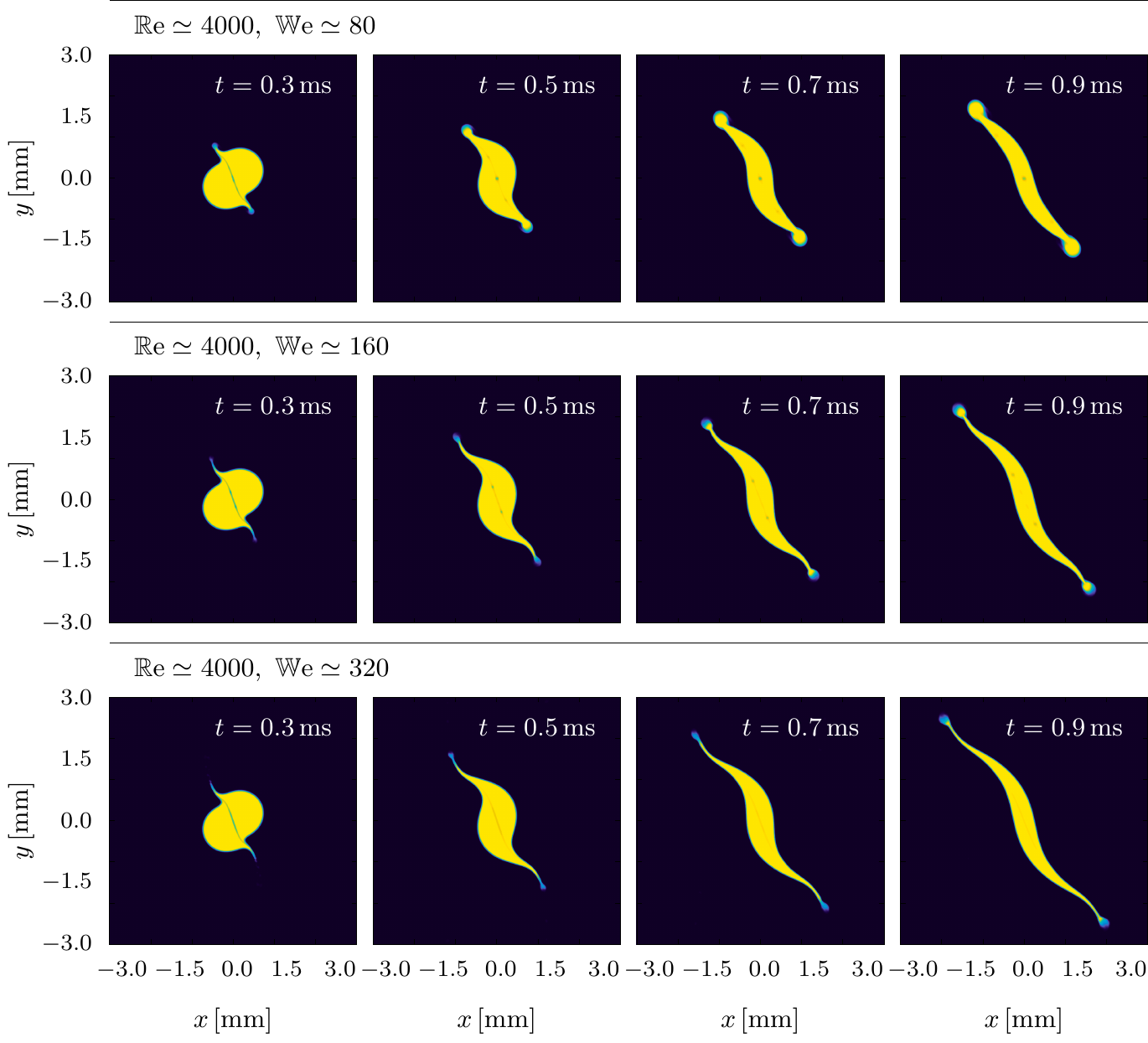}
   \caption{Early stages of the binary droplet collision simulation for three different values of the Weber number.
   The effects of surface tension penalising the formation of sharp features is evident: the lower the Weber number (corresponding to higher strength of capillarity forces), 
   the stronger the tendency towards minimizing the interface curvature.}   
   \label{fig:collision1}
\end{figure}
\begin{figure}[!bp]
   \includegraphics[width=\textwidth]{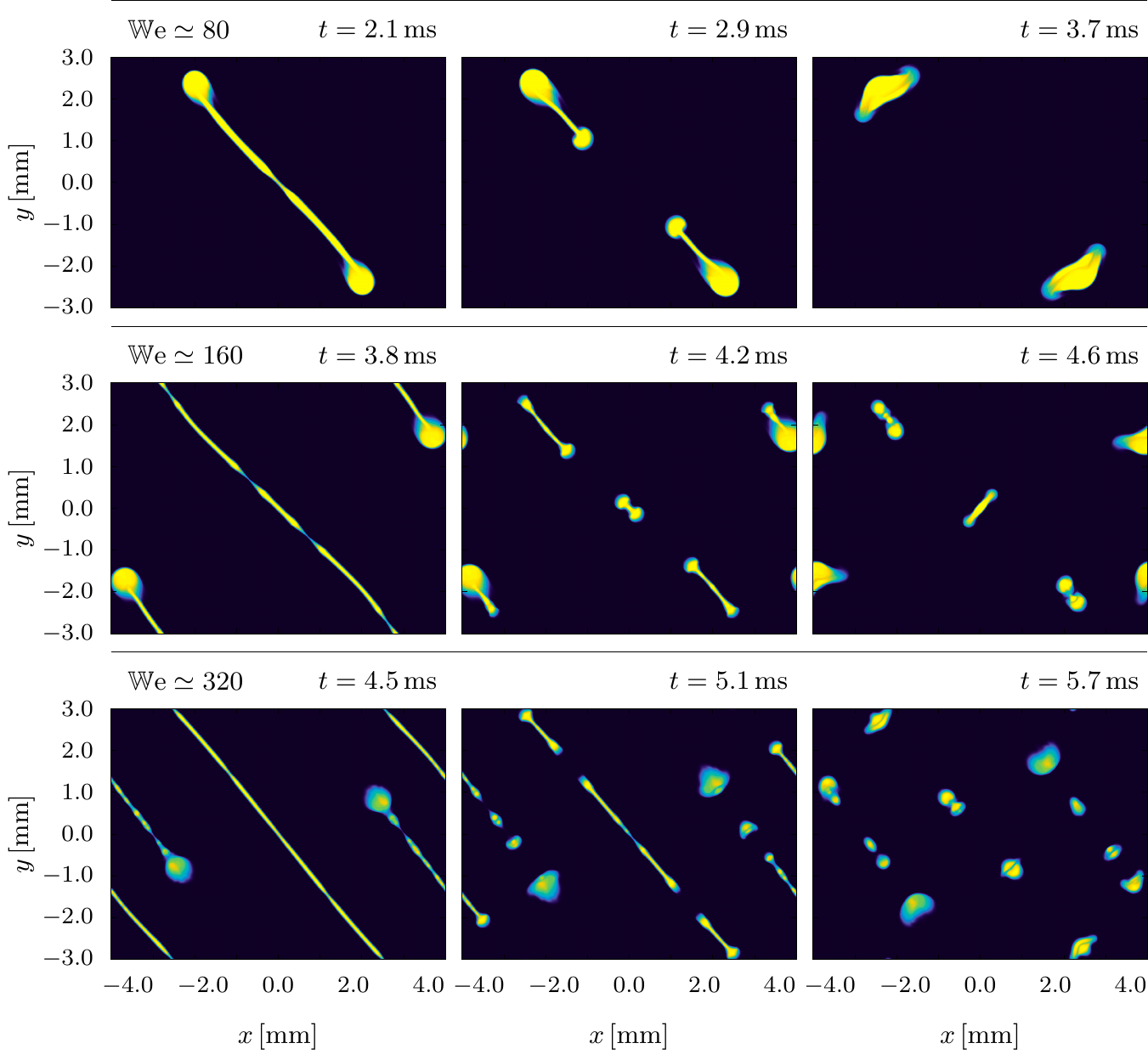}
   \caption{Late stages of the binary droplet collision simulation for three different values of the Weber number.
   At low Weber number ($\weber \simeq 80$) the collision results in stretching and separation of the two droplets. At
   all values of Weber number the effects of Rayleigh--Plateau instability are evident and in particular for the higher Weber numbers 
   they result in the formation of multiple small secondary droplets.}   
   \label{fig:collision2}
\end{figure}
Next, we continue our validation of the semi-implicit curl-preserving scheme for viscous two-phase flow with surface tension 
with an application to binary droplet collisions, motivated by \cite{brenn1, brenn2}.

Two circular droplets are initialised according to the equilibrium solution derived in \cite{hstglm}, 
with 
the centers of the droplets being 
\begin{equation}
\begin{aligned}
   &\vec{x}_\up{cl} = \left(-8\!\times\!10^{-4},\ -2\!\times\!10^{-4},\ 0\right)^\transpose,\\
   &\vec{x}_\up{cr} = \left(+8\!\times\!10^{-4},\ +2\!\times\!10^{-4},\ 0\right)^\transpose.
\end{aligned}
\end{equation}
Their radii are $R_\up{l} = R_\up{r} = 5\times10^{-4}$, and the interface thickness is $\epsilon = 10^{-5}$.
The initial volume fractions are $\alpha_1 = \alpha_\up{min} = 10^{-4}$ for the gas phase
and $\alpha_1 = \alpha_\up{max} = 1 - 10^{-4}$ for the liquid phase.
The atmospheric pressure is $p_\up{atm} = 10^5$ and 
gravity effects are not present ($\vec{g} = \vec{0}$).
The droplets are set on an off-center collision path by superimposing, to each droplet,
a diffuse circular region 
in which the velocity field is $\vec{u} = (1,\ 0,\ 0)^\transpose$
for the left droplet
and $\vec{u} = (-1,\ 0,\ 0)^\transpose$ for the right droplet, 
while the surrounding fluid is at rest. These circular regions
are defined by the same smooth profile used for the droplets but with radius larger by 
a factor $k=1.1$ with respect the droplet itself, and with interface thickness $\epsilon=5\times10^{-5}$.
The strain relaxation times are $\tau_\up{w} = 9.3750\times10^{-8}$ and $\tau_\up{a} = 1.4064\times10^{-6}$ for the liquid and the gas respectively, 
which translates to kinematic viscosities $\nu_\up{w} = 10^{-6}$ and $\nu_\up{a} = 1.5\times10^{-5}$.
The densities are $\rho_1 = 10^{3}$ for the liquid phase and  $\rho_2 = 1$ for the gas phase.
The remaining material parameters are $\gamma_1 = 8.0$, $\gamma_2 = 1.4$, $\Pi_1 = 10^6$, $\Pi_2 = 0$ for the stiffened gas equation of state
and $\cshear=8$ for the mesoscale strain energy closure.
Three separate numerical experiments are carried out with different values of the Weber number $\weber$, defined
by changing the surface tension coefficient from $\sigma = 0.2\times10^{-3}$ (corresponding to $\weber\simeq 80$), to $\sigma = 0.1\times10^{-3}$ (i.e. $\weber\simeq 160$), to 
$\sigma = 0.05\times10^{-3}$ ($\weber\simeq  320$).
The computational domain is $\Omega = [-4\times10^{-3},\ 4\,\times10^{-3}]\times[-3\,\times10^{-3},\ 3\,\times10^{-3}]$ and 
the mesh is of Cartesian type with $2048$ and $1536$ cells in the $x$ and $y$ directions respectively.

The qualitative behaviour is in agreement with experimental findings \cite{dropletcollision1, dropletcollision2, dropletcollision3}
about collision regimes in liquid droplets.
In particular in Figure~\ref{fig:collision2} one can clearly distinguish different separation modes 
taking place following the collision: at low Weber number ($\weber \simeq 80$) the collision results in stretching and separation of the two droplets. At
all values of Weber number, the effects of Rayleigh--Plateau instability are evident and in particular for the higher Weber numbers 
they result in the formation of multiple small secondary droplets.
In the early stages of the simulations (Figure~\ref{fig:collision1}) one can also clearly see how 
sharp interfacial features are penalized by surface tension, which tends to reduce the curvature of interfaces, 
the more the stronger the capillarity forces with respect to convective effects.
The same test, with $\weber \simeq 80$ and a smaller domain $\Omega = [-3\times10^{-3},\ 3\,\times10^{-3}]\times[-1.5\,\times10^{-3},\ 1.5\,\times10^{-3}]$ 
is employed for a study of the strong scaling performance of the computational code. In this case the grid counts $8192$ by $4096$ cells
and the simulations have been carried out with the aid of the HPE--Hawk supercomputer at the HLRS in Stuttgart, Germany, 
in order to test the scaling capabilities of our semi-implicit computational code on massively parallel distributed memory supercomputer architectures.
The results of this latter test are summarised in Figure~\ref{fig:scaling}: the strong scaling tests starts 
from 
64 CPU cores (half node) and extends up to 65\,536 CPU cores (512 nodes) of the
HPE--Hawk supercomputer at HLRS in Stuttgart. Our semi-implicit scheme 
achieves excellent scaling performance 
(more than $\sim95\%$)
up to about 16k CPU cores, also thanks to cache effects mitigating the main bottleneck of the
scheme which is the solution of the pressure system via matrix-free conjugate gradient method and only on 
65k CPU cores, we see the speedup efficiency starting to drop significantly to about $50\%$.

\subsection{Multiphase Rayleigh--Taylor instability}
\begin{figure}[!bp]
   \includegraphics[width=\textwidth]{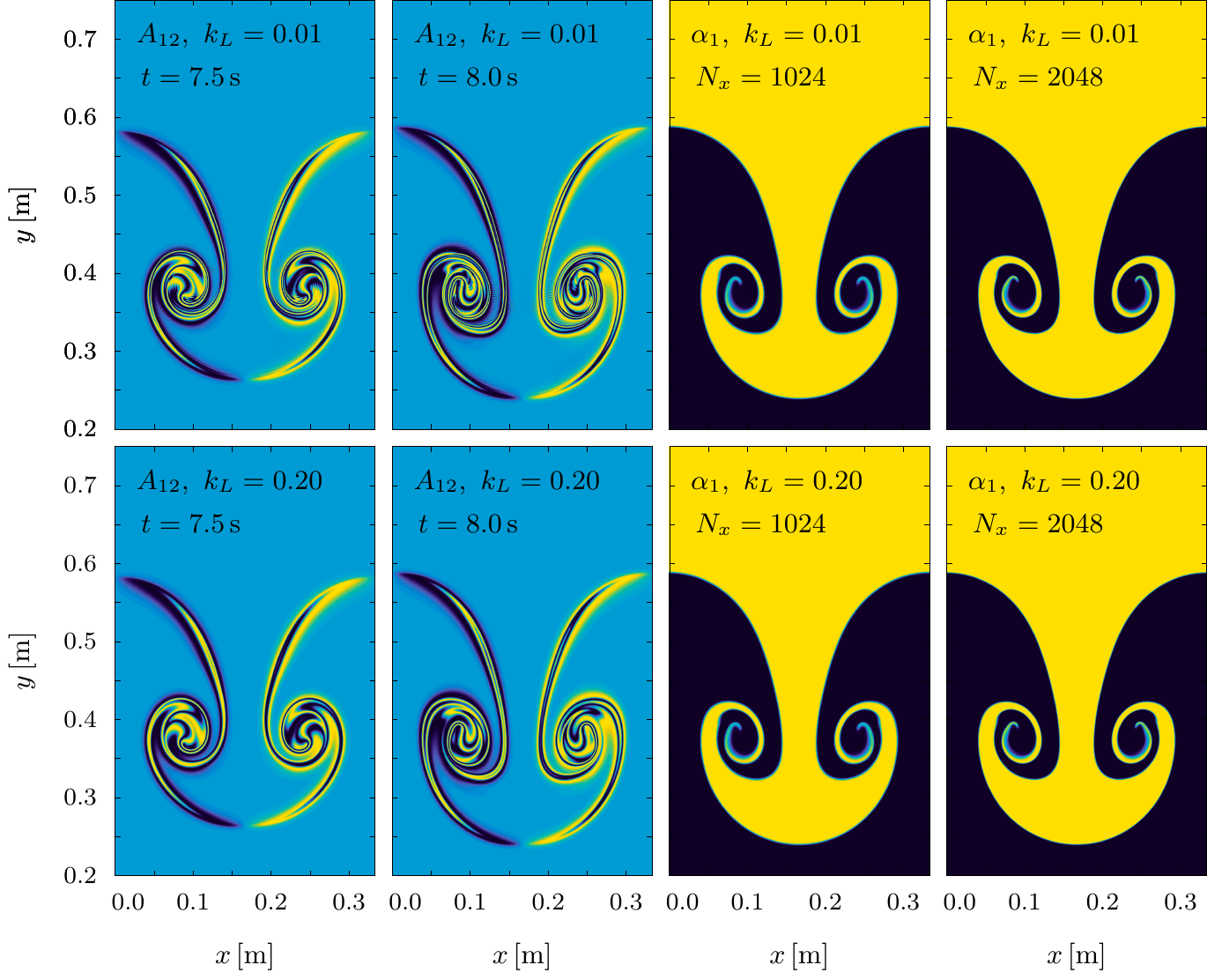}
   \caption{Mesh convergence test for the Rayleigh--Taylor instability problem: the rightmost panels show that mesh convergence has been achieved with a uniform 
   Cartesian mesh counting 1024 by 3072 cells, since no new flow features appear if the mesh resolution is doubled. In the left panels we show the effects 
   of choosing a different reduction coefficient $k_L$ for the compatible vector Laplacian diffusion operator applied to the distortion matrix $\vec{A}$.
   By comparing the top and bottom rows one can see that increasing $k_L$ has a visible effect on the finer features of the distortion field, but this
   does not translate to comparably visible effects in the volume fraction contours on the right.}   
   \label{fig:rt3}
\end{figure}
\begin{figure}[!bp]
   \includegraphics[width=\textwidth]{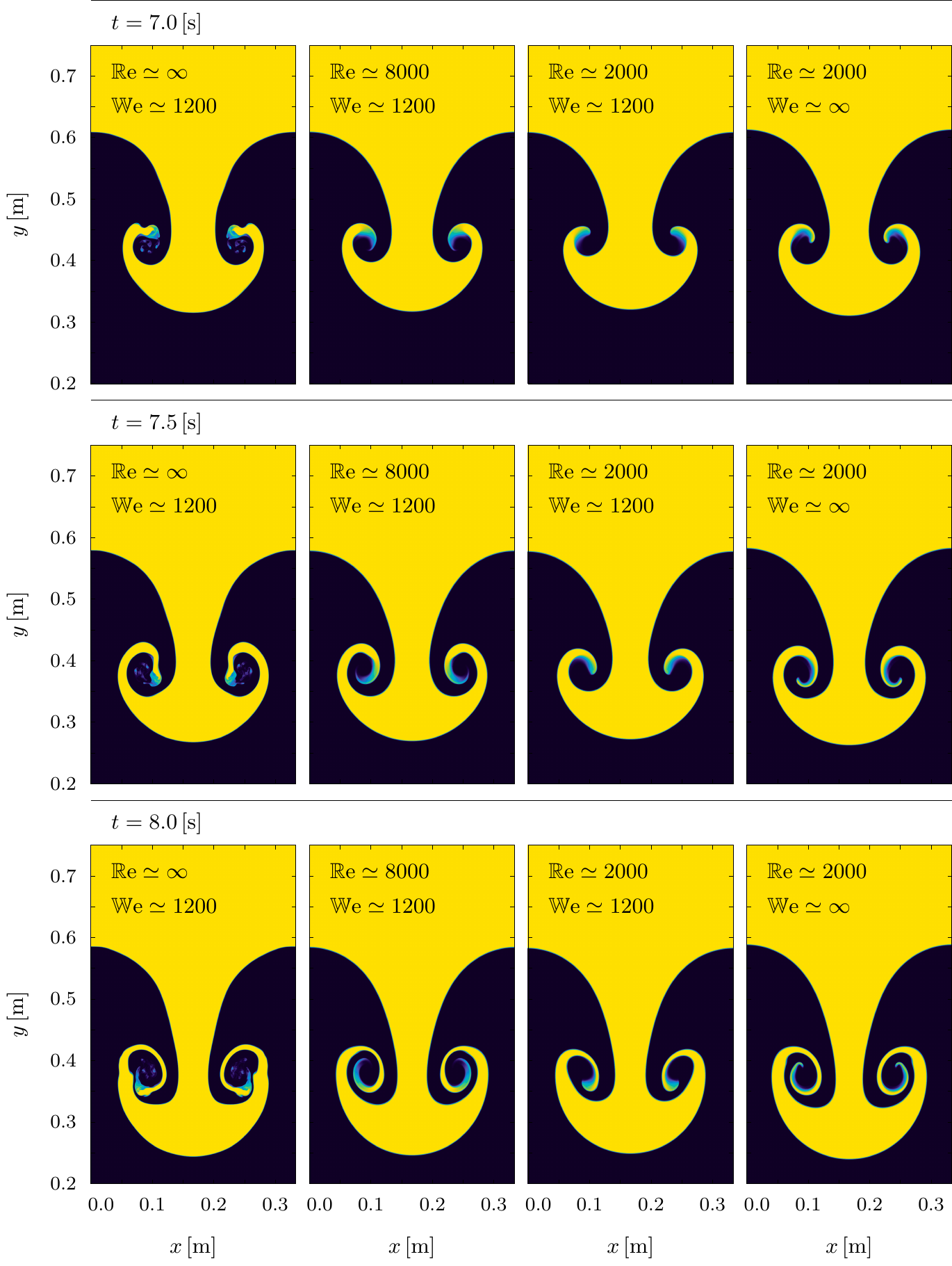}
   \caption{Volume fraction contours for the Rayleigh--Taylor instability problem with varying viscosity and surface tension at different times.}   
   \label{fig:rt2}
\end{figure}
\begin{figure}[!bp]
   \includegraphics[width=\textwidth]{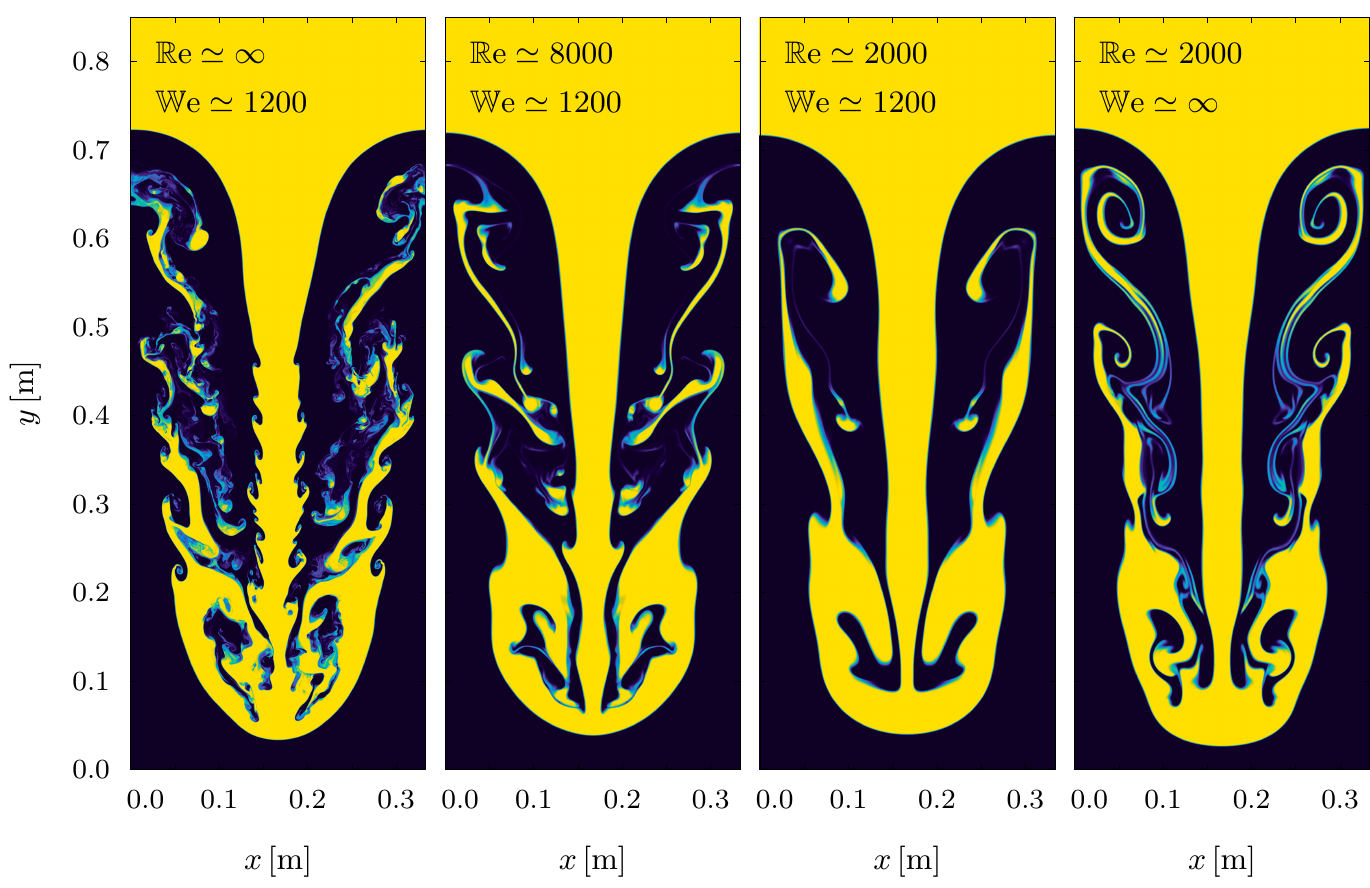}
   \caption{Development of the Rayleigh--Taylor instability at time $t=12.5\,\up{s}$, with varying viscosity and surface tension.
   From left to right: the first panel shows the solution with Weber number $\weber \simeq 1200$ in 
   inviscid flow; the second panel shows how the solution changes if mild viscosity is introduced ($\reynolds\simeq8000$);
   the third panels is relative to a simulation with increased viscosity ($\reynolds\simeq2000$); the 
   fourth and last panel depicts the results of a simulation of viscous flow without surface tension.
   The curl-preserving semi-implicit scheme with grid size 4096 by $12\,288$ has been employed for all simulations.}   
   \label{fig:rt4}
\end{figure}

\begin{figure}[!bp]
   \includegraphics[width=1.0\textwidth]{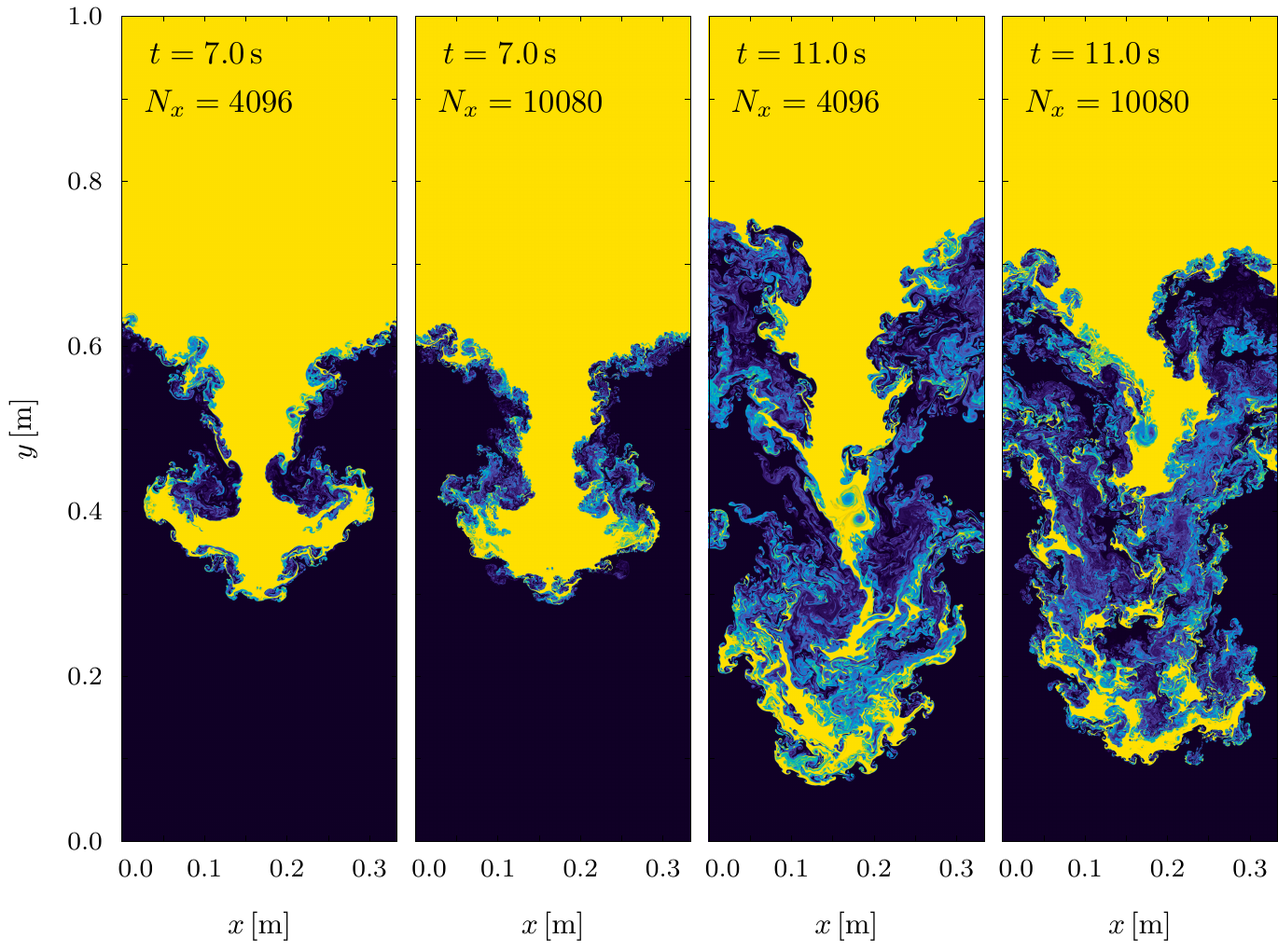}\\[3mm]
   \includegraphics[width=1.0\textwidth]{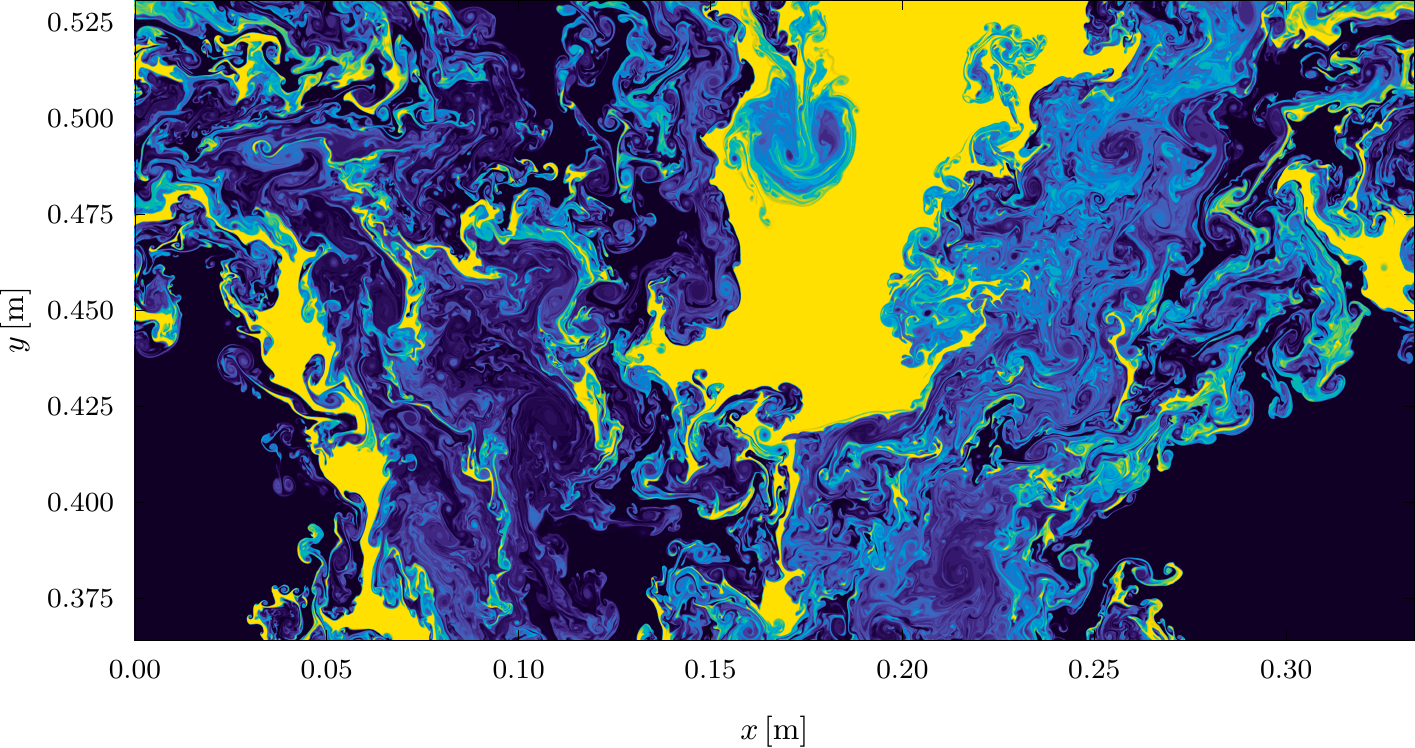}
   \caption{Volume fraction plots for the inviscid two-phase Rayleigh--Taylor instability 
   solved with the Structure Preserving Semi-Implicit Finite Volume scheme applied to the unified model of continuum mechanics
   in the inviscid limit $(\tau=10^{-14})$, at different mesh resolutions ($4096\!\times\!12\,288$ and $10\,080\!\times\!30\,240$).}   
   \label{fig:rt1}
\end{figure}
Finally, we put all elements of the scheme to the test at the same time, by simulating a low Mach, genuinely two-phase, Rayleigh--Taylor instability
with viscosity and surface tension. The setup follows \cite{ReALE2010}, but with the notable modification that, 
instead of initialising a single fluid with two different densities (one above a horizontal material interface, one below), 
we define two separate density fields, each one with \textit{constant} phase densities $\rho_1 = \rho_\up{t}$ and $\rho_2 = \rho_\up{b}$,
then 
distinguishing the two fluids by means of a jump from $\alpha_1 = \alpha_\up{min} = 10^{-4}$
to $\alpha_1 = \alpha_\up{max} = 1 - 10^{-4}$ in the volume fraction field.
This renders the problem much more challenging because near-vacuum states of one of the two phases are introduced almost throughout the computational domain.
The curved material interface location is 
\begin{equation}
   y_I = 0.5 + 0.01\,\cos\left(6\,\pi\,x\right)
\end{equation}
and we impose the transition between the two states by means of a smooth switch function 
\begin{equation}
   s = \frac{1}{2} + \frac{1}{2}\,\up{erf}\left(\frac{y - y_I}{\delta}\right),
\end{equation}
with interface thickness $\delta = \max(0.004,\ 6\,\Delta x)$. This is intended to suppress 
spurious instabilities that would be triggered by inaccurate representation of the initial 
condition on a discrete Cartesian grid (stairstepping), thus allowing only the physical instabilities 
to develop.
The volume fraction field is then 
$\alpha_1 = s\,\alpha_\up{min} + (1 - s)\,\alpha_\up{max}$, while
the top and bottom pressures are 
\begin{equation}
\begin{aligned}
   &p_\up{t} = 1 + \rho_\up{t}\,\vec{g}\cdot\uvec{e}_y\,\left(1 - y\right),\\[1mm]
   &p_\up{b} = 1 + 0.6\,\rho_\up{t}\,\vec{g}\cdot\uvec{e}_y + \rho_\up{b}\,\left(1/2 - y\right),
\end{aligned}
\end{equation}
and, just like the volume fraction $\alpha_1$, the pressure field is given by 
$p = s\,p_\up{t} + (1 - s)\,p_\up{b}$.
We initially set $\vec{A} = \vec{I}$, $\vec{u} = \vec{0}$, and the remaining parameters common to all simulations carried out with 
regard to this test problem are 
$\rho_\up{t} = 2$,
$\rho_\up{b} = 1$,
the gravity vector is 
$\vec{g} = (0,\ -0.1,\ 0)^\transpose$, and finally 
$\vec{b}$ is initialised as a compatible discrete gradient of an auxiliary colour function field $s_c$ 
given by $s_c = 1/2 + \up{erf}\left[(y - y_I)/(2\,\delta)\right]/2$.
With this setup, we carry out a parametric study of the behaviour of the instability.

First, in order to verify mesh convergence of the solution algorithm, 
we set the Reynolds number $\reynolds \simeq 2000$ (which translates to $\cshear=0.3$ and $\tau = 2\times 10^{-3}$) 
and $\weber \simeq \infty$ (i.e. we neglect surface tension), 
and we carry out two simulation on two different meshes, one composed of $1024 \times 3072$ elements, and one of $2048 \times 6144$
elements, and show, in Figure~\ref{fig:rt3}, that indeed the structure of the solution does not depend on 
mesh effects. Furthermore, again in Figure~\ref{fig:rt3}, we also show that the method is robust with respect 
to the choice of the scaling factor for the compatible numerical viscosity $k_L$: by varying such
a factor from $k_L = 0.01$ to $k_L = 0.2$ we see that some of the fine scale structures in the distortion 
field $\vec{A}$ are lost to numerical dissipation. However, this does not translate into a visible effect in the 
shape of the interface separating the two fluids.
This can be explained by the fact that the rotational component $\vec{R}$ of the distortion field $\vec{A}$
is significantly affected by numerical viscosity, but instead this is not the case for the stress $\vec{\sigma}_\up{s} = -\rho\,\cshear^2\,\vec{G}\,\dev{\vec{G}}$, 
which does \textit{not} depend on these rotations, thus leading to the same global dynamics for 
both choices of the numerical viscosity coefficient $k_L$.

In Figure~\ref{fig:rt2} we report a variety of snapshots for this test problem at different times between $t = 7\,\up{s}$ and $t = 8\,\up{s}$, 
at different values of the Weber number $\weber = \rho\,U^2\,L/\sigma$
and at different values of the Reynolds number $\reynolds = U\,L/\nu$.
These runs employ a uniform Cartesian grid counting $4096$ by $12\,288$ elements.
With $\weber \simeq 1200$ we indicate simulations carried out with $\sigma = 2 \times 10^{-5}$, and $\weber \simeq \infty$ corresponds to $\sigma = 0$.
With $\reynolds \simeq 2000$ we indicate simulations carried out with $\cshear=0.3$ and $\tau = 2\times 10^{-3}$, and the label
$\reynolds \simeq 8000$ corresponds to $\cshear=0.3$ and $\tau = 5\times 10^{-4}$, while the inviscid limit 
$\reynolds \simeq \infty$ is given by $\cshear=0.3$ and $\tau = 10^{-14}$.
The results show the stabilising effects of surface tension and viscosity and the characteristic morphology they
regulate. See \cite{experimentrm} for striking experimental results on the Richtmyer--Meshkov instability, 
featuring similar flow structures.

In Figure~\ref{fig:rt4}, we plot the results of the same battery of simulations at a later time $t = 12.5\,\up{s}$. 
Again, the distinctive features determined by the presence of viscosity and surface tension, 
as well as Rayleigh--Plateau secondary instabilities, 
can be clearly identified.

In a final set of runs we study the $\reynolds\to\infty$, $\weber\to\infty$ limit of the governing equations,
by setting $\tau = 10^{-14}$, $\sigma = 0$, $\cshear = 0.1$. With this test we intend to verify
the applicability of our MPI-parallel computational code to large scale simulations, its resolution properties, 
and its robustness in turbulent multiphase flows.
We report, in Figure~\ref{fig:rt1}, the results of a large scale simulation of turbulent multiphase flow, 
run on $32\up{k}$ CPU cores of the HPE--Hawk supercomputer at the HLRS in Stuttgart, with a grid 
counting $10\,080$ by $30\,240$ uniform Cartesian elements. For comparison, in the same Figure are
also plotted the results of the same setup on a coarser grid of $4096$ by $12\,288$ cells.

%% file: SISurfaceTension.bbl
\begin{thebibliography}{10}
\providecommand{\url}[1]{{#1}}
\providecommand{\urlprefix}{URL }
\expandafter\ifx\csname urlstyle\endcsname\relax
  \providecommand{\doi}[1]{DOI~\discretionary{}{}{}#1}\else
  \providecommand{\doi}{DOI~\discretionary{}{}{}\begingroup
  \urlstyle{rm}\Url}\fi

\bibitem{abgrallcondition}
Abgrall, R.: How to prevent pressure oscillations in multicomponent flow
  calculations: A quasi conservative approach.
\newblock Journal of Computational Physics \textbf{125}(1), 150--160 (1996)

\bibitem{dropletcollision2}
Al-Dirawi, K.H., Bayly, A.E.: An experimental study of binary collisions of
  miscible droplets with non-identical viscosities.
\newblock Experiments in Fluids \textbf{61} (2020)

\bibitem{imex1}
Ascher, U.M., Ruuth, S.J., Spiteri, R.J.: Implicit-explicit runge-kutta methods
  for time-dependent partial differential equations.
\newblock Applied Numerical Mathematics \textbf{25}(2), 151--167 (1997).
\newblock Special Issue on Time Integration

\bibitem{BaerNunziato1986}
Baer, M.R., Nunziato, J.W.: A two--phase mixture theory for the
  deflagration-to-detonation transition {(DDT)} in reactive granular materials.
\newblock J. Multiphase Flow \textbf{12}, 861--889 (1986)

\bibitem{Barton2019}
Barton, P.: {An interface-capturing Godunov method for the simulation of
  compressible solid-fluid problems}.
\newblock Journal of Computational Physics \textbf{390}, 25--50 (2019)

\bibitem{Bell1989}
Bell, J.B., Coletta, P., Glaz, H.M.: A second-order projection method for the
  incompressible {Navier-Stokes} equations.
\newblock Journal of Computational Physics \textbf{85}, 257--283 (1989)

\bibitem{Berry2008a}
Berry, R., Saurel, R., Petitpas, F., Daniel, E., Le~M{\'{e}}tayer, O.,
  Gavrilyuk, S.: {Progress in the Development of Compressible , Multiphase Flow
  Modeling Capability for Nuclear Reactor Flow Applications}.
\newblock Idaho National Laboratory (October) (2008)

\bibitem{Besseling1968}
Besseling, J.: A thermodynamic approach to rheology.
\newblock In: H.~Parkus, L.~Sedov (eds.) Irreversible Aspects of Continuum
  Mechanics and Transfer of Physical Characteristics in Moving Fluids, IUTAM
  Symposia, pp. 16--53. Springer Vienna (1968)

\bibitem{imex4}
Boscarino, S.: Error analysis of imex runge–kutta methods derived from
  differential-algebraic systems.
\newblock SIAM Journal on Numerical Analysis \textbf{45}(4), 1600--1621 (2007)

\bibitem{imex6}
Boscarino, S.: On an accurate third order implicit-explicit runge–kutta
  method for stiff problems.
\newblock Applied Numerical Mathematics \textbf{59}(7), 1515--1528 (2009)

\bibitem{imex7}
Boscarino, S., Pareschi, L., Russo, G.: Implicit-explicit runge--kutta schemes
  for hyperbolic systems and kinetic equations in the diffusion limit.
\newblock SIAM Journal on Scientific Computing \textbf{35}(1), A22--A51 (2013)

\bibitem{imex5}
Boscarino, S., Russo, G.: On a class of uniformly accurate imex runge–kutta
  schemes and applications to hyperbolic systems with relaxation.
\newblock SIAM Journal on Scientific Computing \textbf{31}(3), 1926--1945
  (2009)

\bibitem{boscheriulaggpr}
Boscheri, W., Chiocchetti, S., Peshkov, I.: A cell-centered implicit-explicit
  lagrangian scheme for a unified model of nonlinear continuum mechanics on
  unstructured meshes.
\newblock Journal of Computational Physics \textbf{451}, 110852 (2022)

\bibitem{sigpr}
Boscheri, W., Dumbser, M., Ioriatti, M., Peshkov, I., Romenski, E.: A
  structure-preserving staggered semi-implicit finite volume scheme for
  continuum mechanics.
\newblock Journal of Computational Physics \textbf{424}, 109866 (2021)

\bibitem{boscherisigpr}
Boscheri, W., Dumbser, M., Ioriatti, M., Peshkov, I., Romenski, E.: A
  structure-preserving staggered semi-implicit finite volume scheme for
  continuum mechanics.
\newblock Journal of Computational Physics \textbf{424}, 109866 (2021)

\bibitem{Iollo2017}
de~Brauer, A., Iollo, A., Milcent, T.: {A Cartesian Scheme for Compressible
  Multimaterial Hyperelastic Models with Plasticity}.
\newblock Communications in Computational Physics \textbf{22}, 1362--1384
  (2017)

\bibitem{sarayhtc}
Busto, S., Dumbser, M., Peshkov, I., Romenski, E.: On thermodynamically
  compatible finite volume schemes for continuum mechanics.
\newblock SIAM Journal on Scientific Computing \textbf{44}(3), A1723--A1751
  (2022)

\bibitem{castro2006}
Castro, M.J., Gallardo, J.M., Par\'es, C.: High-order finite volume schemes
  based on reconstruction of states for solving hyperbolic systems with
  nonconservative products. applications to shallow-water systems.
\newblock Mathematics of Computations \textbf{75}, 1103--1134 (2006)

\bibitem{CasulliCompressible}
Casulli, V., Greenspan, D.: Pressure method for the numerical solution of
  transient, compressible fluid flows.
\newblock International Journal for Numerical Methods in Fluids \textbf{4}(11),
  1001--1012 (1984)

\bibitem{CasulliZanolli2010}
Casulli, V., Zanolli, P.: {A nested Newton--type algorithm for finite volume
  methods solving Richards' equation in mixed form}.
\newblock SIAM Journal on Scientific Computing \textbf{32}, 2255--2273 (2009)

\bibitem{CasulliZanolli2012}
Casulli, V., Zanolli, P.: Iterative solutions of mildly nonlinear systems.
\newblock Journal of Computational and Applied Mathematics \textbf{236},
  3937--3947 (2012)

\bibitem{chiocchettithesis}
Chiocchetti, S.: High order numerical methods for a unified theory of fluid and
  solid mechanics (2022)

\bibitem{chiocchettimueller}
Chiocchetti, S., M{\"{u}}ller, C.: {A Solver for Stiff Finite-Rate Relaxation
  in Baer-Nunziato Two-Phase Flow Models}.
\newblock Fluid Mechanics and its Applications \textbf{121}, 31--44 (2020)

\bibitem{hstglm}
Chiocchetti, S., Peshkov, I., Gavrilyuk, S., Dumbser, M.: High order ader
  schemes and glm curl cleaning for a first order hyperbolic formulation of
  compressible flow with surface tension.
\newblock Journal of Computational Physics \textbf{426}, 109898 (2021)

\bibitem{Dedneretal}
Dedner, A., Kemm, F., Kr\"oner, D., Munz, C.D., Schnitzer, T., Wesenberg, M.:
  Hyperbolic divergence cleaning for the {MHD} equations.
\newblock Journal of Computational Physics \textbf{175}, 645--673 (2002)

\bibitem{SIMHD}
Dumbser, M., Balsara, D., Tavelli, M., Fambri, F.: {A divergence-free
  semi-implicit finite volume scheme for ideal, viscous and resistive
  magnetohydrodynamics}.
\newblock International Journal for Numerical Methods in Fluids \textbf{89},
  16--42 (2019)

\bibitem{DumbserCasulli2016}
Dumbser, M., Casulli, V.: {A conservative, weakly nonlinear semi-implicit
  finite volume method for the compressible Navier-Stokes equations with
  general equation of state}.
\newblock Applied Mathematics and Computation \textbf{272}, 479--497 (2016)

\bibitem{godbook}
Dumbser, M., Chiocchetti, S., Peshkov, I.: On Numerical Methods for Hyperbolic
  PDE with Curl Involutions, pp. 125--134.
\newblock Springer International Publishing, Cham (2020)

\bibitem{GPRmodel}
Dumbser, M., Peshkov, I., Romenski, E., Zanotti, O.: {High order ADER schemes
  for a unified first order hyperbolic formulation of continuum mechanics:
  Viscous heat-conducting fluids and elastic solids}.
\newblock Journal of Computational Physics \textbf{314}, 824--862 (2016)

\bibitem{GPRmodelMHD}
Dumbser, M., Peshkov, I., Romenski, E., Zanotti, O.: {H}igh order {ADER}
  schemes for a unified first order hyperbolic formulation of {N}ewtonian
  continuum mechanics coupled with electro-dynamics.
\newblock Journal of Computational Physics \textbf{348}, 298--342 (2017)

\bibitem{browneinstein}
Einstein, A.: Über die von der molekularkinetischen theorie der wärme
  geforderte bewegung von in ruhenden flüssigkeiten suspendierten teilchen.
\newblock Annalen der Physik \textbf{322}(8), 549--560 (1905)

\bibitem{FavrGavr2012}
Favrie, N., Gavrilyuk, S.: Diffuse interface model for compressible
  fluid-compressible elastic-plastic solid interaction.
\newblock Journal of Computational Physics \textbf{231}, 2695--2723 (2012)

\bibitem{FavrieGavrilyukSaurel}
Favrie, N., Gavrilyuk, S., Saurel, R.: {Solid--fluid diffuse interface model in
  cases of extreme deformations}.
\newblock Journal of Computational Physics \textbf{228}, 6037--6077 (2009)

\bibitem{favrie2009solid}
Favrie, N., Gavrilyuk, S.L., Saurel, R.: Solid--fluid diffuse interface model
  in cases of extreme deformations.
\newblock Journal of computational physics \textbf{228}(16), 6037--6077 (2009)

\bibitem{feynmanbook}
Feynman, R.P.: The Feynman lectures on physics.
\newblock Addison-Wesley Pub. Co., Reading, Mass. (1963)

\bibitem{dropletcollision3}
Finotello, G., Kooiman, R.F., Padding, J.T., Buist, K.A., Jongsma, A., Innings,
  F., Kuipers, J.A.M.: The dynamics of milk droplet–droplet collisions.
\newblock Experiments in Fluids \textbf{59} (2017)

\bibitem{ptrsa}
Gabriel, A.A., Li, D., Chiocchetti, S., Tavelli, M., Peshkov, I., Romenski, E.,
  Dumbser, M.: A unified first-order hyperbolic model for nonlinear dynamic
  rupture processes in diffuse fracture zones.
\newblock Philosophical Transactions of the Royal Society A \textbf{379} (2021)

\bibitem{Gavrilyuk2008}
Gavrilyuk, S., Favrie, N., Saurel, R.: {Modelling wave dynamics of compressible
  elastic materials}.
\newblock Journal of Computational Physics \textbf{227}(5), 2941--2969 (2008)

\bibitem{God1972MHD}
Godunov, S.: Symmetric form of the magnetohydrodynamic equation.
\newblock Numerical Methods for Mechanics of Continuum Medium \textbf{3}(1),
  26--34 (1972)

\bibitem{God1978}
Godunov, S.K.: Elements of mechanics of continuous media.
\newblock Nauka (1978)

\bibitem{Godunov1996}
Godunov, S.K., Mikha{\^{i}}lova, T.Y., Romenski{\^{i}}, E.I.: {Systems of
  thermodynamically coordinated laws of conservation invariant under
  rotations}.
\newblock Siberian Mathematical Journal \textbf{37}(4), 690--705 (1996)

\bibitem{Pesh2010}
Godunov, S.K., Peshkov, I.: Thermodynamically consistent nonlinear model of
  elastoplastic maxwell medium.
\newblock Computational Mathematics and Mathematical Physics \textbf{50(8)},
  1409--1426 (2010)

\bibitem{GodRom1998}
Godunov, S.K., Romenski, E.: Elements of mechanics of continuous media.
\newblock Nauchnaya Kniga (1998.)

\bibitem{GodRom2003}
Godunov, S.K., Romenski, E.: Elements of Continuum Mechanics and Conservation
  Laws.
\newblock Kluwer Academic/Plenum Publishers (2003)

\bibitem{GodunovRomenski72}
Godunov, S.K., Romenski, E.I.: Nonstationary equations of the nonlinear theory
  of elasticity in {Euler} coordinates.
\newblock Journal of Applied Mechanics and Technical Physics \textbf{13},
  868--885 (1972)

\bibitem{hankplasticity}
Hank, S., Gavrilyuk, S., Favrie, N., Massoni, J.: Impact simulation by an
  eulerian model for interaction of multiple elastic-plastic solids and fluids.
\newblock International Journal of Impact Engineering \textbf{109}, 104--111
  (2017)

\bibitem{brenn1}
Hinterbichler, H., Planchette, C., Brenn, G.: Ternary drop collisions.
\newblock Experiments in Fluids \textbf{56}, 190/1--190/12 (2015)

\bibitem{Jackson2019a}
Jackson, H., Nikiforakis, N.: {A numerical scheme for non-Newtonian fluids and
  plastic solids under the GPR model}.
\newblock Journal of Computational Physics \textbf{387}, 410--429 (2019)

\bibitem{Jackson2019}
Jackson, H., Nikiforakis, N.: {A unified Eulerian framework for multimaterial
  continuum mechanics}.
\newblock Journal of Computational Physics \textbf{401}, 109022 (2019)

\bibitem{kapila2001}
Kapila, A.K., Menikoff, R., Bdzil, J.B., Son, S.F., Stewart, D.S.: Two-phase
  modelling of deflagration-to-detonation in granular materials: reduced
  equations.
\newblock Physics of Fluids \textbf{13}, 3002--3024 (2001)

\bibitem{imex2}
Kennedy, C.A., Carpenter, M.H.: Additive runge–kutta schemes for
  convection–diffusion–reaction equations.
\newblock Applied Numerical Mathematics \textbf{44}(1), 139--181 (2003)

\bibitem{KlaMaj}
Klainermann, S., Majda, A.: Singular limits of quasilinear hyperbolic systems
  with large parameters and the incompressible limit of compressible fluid.
\newblock Communications on Pure and Applied Mathematics \textbf{34}, 481--524
  (1981)

\bibitem{KlaMaj82}
Klainermann, S., Majda, A.: Compressible and incompressible fluids.
\newblock Communications on Pure and Applied Mathematics \textbf{35}, 629--651
  (1982)

\bibitem{Klein2001}
Klein, R., Botta, N., Schneider, T., Munz, C., S.Roller, Meister, A., Hoffmann,
  L., Sonar, T.: Asymptotic adaptive methods for multi-scale problems in fluid
  mechanics.
\newblock Journal of Engineering Mathematics \textbf{39}, 261--343 (2001)

\bibitem{ReALE2010}
Loub{\`e}re, R., Maire, P.H., Shashkov, M., Breil, J., Galera, S.: {ReALE: A
  reconnection-based arbitrary-Lagrangian–Eulerian method}.
\newblock Journal of Computational Physics \textbf{229}, 4724--4761 (2010)

\bibitem{MunzDumbserRoller}
Munz, C., Dumbser, M., Roller, S.: Linearized acoustic perturbation equations
  for low {Mach} number flow with variable density and temperature.
\newblock Journal of Computational Physics \textbf{224}, 352--364 (2007)

\bibitem{MunzCleaning}
Munz, C., Omnes, P., Schneider, R., Sonnendr\"ucker, E., Voss, U.: {Divergence
  Correction Techniques for Maxwell Solvers Based on a Hyperbolic Model}.
\newblock Journal of Computational Physics \textbf{161}, 484--511 (2000)

\bibitem{Munz2003}
Munz, C., Roller, S., Klein, R., Geratz, K.: The extension of incompressible
  flow solvers to the weakly compressible regime.
\newblock Computers and Fluids \textbf{32}, 173--196 (2003)

\bibitem{Ndanou2014}
Ndanou, S., Favrie, N., Gavrilyuk, S.: Criterion of hyperbolicity in
  hyperelasticity in the case of the stored energy in separable form.
\newblock Journal of Elasticity \textbf{115}, 1--25 (2014)

\bibitem{NdanouFavrieGavrilyuk}
Ndanou, S., Favrie, N., Gavrilyuk, S.: {Multi--solid and multi--fluid diffuse
  interface model: Applications to dynamic fracture and fragmentation}.
\newblock Journal of Computational Physics \textbf{295}, 523--555 (2015)

\bibitem{experimentrm}
Niederhaus, C.E., Jacobs, J.W.: Experimental study of the richtmyer–meshkov
  instability of incompressible fluids.
\newblock Journal of Fluid Mechanics \textbf{485}, 243–277 (2003)

\bibitem{pares2006}
Par\'es, C.: Numerical methods for nonconservative hyperbolic systems: a
  theoretical framework.
\newblock SIAM Journal on Numerical Analysis \textbf{44}, 300--321 (2006)

\bibitem{imex3}
Pareschi, L., Russo, G.: Implicit-explicit runge-kutta schemes and applications
  to hyperbolic systems with relaxation.
\newblock Journal of Scientific Computing \textbf{25}(1), 129--155 (2005)

\bibitem{MunzPark}
Park, J., Munz, C.: Multiple pressure variables methods for fluid flow at all
  mach numbers.
\newblock International Journal for Numerical Methods in Fluids \textbf{49},
  905--931 (2005)

\bibitem{ilyanonnewtonian}
Peshkov, I., Dumbser, M., Boscheri, W., Romenski, E., Chiocchetti, S.,
  Ioriatti, M.: Simulation of non-newtonian viscoplastic flows with a unified
  first order hyperbolic model and a structure-preserving semi-implicit scheme.
\newblock Computers {\&} Fluids \textbf{224}, 104963 (2021)

\bibitem{SHTC-GENERIC-CMAT}
Peshkov, I., Pavelka, M., Romenski, E., Grmela, M.: {Continuum mechanics and
  thermodynamics in the Hamilton and the Godunov-type formulations}.
\newblock Continuum Mechanics and Thermodynamics \textbf{30}(6), 1343--1378
  (2018)

\bibitem{PeshRom2014}
Peshkov, I., Romenski, E.: A hyperbolic model for viscous {{N}ewtonian} flows.
\newblock Continuum Mechanics and Thermodynamics \textbf{28}, 85--104 (2016)

\bibitem{brenn2}
Planchette, C., Hinterbichler, H., Liu, M., Bothe, D., Brenn, G.: Colliding
  drops as coalescing and fragmenting liquid springs.
\newblock Journal of Fluid Mechanics \textbf{814}, 277--300 (2017)

\bibitem{Powell1997}
Powell, K.: {An Approximate Riemann Solver for Magnetohydrodynamics}.
\newblock In: v.L.B. M.Y., V.R. J. (eds.) Upwind and High-Resolution Schemes,
  pp. 570--583. Springer Berlin Heidelberg, Berlin, Heidelberg (1997)

\bibitem{powell1999}
Powell, K., Roe, P., Linde, T., Gombosi, T., Zeeuw, D.D.: A solution-adaptive
  upwind scheme for ideal magnetohydrodynamics.
\newblock Journal of Computational Physics \textbf{154}(2), 284--309 (1999)

\bibitem{Rom2001}
Romensky, E.I.: {Thermodynamics and hyperbolic systems of balance laws in
  continuum mechanics}.
\newblock In: E.~Toro (ed.) Godunov Methods: Theory and Applications, pp.
  745--761. Springer US, New York (2001)

\bibitem{Schmidmayer2017}
Schmidmayer, K., Petitpas, F., Daniel, E., Favrie, N., Gavrilyuk, S.: A model
  and numerical method for compressible flows with capillary effects.
\newblock Journal of Computational Physics \textbf{334}, 468--496 (2017)

\bibitem{dropletcollision1}
Sommerfeld, M., Pasternak, L.: Advances in modelling of binary droplet
  collision outcomes in sprays: A review of available knowledge.
\newblock International Journal of Multiphase Flow \textbf{117}, 182--205
  (2019)

\bibitem{tavellicrack}
Tavelli, M., Chiocchetti, S., Romenski, E., Gabriel, A.A., Dumbser, M.:
  Space-time adaptive ader discontinuous galerkin schemes for nonlinear
  hyperelasticity with material failure.
\newblock Journal of Computational Physics \textbf{422}, 109758 (2020)

\bibitem{Tavelli2015}
Tavelli, M., Dumbser, M.: {A staggered space--time discontinuous Galerkin
  method for the incompressible Navier--Stokes equations on two--dimensional
  triangular meshes}.
\newblock Computers and Fluids \textbf{119}, 235--249 (2015)

\bibitem{ToroVazquez}
Toro, E., V\'azquez-Cend\'on, M.: {Flux splitting schemes for the Euler
  equations}.
\newblock Computers and Fluids \textbf{70}, 1--12 (2012)

\bibitem{torobook}
Toro, E.F.: Riemann Solvers and Numerical Methods for Fluid Dynamics. A
  Practical Introduction, Third edition.
\newblock Springer-Verlag, Berlin (2009)

\bibitem{vanleer1974}
{van Leer}, B.: Towards the ultimate conservative difference scheme. ii.
  monotonicity and conservation combined in a second-order scheme.
\newblock Journal of Computational Physics \textbf{14}, 361--370 (1974)

\bibitem{vanleer1979}
{van Leer}, B.: Towards the ultimate conservative difference scheme. v. a
  second-order sequel to godunov's method.
\newblock Journal of Computational Physics \textbf{32}, 101--136 (1979)

\bibitem{wood1930}
Wood, A.: A Textbook of Sound.
\newblock B. Bell and Sons LTD, London (1930)

\end{thebibliography}
